\documentclass[12pt]{amsart}

\usepackage{amssymb}
\usepackage{fullpage}
\usepackage{amscd}   
\usepackage[all]{xy} 
\usepackage{comment} 

\DeclareFontEncoding{OT2}{}{} 
\newcommand{\textcyr}[1]{%
 {\fontencoding{OT2}\fontfamily{wncyr}\fontseries{m}\fontshape{n}\selectfont #1}}
\newcommand{\Sha}{{\mbox{\textcyr{Sh}}}}

\usepackage{color}

\newcommand{\bjorn}[1]{}
\newcommand{\ed}[1]{}
\newcommand{\michael}[1]{}

\newcommand{\known}{{\operatorname{known}}}
\newcommand{\Abar}{{\overline{A}}}
\newcommand{\To}{\longrightarrow}
\newcommand{\nichts}{\ensuremath{\left.\right.}}
\newcommand{\ssigma}{{\! \nichts^\sigma \!}}
\newcommand{\GQ}{{{\mathcal G}_\Q}}
\newcommand{\GQp}{{{\mathcal G}_{\Q_p}}}
\newcommand{\Gk}{{{\mathcal G}_k}}
\newcommand{\GK}{{{\mathcal G}_K}}
\newcommand{\YY}{{\mathcal Y}}

\newcommand{\fake}{{\operatorname{fake}}}
\newcommand{\conts}{{\operatorname{conts}}}
\newcommand{\wedgeisom}{{\;\stackrel{\wedge}\isom}\;}
\newcommand{\xbar}{{\bar{x}}}
\newcommand{\ybar}{{\bar{y}}}
\newcommand{\zbar}{{\bar{z}}}
\newcommand{\Xtilde}{{\tilde{X}}}
\newcommand{\pitilde}{{\tilde{\pi}}}
\newcommand{\subs}{{\operatorname{subset}}}

\newcommand{\Aff}{{\mathbb A}}
\newcommand{\C}{{\mathbb C}}
\newcommand{\F}{{\mathbb F}}
\newcommand{\G}{{\mathbb G}}

\newcommand{\PP}{{\mathbb P}}
\newcommand{\Q}{{\mathbb Q}}
\newcommand{\R}{{\mathbb R}}
\newcommand{\Z}{{\mathbb Z}}
\newcommand{\Qbar}{{\overline{\Q}}}

\newcommand{\kbar}{{\overline{k}}}
\newcommand{\Kbar}{{\overline{K}}}

\newcommand{\Fbar}{{\overline{\F}}}

\newcommand{\pp}{{\mathfrak p}}
\newcommand{\frakq}{{\mathfrak q}}

\newcommand{\calC}{{\mathcal C}}

\newcommand{\calE}{{\mathcal E}}

\newcommand{\calJ}{{\mathcal J}}

\newcommand{\calQ}{{\mathcal Q}}

\newcommand{\OO}{{\mathcal O}}

\DeclareMathOperator{\End}{End}
\DeclareMathOperator{\Jac}{Jac}

\DeclareMathOperator{\Hom}{Hom}

\DeclareMathOperator{\Aut}{Aut}
\DeclareMathOperator{\Gal}{Gal}

\DeclareMathOperator{\Br}{Br}

\DeclareMathOperator{\Sel}{Sel}

\DeclareMathOperator{\divv}{div}

\DeclareMathOperator{\Sym}{Sym}
\DeclareMathOperator{\Div}{Div}
\DeclareMathOperator{\Pic}{Pic}

\DeclareMathOperator{\Spec}{Spec}

\DeclareMathOperator{\rank}{rank}

\newcommand{\smooth}{{\operatorname{smooth}}}

\newcommand{\unr}{{\operatorname{unr}}}

\newcommand{\HH}{{\operatorname{H}}}

\newcommand{\M}{\operatorname{M}}
\newcommand{\GL}{\operatorname{GL}}
\newcommand{\SL}{\operatorname{SL}}
\newcommand{\PGL}{\operatorname{PGL}}
\newcommand{\PSL}{\operatorname{PSL}}

\newcommand{\surjects}{\twoheadrightarrow}
\newcommand{\injects}{\hookrightarrow}
\newcommand{\isom}{\simeq}
\newcommand{\notdiv}{\nmid}

\newcommand{\Union}{\bigcup}

\newcommand{\tensor}{\otimes}

\numberwithin{equation}{section}
\newtheorem{theorem}[equation]{Theorem}
\newtheorem{lemma}[equation]{Lemma}
\newtheorem{corollary}[equation]{Corollary}
\newtheorem{proposition}[equation]{Proposition}

\theoremstyle{definition}

\theoremstyle{remark}
\newtheorem{remark}[equation]{Remark}
\newtheorem{remarks}[equation]{Remarks}
\newtheorem{warning}[equation]{Warning}

\renewcommand{\arraystretch}{1.4}

\newif\ifpdf
\ifx\pdfoutput\undefined
  \pdffalse
\else
\ifnum\pdfoutput>0 
  \pdfoutput=1
  \pdftrue
\fi 
\fi

\ifpdf
  \usepackage[pdftex]{epsfig}
  \usepackage[pdftex,colorlinks=true,%
              pdftitle={Twists of X(7) and primitive solutions to x^2+y^3=z^7},%
              pdfauthor={Bjorn Poonen, Ed Schaefer, Michael Stoll}]{hyperref}
  \newcommand{\Graph}[2]{\psfig{file=#1.pdf,width=#2}}
\else
  \usepackage{epsfig}
  \newcommand{\Graph}[2]{\psfig{file=#1.eps,width=#2}}
\fi

\usepackage[alphabetic,lite]{amsrefs} 

\begin{document}

\title[$x^2+y^3=z^7$]{Twists of $X(7)$ and primitive solutions to $x^2+y^3=z^7$}
\subjclass[2000]{Primary 11D41; Secondary 11G10, 11G18, 11G30, 14G05}
\keywords{generalized Fermat equation, Jacobian, Mordell-Weil sieve, Chabauty, modular curve, Klein quartic, descent}
\author{Bjorn Poonen}
\address{Department of Mathematics, University of California, Berkeley, CA  
94720-3840, USA}
\email{poonen@math.berkeley.edu}

\author{Edward F. Schaefer}
\address{Department of Mathematics and Computer Science,
Santa Clara University, Santa Clara, CA 95053, USA}
\email{eschaefe@math.scu.edu}

\author{Michael Stoll}
\address{School of Engineering and Science, International University Bremen,
         P.O.Box 750561, 28725 Bremen, Germany.}
\email{m.stoll@iu-bremen.de}

\date{July 30, 2005}

\begin{abstract}
We find the primitive integer solutions to $x^2+y^3=z^7$.
A nonabelian descent argument involving the simple group of order $168$
reduces the problem to the determination of the set of rational points
on a finite set of twists of the Klein quartic curve $X$.
To restrict the set of relevant twists, we exploit
the isomorphism between $X$ and the modular curve $X(7)$,
and use modularity of elliptic curves and level lowering.
This leaves $10$ genus-$3$ curves, whose rational points 
are found by a combination of methods.
\end{abstract}

\maketitle

\section{Introduction}
\label{introduction}

\subsection{Main result}

Fix integers $p,q,r \ge 1$.
An integer solution $(x,y,z)$ to $x^p+y^q=z^r$ is called {\em primitive}
if $\gcd(x,y,z)=1$.
The main goal of this paper is to prove:

\begin{theorem} \label{goal}
The primitive integer solutions to $x^2+y^3=z^7$
are the $16$ triples
\begin{gather*}
 (\pm 1, -1, 0),\quad (\pm 1, 0, 1),\quad \pm (0, 1, 1), 
  \quad (\pm 3, -2, 1),\quad (\pm 71, -17, 2), \\
 (\pm 2213459, 1414, 65),\quad (\pm 15312283, 9262, 113), 
  \quad (\pm 21063928, -76271, 17)\,.
\end{gather*}
\end{theorem}

\subsection{Previous work on generalized Fermat equations}

We temporarily return to the discussion of $x^p + y^q = z^r$
for a fixed general triple $(p,q,r)$.

Let $\chi = \frac{1}{p} + \frac{1}{q} + \frac{1}{r} - 1$.
There is a na\"{\i}ve heuristic that for each triple of integers $(x,y,z)$ with
$|x|\le B^{1/p}$, $|y| \le B^{1/q}$, $|z| \le B^{1/r}$,
the probability that $x^p + y^q = z^r$ is $1/B$:
this leads to the guesses that the set of primitive integer solutions
is finite when $\chi<0$ and infinite when $\chi>0$.
In fact, these statements are theorems, proved by 
Darmon and Granville~\cite{Darmon-Granville1995} and
Beukers~\cite{Beukers1998},
respectively.\footnote{
        In fact, they prove results more generally for equations
        $Ax^p+ By^q=C z^r$ where $A,B,C$ are nonzero integer constants.}

The method of {\em descent}, in one form or another,
is the key ingredient used in the proofs of these theorems.
Descent relates the primitive integer solutions to the rational points
(possibly over some number field) on one or more auxiliary curves,
whose Euler characteristic $2-2g$ is a positive integer multiple of $\chi$.
This connection is made clear in~\cite{Darmon1997Faltings}, for example.
In Section~\ref{S:descent theory},
we will explain it in detail 
in the context of the special case $(p,q,r)=(2,3,7)$.

Suppose $\chi>0$.
Then $2-2g>0$, so $g=0$.
Thus the auxiliary curves that have rational
points are isomorphic to $\PP^1$.
The conclusion is that the primitive integer solutions
fall into a finite number of parametrized families.
These parametrizations have been found explicitly for every $(p,q,r)$
with $\chi>0$, as we now indicate.
The case where some exponent is $1$ is trivial,
and the case where the multiset $\{p,q,r\}$ equals $\{2,2,n\}$
for some $n \ge 2$ is elementary.
The only remaining possibilities are $\{2,3,3\}$, $\{2,3,4\}$, $\{2,3,5\}$;
these were completed by Mordell, Zagier, and Edwards, respectively.
See~\cite{Edwards2004} for details.

If $\chi=0$, then $2-2g=0$, so $g=1$.
There are only the cases $\{2,3,6\}$, $\{2,4,4\}$, $\{3,3,3\}$,
and the auxiliary curves turn out to be elliptic curves $E$ over $\Q$,
and descent proves that $E(\Q)$ is finite in each case.
(The finiteness should not be surprising, 
since these elliptic curves have very low conductor.)
One then easily finds that the only nontrivial primitive solution,
up to permutations and sign changes, is $1^6 + 2^3 = 3^2$.
Many, if not all, of these cases are classical: 
Fermat invented the method of descent to do $(4,4,2)$,
and Euler
 did $(3,3,3)$.\footnote{
       Euler's argument had a flaw, but it was later fixed by others.}

If $\chi<0$, then $g>1$.
Each auxiliary curve has finitely many rational points,
by~\cite{Faltings1983},
and hence we get the Darmon-Granville result
that for each $(p,q,r)$ there are at most finitely many 
nontrivial primitive solutions.
For many $(p,q,r)$ these solutions have been determined explicitly.
These are summarized in Table~\ref{Table:pqr}.

\begin{table}
\begin{tabular}{ll}
$\{2,3,7\}$ & See Theorem~\ref{goal}. \\
$\{2,3,8\}$ & \cites{Bruin1999,Bruin2003}: 
                $1549034^2 + 33^8 = 15613^3$, 
                $96222^3 + 43^8 = 30042907^2$ \\
$\{2,3,9\}$ & \cite{Bruin2004}: $13^2 + 7^3 = 2^9$ \\
$\{2,2\ell,3\}$ & \cite{Chen2004preprint}, for prime $7 < \ell < 1000$ with $\ell \ne 31$ \\
$\{2,4,5\}$ & \cite{Bruin2003}: $7^2 + 2^5 = 3^4$, $11^4 + 3^5 = 122^2$ \\
$\{2,4,6\}$ & \cite{Bruin1999} \\
$(2,4,7)$ & \cite{Ghioca2002} \\
$(2,4,\ell)$ & \cite{Ellenberg2004Fermat}, for prime $\ell \ge 211$ \\
$(2,n,4)$ & follows easily from \cite{Bennett-Skinner2004}; \\ 
          & $(4,n,4)$ was done earlier in~\cite{Darmon1993CR} \\
$\{2,n,n\}$ & \cite{Darmon-Merel1997}; others for small $n$ \\
$\{3,3,4\}$ & \cite{Bruin2000} \\
$\{3,3,5\}$ & \cite{Bruin2000} \\
$\{3,3,\ell\}$ & \cite{Kraus1998}, for prime $17 \le \ell \le 10000$ \\
$\{3,n,n\}$ & \cite{Darmon-Merel1997}; others for small $n$ \\
$\{2n,2n,5\}$ & \cite{Bennett2004preprint} \\
$\{n,n,n\}$ & \cites{Wiles1995,Taylor-Wiles1995}, building on work of many others\\
\end{tabular}
\caption{
We list exponent triples $(p,q,r)$ such that $\chi<0$
and the primitive integer solutions 
to $x^p+y^q=z^r$ have been determined.  
The notation $\{p,q,r\}$ means that the solutions have been determined 
for every permutation of $(p,q,r)$: this makes a difference only if
at least two of $p,q,r$ are even.  
Solutions are listed only up to sign changes and permutations.
Solutions of the form $1^n + 2^3 = 3^2$ are omitted.
Heuristics suggest that probably there are no more solutions
for $(p,q,r)$ with $\chi<0$: in fact a cash prize awaits 
anyone who finds a new solution~\cite{Darmon1997Faltings}.
}
\hrule
\label{Table:pqr}
\end{table}

\subsection{Why $x^2+y^3=z^7$ is particularly difficult}

We give several answers.

\subsubsection{It corresponds to the negative value of $\chi$ closest to $0$}

When $(p,q,r)$ is such that $\chi$ is just barely negative,
the na\"{\i}ve heuristic mentioned earlier predicts
that the set of primitive integer solutions
should be ``just barely finite''
in the sense that it becomes likely that there are several
even of moderately large size.
The existence of solutions, especially large solutions,
naturally tends to make proving the nonexistence of others
a difficult task.

{}From this point of view, the triple $\{p,q,r\}$ with $\chi<0$
that should be most difficult
is the one with $\chi$ closest to $0$: this is $\{2,3,7\}$.
Indeed the equation $x^2+y^3=z^7$ has several solutions,
and some of them are large.
Now that this equation is finished, we know the complete list 
of solutions for every $(p,q,r)$ for which $\chi<0$
and for which nontrivial primitive solutions (excluding $1^n + 2^3 = 3^2$)
are known to exist.

\subsubsection{The exponents are prime}

The triple $\{2,3,8\}$ analyzed by Bruin has almost the same value of $\chi$,
and also has large solutions, so it also is difficult.
But at least it could be studied by using Zagier's parametrizations
for the case $\{2,3,4\}$ (with $\chi>0$).

\subsubsection{The exponents are pairwise relatively prime}

For descent, one requires finite \'etale covers,
and (geometrically) abelian \'etale covers can be easily constructed
in cases where two of the exponents share a common factor.
In contrast, for $\{2,3,7\}$,
the smallest nontrivial Galois \'etale covering
has a nonabelian Galois group of size $168$.
Before the present work, no pairwise coprime $(p,q,r)$ with $\chi<0$
had been studied successfully.

\section{Notation}

If $k$ is a field, then $\kbar$ denotes an algebraic closure,
and $\Gk = \Gal(\kbar/k)$.
Let $\Qbar$ be the algebraic closure of $\Q$ in $\C$.
Let $\zeta = e^{2 \pi i/7} \in \Qbar$.

If $X$ is a scheme over a ring $R$, and $R'$ is an $R$-algebra,
then $X_{R'}$ denotes $X \underset{\Spec R}\times \Spec R'$.
If $X$ is a variety over $\Q$,
and $R$ is a localization of $\Z$,
then $X_R$ will denote a smooth model of $X$ over $R$ (assuming it exists),
and $X(R)$ means $X_R(R)$.
Conversely, if $X_R$ is a scheme over such a ring $R$,
then $X$ will denote the generic fiber, a scheme over $\Q$.

We say that a $k$-variety $Y$ is a {\em twist} of a $k$-variety $X$
if $Y_\kbar \isom X_\kbar$.
Twists of morphisms are defined similarly.

If $M$ is a $\GQ$-module, and $k$ is an extension of $\Q$,
then $M_k$ denotes the same group with the action
of $\Gk$ via the restriction map $\Gk \rightarrow \GQ$.
In particular, $M_\Qbar$ denotes the abelian group with no Galois action.

\section{The theory behind the initial descent}
\label{S:descent theory}

This section explains the theory that reduces the original problem
to the problem of determining the rational points
on a finite list of auxiliary curves, which turn out to have genus $3$.

The argument is essentially the same as that in~\cite{Darmon-Granville1995}
(see also~\cite{Darmon1997Faltings}),
though we do have one small theoretical innovation:
By working with punctured affine surfaces
instead of ``M-curves'' 
(which are essentially the stack quotients of the 
surfaces by a $\G_m$-action), 
we can use the standard theory of descent
for finite \'etale covers of varieties \`a la Chevalley-Weil
and avoid having to use or develop more general statements
for ramified coverings or for stacks.
Keeping the stacks in the back of one's mind, however, still
simplifies the task of {\em finding} the \'etale covers 
of the punctured surfaces; we will use them for this purpose.

Another difference in our presentation is that we word the argument 
so as to obtain a finite list
of curves over $\Q$ instead of a single curve over a larger number field;
this is important for the practical computations.

Equations for the varieties and morphisms of this section
will appear in later sections of our paper.

\subsection{Scheme-theoretic interpretation of the problem}

Start with the subscheme $x^2+y^3=z^7$ in $\Aff^3_\Z$
and let $S_\Z$ be the (quasi-affine) scheme obtained by deleting
the subscheme $x=y=z=0$.
Then $S(\Z)$ is the set of primitive integer solutions to $x^2+y^3=z^7$.

\subsection{An \'etale cover of $S_\C$}

In order to apply the method of descent, 
we need a finite \'etale cover of $S_\Z$, 
or at least of $S_{\Z[1/N]}$ for some $N \ge 1$.
We first find a finite \'etale cover of $S_\C$.

The multiplicative group $\G_m$ acts on $S_\C$:
namely, $\lambda$ acts by 
$(x,y,z) \mapsto (\lambda^{21}x,\lambda^{14}y,\lambda^{6}z)$.
The subfield of the function field $\C(S_\C)$ invariant under $\G_m$
is $\C(j)$ where $j:=1728\, y^3/z^7$,
so the quotient of $S_\C$ by $\G_m$ should be birational to $\PP^1$
with coordinate $j$.
(The reason for the factor $1728$ and for the name $j$ will become clear
later on.)
But $\G_m$ does not act freely on $S_\C$:
the points where one of $x$, $y$, $z$ vanishes have nontrivial stabilizers.
Therefore, it is best to consider the stack quotient $[S_\C/\G_m]$,
which looks like $\PP^1$ except that the points $j=1728$, $j=0$, $j=\infty$
are replaced by a $\frac12$-point, a $\frac13$-point, and a $\frac17$-point,
respectively.
(This stack is essentially the same as the ``M-curve'' 
of~\cite{Darmon1997Faltings}.)
The induced morphism $[S_\C/\G_m] \stackrel{j}\to \PP^1$ 
is a ramified map of degree~1,
with ramification indices $2$, $3$, $7$ above $1728$, $0$, $\infty$,
respectively.

The strategy for constructing a finite \'etale cover $\Xtilde_\C$ of $S_\C$
will be first to construct a finite \'etale cover $X_\C$ of $[S_\C/\G_m]$
and then to base extend by $S_\C \to [S_\C/\G_m]$.
A finite \'etale cover of $[S_\C/\G_m]$
is a finite cover of $\PP^1$
whose ramification is accounted for by the ramification of
$[S_\C/\G_m] \to \PP^1$: in other words,
a finite cover of $\PP^1$ ramified only above $1728$, $0$, and $\infty$,
and with ramification indices $2$, $3$, $7$, respectively.
In fact, we will seek such a cover that is generically Galois.
Such covers can be specified analytically by giving
the monodromy above the three branch points,
and are automatically algebraic, by the Riemann Existence Theorem.
Specifically, the Galois group should be a {\em Hurwitz group}:
a finite group generated by elements $a$, $b$, $c$ 
of orders $2$, $3$, $7$, respectively, such that $abc=1$.
In principle we could use any Hurwitz group
(even the monster group is a possibility~\cite{Wilson2001}),
but to make future computations practical,
we choose the smallest Hurwitz group: 
the simple group $G:=\PSL_2(\F_7)$ of order $168$, with generators
$a = \left( \begin{smallmatrix} 0 & -1 \\ 1 & 0 \end{smallmatrix} \right)$,
$b = \left( \begin{smallmatrix} 0 & 1 \\ -1 & 1 \end{smallmatrix} \right)$,
$c = \left( \begin{smallmatrix} 1 & 1 \\ 0 & 1 \end{smallmatrix} \right)$.

We next identify the covering curve $X_\C$.
Since $[S_\C/\G_m]$ is obtained from $\PP^1$ by deleting three points
and reinserting a $\frac12$-point, a $\frac13$-point, and a $\frac17$-point,
its topological Euler characteristic satisfies
\begin{align*}
        \chi([S_\C/\G_m]) 
        &= \chi(\PP^1) - (1+1+1) + \left(\frac12 + \frac13 + \frac17 \right)\\
        &= \frac12 + \frac13 + \frac17 - 1 \\
        &= -\frac{1}{42}.
\end{align*}
The genus $g$ of its degree-$168$ cover $X_\C$ satisfies
\begin{align*}
        2-2g = \chi(X_\C) = 168\, \chi([S_\C/\G_m]) = -4.
\end{align*}
Hence $X_\C$ is a curve of genus $3$ with $168$ automorphisms;
the only such curve over $\C$ is the {\em Klein quartic curve}
defined in $\PP^2$ by the homogeneous equation
\[
        x^3 y + y^3 z + z^3 x = 0.
\]
See~\cite{Elkies1999} for many facts about $X$ and $G$.
The base extension of $X_\C$ by $S_\C \to [S_\C/\G_m]$
is the punctured cone $\Xtilde_\C$ above $X_\C$,
obtained by removing the origin from the surface in $\Aff^3$ defined by
the same quartic polynomial.

\subsection{An \'etale cover of $S_{\Z[1/42]}$}

Let $X$ and $\Xtilde$ be the varieties over $\Q$ defined by the same
equations, and let $X_R$ and $\Xtilde_R$ be the corresponding smooth models
over $R:=\Z[1/42]$.
We identify $G$ with $\Aut(X_\C)=\Aut(X_\Qbar)$,
so $\GQ$ acts on $G$.
The degree-$168$ \'etale map $\Xtilde_\C \to S_\C$
arises from an \'etale map $\pitilde\colon \Xtilde_R \to S_R$;
the associated morphism $X \to \PP^1$
is given explicitly in Section~\ref{S:canonical morphism}.

\begin{figure}
$$\xymatrix{
\Xtilde \ar[r] \ar@{~>}[d]_{\pitilde} & X \ar@{=}[r] \ar@{~>}[d] \ar[rd]^{\pi} & X(7) \ar[rd] & & & \hbox{\phantom{m}} \\
S \ar[r] & [S/\G_m] \ar[r]^{\text{birational}} & \PP^1 \ar@{=}[r] & X(1) & & \hbox{\phantom{m}} \\
(a,b,c) \ar@{|->}[rr] && j:= \dfrac{1728\, b^3}{c^7} \ar@{|->}[r] & E_{(a,b,c)}\colon \lefteqn{Y^2 = X^3 + 3b\,X - 2a} \\
}$$
\caption{The squiggly arrows denote \'etale morphisms.  All vertical and diagonal morphisms are finite of degree $168$.  The schemes may be considered over $\C$ or over $R=\Z[1/42]$, or anything in between.}
\label{F:etale covers}
\bigskip
\hrule
\end{figure}

\subsection{Descent}
\label{S:descent}

Descent theory 
(cf.~\cite{Chevalley-Weil1930} and \cite{Skorobogatov2001}*{\S5.3})
tells us the following:
\begin{itemize}
\item The set of classes in the Galois cohomology set $\HH^1(\Q,G)$
that are unramified outside $2$, $3$, and $7$ is finite.
(``Unramified outside $\{2,3,7\}$''
means that the image in $\HH^1(\Q_p^\unr,G)$
is the neutral element for every {\em finite} prime $p \notin \{2,3,7\}$.)
\item Each such class corresponds to an isomorphism class
of twists $\pitilde'\colon \Xtilde' \to S$
of $\pitilde\colon \Xtilde \to S$.
Choosing a cocycle in the class lets one construct a twist 
within the isomorphism class.
\item The set $S(R)$ is the disjoint union of the sets 
$\pitilde'(\Xtilde'(R))$ so obtained.
\end{itemize}
In fact, each such cocycle twists the upper half of the square
\[
\begin{CD}
        \Xtilde @>>> X \\
        @V{\pitilde}VV @VV{\pi}V \\
        S @>j>> \PP^1
\end{CD}
\qquad\text{ to yield }\qquad
\begin{CD}
        \Xtilde' @>>> X' \\
        @V{\pitilde'}VV @VV{\pi'}V \\
        S @>j>> \PP^1\lefteqn{.}
\end{CD}
\]
Since $\Pic R$ is trivial,
each map $\Xtilde'(R) \to X'(R) = X'(\Q)$ is surjective.
Hence $S(\Z) \subseteq S(R) = \Union j^{-1}(\pi'(X'(\Q)))$,
where the union is over our finite list of twists $X'$ of $X$.
So to solve the original equation, it would suffice to:
\begin{enumerate}
\item Find equations for each twist $X'$ and map $\pi'$
arising from a cocycle unramified outside $2$, $3$, $7$.
\item Determine $X'(\Q)$ for each $X'$.  (Since $X'$ is a genus-$3$ curve,
$X'(\Q)$ is finite by~\cite{Faltings1983}.)
\end{enumerate}
The first task reduces to a finite computation in principle,
because Hermite's theorem on number fields with bounded ramification
is effective.
The second task is not known to reduce to a finite computation,
and in fact for one of the curves $X'$ that arises,
we are unable to complete this task!
(See Section~\ref{S:Mordell-Weil sieve}.)

But since we want $S(\Z)$ and not all of $S(R)$,
we can get by with computing only the rational points inside
certain ``residue classes'' on each curve.
Say that $X'$ {\em passes the local test}
if the subset $\pi'(X'(\Q_p)) \cap j(S(\Z_p))$ of $\PP^1(\Q_p)$ 
is nonempty for all primes $p$.
We will do the following:
\begin{enumerate}
\item For each $X'$ passing the local test, 
find an equation for $X'$ and $\pi'$.
\item For each such $X'$, determine 
$\{ P \in X'(\Q) : \pi'(P) \in j(S(\Z_p)) \text{ for all primes $p$}\}$.
\end{enumerate}

\section{The modular interpretation}

In Section~\ref{S:X=X(7)} we observe that the map $X \to \PP^1$ 
is isomorphic to the map of modular curves $X(7) \to X(1)$.
In later sections this will be used to help catalogue the needed twists of $X$.
Indeed, each twist of $X(7)$ is a moduli space classifying elliptic
curves with an exotic level-$7$ structure,
so we reduce to a problem of enumerating level structures.

\subsection{Definition of $X(7)$}

If $E$ is an elliptic curve over a field $k$ of characteristic zero,
let $E[7]$ be the $\Gk$-module of $\kbar$-points
killed by the multiplication-by-$7$ map $[7]\colon E \rightarrow E$.
The Weil pairing gives an isomorphism $\bigwedge^2 E[7] \isom \mu_7$.
Let $M$ be the $\GQ$-module $\mu_7 \times \Z/7\Z$.\footnote{
Writing $M=\mu_7 \times \Z/7\Z$ instead of $M=\Z/7\Z \times \mu_7$
does not change it (even considering symplectic structure).
Our reason for this is so that when $M_\Qbar$ is identified with
column vectors, the action on $\SL_2(\F_7)$ on $M_\Qbar$
agrees with its usual action on $X(7)$, as groups
{\em equipped with a $\GQ$-action}.}
There is a canonical isomorphism $\bigwedge^2 M \isom \mu_7$
mapping $(\zeta^1,0) \wedge (\zeta^0,1)$ to $\zeta$.
When we write $\phi\colon M_k \wedgeisom E[7]$,
we mean that $\phi$ is a {\em symplectic isomorphism},
that is, an isomorphism of $\Gk$-modules
such that $\bigwedge^2 \phi\colon \bigwedge^2 M_k \rightarrow \bigwedge^2 E[7]$
is the identity $\mu_7 \rightarrow \mu_7$.

For a field $k$ of characteristic zero,
let $\YY(7)(k)$ be the set of isomorphism classes of pairs $(E,\phi)$
where $E$ is an elliptic curve over $k$
and $\phi\colon M_k \wedgeisom E[7]$.
It is possible to extend $\YY(7)$ to a functor
from $\Q$-schemes to sets, with essentially the same definition.
The functor $\YY(7)$ is representable 
by a smooth affine curve $Y(7)$ over $\Q$.
The curve $Y(7)$ can be completed to a smooth projective curve $X(7)$ 
over $\Q$.

\begin{remark}
The $\Q$-variety $X(7)$ 
extends to a smooth projective curve over $\Z[1/7]$
with geometrically integral fibers.
It too is a fine moduli space, representing a certain functor.
There is even a moduli scheme over $\Z$: see~\cite{Katz-Mazur1985}.
\end{remark}

\begin{warning}
Some authors use instead the set of isomorphism classes of pairs $(E,\phi)$
where $\phi\colon (\Z/7\Z)^2 \isom E[7]$ is an isomorphism of Galois modules
not necessarily respecting the symplectic structures.
This leads to an irreducible scheme over $\Q$
that over $\Q(\zeta)$ decomposes into $6$ disjoint components,
each isomorphic to our $X(7)_{\Q(\zeta)}$.
\end{warning}

\subsection{Automorphisms of $X(7)$ over $\Qbar$}

Let $Y(1)$ be the (coarse) moduli space of elliptic curves over $\Q$.
The $j$-invariant gives an isomorphism $Y(1) \isom \Aff^1$.
Let $X(1) \isom \PP^1$ be the usual compactification of $Y(1)$.
The forgetful functor mapping a pair $(E,\phi)$ to $E$
induces a morphism $Y(7) \rightarrow Y(1)$
and hence also a morphism $X(7) \rightarrow X(1)$.

The isomorphism $M_\Qbar \isom (\Z/7\Z)^2$
mapping $(\zeta^a,b)$ to
the column vector
$\left( \begin{smallmatrix} a \\ b \end{smallmatrix} \right)$
identifies $\bigwedge^2 M_\Qbar \isom (\mu_7)_\Qbar$
with $\bigwedge^2 (\Z/7\Z)^2 \, \stackrel{\det}\isom \, \Z/7\Z$.
Therefore the group of symplectic automorphisms of $M_\Qbar$ equals
        $$\Aut_\wedge(M_\Qbar) = \Aut_\wedge((\Z/7\Z)^2) = \SL_2(\F_7),$$
which we view as acting on the left on column vectors.
If $g \in \Aut_\wedge(M_\Qbar)$, then $g$ acts on the left on $Y(7)_\Qbar$
by mapping each pair $(E,\phi)$ to $(E, \phi \circ g^{-1})$;
this automorphism of $Y(7)_\Qbar$ extends to an automorphism of $X(7)_\Qbar$.
Thus we obtain a homomorphism
        $$\SL_2(\F_7) \rightarrow \Aut X(7)_\Qbar.$$
Since $\Aut_\wedge(M_\Qbar)$ acts transitively on the set of $\phi$
corresponding to a given $E$ over $\Qbar$,
the covering $X(7)_\Qbar \rightarrow X(1)_\Qbar$
is generically Galois,
and $\SL_2(\F_7)$ surjects onto the Galois group of the covering.
Pairs $(E,\phi)$ and $(E,\phi')$ are isomorphic if and only if
$\phi' = \alpha \circ \phi$ for some $\alpha \in \Aut(E)$;
but $\Aut(E)$ is generated by $[-1]$ for most $E$,
so the kernel of $\SL_2(\F_7) \rightarrow \Gal(X(7)_\Qbar/X(1)_\Qbar)$
is generated by the scalar matrix $-1 \in \SL_2(\F_7)$.
In other words, there is an isomorphism
$\PSL_2(\F_7) \isom \Gal(X(7)_\Qbar/X(1)_\Qbar)$.

\subsection{The Klein quartic as $X(7)$}
\label{S:X=X(7)}

There exists an isomorphism of $\Q$-varieties $X \to X(7)$ \cite{Elkies1999}.
The group of automorphisms of $X$ over $\Q$ has order $3$,
so the isomorphism $X \isom X(7)$ is not unique;
we will take the the particular isomorphism in~\cite{Elkies1999}
for which the homogeneous coordinates $x$, $y$, $z$
correspond on $X(7)$ to cusp forms whose $q$-expansions begin
\begin{align*}
        {\sf x} &= q^{4/7}(-1+4q-3q^2-5q^3+5q^4+8q^6 -10q^7+4q^9-6q^{10}
                        +\dots), \\
        {\sf y} &= q^{2/7}(1-3q-q^2+8q^3-6q^5-4q^6+2q^8+9q^{10} + \dots),
                        \text{ and} \\
        {\sf z} &= q^{1/7}(1-3q+4q^3+2q^4+3q^5-12q^6-5q^7+7q^9+16q^{10}
                        + \dots).
\end{align*}
Taking the quotient of each curve by its $\Qbar$-automorphism group,
we obtain the commutative parallelogram at the right in Figure~\ref{F:etale covers}.
Comparing the branch points and ramification indices
of $X(7) \to X(1)$ with $X \to \PP^1$
determines the values of the isomorphism $X(1) \to \PP^1$ at $3$ points;
it follows that the isomorphism 
is the standard one given by the $j$-invariant.

\subsection{Twists of $X(7)$}
\label{twistsofX(7)}

Since $\GQ$ acts on $M_\Qbar$,
it acts also on $\Aut_\wedge(M_\Qbar)$.
The homomorphism $\Aut_\wedge(M_\Qbar) \rightarrow \Aut(X(7)_\Qbar)$
is $\GQ$-equivariant, so we obtain a map of cohomology sets
        $$\HH^1(\GQ,\Aut_\wedge(M_\Qbar))
        \rightarrow \HH^1(\GQ,\Aut(X(7)_\Qbar)).$$
This map, which is not a bijection, can be reinterpreted as
\begin{equation}
\label{E:symplectic twists}
        \{\text{symplectic twists of $M$}\} \rightarrow
        \{\text{twists of $X(7)$}\}.
\end{equation}
(A {\em symplectic twist} of $M$ is
a $\GQ$-module $M'$ with an isomorphism
$\bigwedge^2 M' \isom \mu_7$
such that there is an isomorphism $\iota\colon M_\Qbar \rightarrow M'_\Qbar$
such that $\bigwedge^2 \iota$ is the identity $\mu_7 \rightarrow \mu_7$
over $\Qbar$.)
In \eqref{E:symplectic twists}, let $X_{M'}(7)$ be the image of $M'$.
Thus $X_{M'}(7)$ is the smooth projective model
of the smooth affine curve whose points
classify pairs $(E,\phi)$ where $E$ is an elliptic curve
and $\phi\colon M' \wedgeisom E[7]$.

Since $\Aut(X(7)_\Qbar)$
acts as automorphisms of $X(7)_\Qbar$ over $X(1)_\Qbar$,
there is a canonical morphism $X_{M'}(7) \to X(1)$.
In the moduli interpretation, this is the forgetful functor
mapping $(E,\phi)$ to $E$.

If $a \in \F_7^\times$, composing $\bigwedge^2 M' \isom \mu_7$
with the $a^{\operatorname{th}}$-power map $\mu_7 \rightarrow \mu_7$
gives a new isomorphism $\bigwedge^2 M' \isom \mu_7$.
Let ${M'}^a$ be $M'$ with this new symplectic structure.
If $a = b^2 \alpha$ for some $b \in \F_7^\times$,
then multiplication-by-$b$ induces
${M'}^a \wedgeisom {M'}^\alpha$.
Thus only the image of $a$ in $\F_7^\times/\F_7^{\times 2}$ matters:
we lose nothing by considering $a=\pm 1$ only.

If $M'$ is a symplectic twist of $M$,
and $M'' = M' \tensor \chi$ for some quadratic character
$\chi\colon\GQ \rightarrow \{\pm 1\}$,
then the elements of $\HH^1(\GQ,\Aut_\wedge(M_\Qbar))$
corresponding to $M'$ and $M''$
have the same image in $\HH^1(\GQ,\Aut(X(7)_\Qbar))$,
because $\{\pm1\}$ is in the kernel of
$\Aut_\wedge(M_\Qbar) \rightarrow \Aut(X(7)_\Qbar)$.

\subsection{Twists of $X(7)$ associated to elliptic curves}
\label{twistsofX(7) associated to elliptic curves}

Let $E$ be an elliptic curve over $\Q$.
Define $X_E(7):=X_{E[7]}(7)$ and $X_E^-(7):=X_{E[7]^{-1}}(7)$.

If $F$ is a quadratic twist of $E$,
then $X_E(7) \isom X_F(7)$
and $X_E^-(7) \isom X_F^-(7)$.

\begin{lemma}
\label{L:isomorphism of Galois modules}
If $F'$ is a quadratic twist of an elliptic curve $F$
such that $F[7] \isom E[7]$ as $\GQ$-modules 
{\em without symplectic structure},
then the point $j(F') \in X(1)(\Q)$
is in the image of $X_E(7)(\Q) \to X(1)(\Q)$
or $X_E^-(7)(\Q) \to X(1)(\Q)$.
\end{lemma}

\begin{proof}
  By the discussion above, we need to have $F[7] \wedgeisom E[7]$
  or $F[7] \wedgeisom E[7]^{-1}$, and so $F$ gives rise to a rational
  point on $X_E(7)$ or on $X_E^-(7)$. Therefore, $j(F') = j(F)$ is
  in the image of $X_E(7)(\Q) \to X(1)(\Q)$
  or of $X_E^-(7)(\Q) \to X(1)(\Q)$.
\end{proof}

\begin{warning}
Not every twist of $X(7)$ can be written as $X_E(7)$ or $X_E^-(7)$
for an elliptic curve $E$ over $\Q$.
For example, let $M':=\mu_7^{\tensor 2} \times \mu_7^{\tensor -1}$,
equipped with a symplectic structure $\bigwedge^2 M' \isom \mu_7$,
and let $X' = X_{M'}(7)$.
If $X'$ were isomorphic to $X_E^-(7)$, it would also be isomorphic
to $X_E(7)$, since $\Aut(M')$ acts transitively on the symplectic
structures on $M'$.
This is impossible, because $X'(\Q)=\emptyset$:
indeed, $X'(\F_2)=\emptyset$,
since there are no elliptic curves $E$ over $\F_2$
(not even degenerate ones corresponding to cusps)
with $E[7] \isom \mu_7^{\tensor 2} \times \mu_7^{\tensor -1}$,
that is, with the $2$-power Frobenius acting on $E[7]$
as multiplication-by-$4$.

In fact, $X(7)$ itself is another example of a twist
that cannot be written as $X_E(7)$ or $X_E^-(7)$
(at least if one requires $E$ to be a {\em nondegenerate} elliptic curve
over $\Q$):
again $\Aut(M)$ acts transitively on the symplectic structures
on~$M$, and although this time $X(7)$ has rational points,
they are all cusps~\cite{Elkies1999}.
\end{warning}

\subsection{Primitive solutions and elliptic curves}

Suppose $(a,b,c) \in S(\Z)$ and $a,b,c$ are nonzero.
Then $j:=1728\, b^3/c^7$ is not $1728$, $0$, or $\infty$.
The corresponding point on $X(1)$ is represented by the elliptic curve
\[
        E = E_{(a,b,c)} \colon Y^2 = X^3 + 3b\,X - 2a
\]
and its quadratic twists.
The twist $X'$ of $X$ in Section~\ref{S:descent}
for which $j \in \pi'(X'(\Q))$
corresponds to the twist $X_E(7)$ of $X(7)$.
Therefore to find the relevant twists of $X(7)$,
it suffices to catalogue the possibilities for 
the $\GQ$-module $E[7]$,
up to quadratic twist.
Restrictions on $E[7]$ will be deduced from the following.

\begin{lemma}
\label{L:conductor}
If $(a,b,c) \in S(\Z)$ and $a,b,c$ are nonzero,
then $1728\, b^3/c^7$ is the $j$-invariant of an elliptic curve
$E$ over $\Q$ such that
\begin{enumerate}
\item[(a)] The conductor $N$ of $E$ is of the form 
$2^r 3^s \prod_{p \in T} p$ where $r \le 6$, $s \le 3$,
and $T$ is a finite set of primes $\ge 5$.
\item[(b)] $E[7]$ is finite at $7$ in the sense of~\cite{Serre1987}*{p.~189}.
\end{enumerate}
\end{lemma}

\begin{proof}\hfill

(a)
Since $\gcd(a,b,c)=1$ and $a^2+b^3=c^7$,
the integers $a,b,c$ are {\em pairwise} relatively prime.
We need to bound the $p$-adic valuation $v_p(N)$ for each prime $p$.

\medskip
\noindent
{\em Case 1: $p \ge 5$}

Associated to the Weierstrass equation for $E_{(a,b,c)}$ are the quantities
\[ c_4 = -12^2\,b\,, \qquad
   c_6 = 12^3\,a\,, \qquad
   \Delta = -12^3\,c^7 \,.
\]
(See~\cite{SilvermanAEC}*{III.\S1} for definitions.)
Suppose $p \ge 5$ is a prime dividing $N$.
Then $p \mid \Delta$.
But $\gcd(b,c)=1$, so $p \notdiv c_4$.
By \cite{SilvermanAEC}*{VII.5.1(b)}, 
$E$ has multiplicative reduction at $p$, and $v_p(N)=1$.

\medskip
\noindent{\em Case 2: $p=3$}

If $3 \notdiv c$, then $v_3(N) \le v_3(\Delta) \le 3$.
Therefore suppose $3 \mid c$.
Then $3 \notdiv a,b$.
Let $t=a/b \in \Z_3$.
Then $a^2+b^3 \equiv 0 \pmod{3^7}$
implies $(a,b) \equiv (-t^3,-t^2) \pmod{3^7}$.
The right hand side of the Weierstrass equation is now
congruent to $(X-t)(X-t)(X+2t)$ modulo $3^7$,
so the substitution $X \mapsto x+t$
results in a cubic polynomial congruent to $x^3 + 3t x^2$
modulo $3^7$.
Hence we can twist by $1/3$ to obtain an elliptic curve $E'$ such that
\[
        v_3(c_4') = v_3(-12^2 b) -2 = 0 \quad \text{ and } \quad 
        v_3(\Delta') = v_3(-12^3 c^7) - 6 > 0.
\]
Thus $E'$ has multiplicative reduction, and $v_3(N')=1$.

\medskip
\noindent{\em Case 3: $p=2$}

If $2 \notdiv c$, then $v_2(N) \le v_2(\Delta) \le 6$.
Therefore suppose $2 \mid c$.
Twisting by $-1$ if necessary, we may assume $a \equiv 1 \pmod 4$.
As in the previous paragraph,
we have $(a,b) \equiv (-t^3,-t^2) \pmod{2^7}$ for some $t \in \Z_2$
and we can find a model $y^2 = f(x)$
where $f(x) \equiv x^3 + 3t x^2 \pmod{2^7}$.
Now $a \equiv 1 \pmod 4$ implies $t \equiv -1 \pmod 4$,
so the further substitution $y \mapsto y' + x$
shows that the model is not minimal at~$2$.
Thus we arrive at a model $E'$ with
\[
        v_2(c_4') = v_2(-12^2 b) - 4 = 0 \quad \text{ and } \quad 
        v_2(\Delta') = v_2(-12^3 c^7) - 12 > 0,
\]
and hence $v_2(N')=1$.

(b)
If $7 \notdiv c$, then $E$ has good reduction at $7$,
so $E[7]$ is finite at $7$.
If $7 \mid c$, then $v_7(j) = -7 v_7(c)$ is a negative multiple of $7$,
so by the theory of Tate curves, $E[7]$ is finite at $7$.
(The twists done in part (a) affect the behavior
at $2$ and $3$ independently, 
and do not interfere with the behavior at other primes.)
\end{proof}

\section{Reducible $7$-torsion}

By Lemma~\ref{L:conductor} 
we need only catalogue the possibilities for $E[7]$ (up to quadratic twist)
such that $E[7]$ is unramified outside $\{2,3,7\}$.
In this section we assume moreover that $E[7]$ is reducible as a $\GQ$-module.

\subsection{Intermediate modular curves}

Let $Y_\mu(7)$ over $\Q$ be the affine curve parametrizing
pairs $(E,\psi)$ where $E$ is an elliptic curve and
$\psi\colon \mu_7 \hookrightarrow E[7]$ is an injection.
Let $X_\mu(7)$ be the smooth projective curve having $Y_\mu(7)$ as
an affine subset.

\medskip

\begin{remark}
The curve $X_\mu(7)$ is a twist of the classical modular curve $X_1(7)$,
the completion of the curve $Y_1(7)$ parametrizing $(E,P)$
where $P$ is a point of exact order~$7$ on $E$:
giving $P$ is equivalent to giving an injection
$\Z/7\Z \hookrightarrow E[7]$.
In a less obvious way, $X_\mu(7)$ is also isomorphic to $X_1(7)$:
if $(E,\psi)$ as above is defined over $k$,
one can map it to the pair $(E',P')$
where $E'$ is the elliptic curve $E/\psi(\mu_7)$ over $k$
and $P' \in E'[7](k)$ is the image of a point $P \in E[7](\kbar)$
whose Weil pairing with $\psi(\zeta)$ gives $\zeta$.
\end{remark}

\medskip

Define $Y_0(7)$ as the (coarse) moduli space
parametrizing pairs $(E,C)$
where $C$ is a (cyclic) order-$7$ subgroup of $E$,
and define $X_0(7)$ as the completed curve.

\begin{figure}
$$\xymatrix{
{X(7)} \ar_7[d] & {(E,\phi)} \ar@{|->}[d] 
  & {\{1\} = \begin{pmatrix} 1 & 0 \\ 0 & 1 \end{pmatrix}} & 1 \\
{X_\mu(7)} \ar_3[d] & {(E,\phi|_{\mu_7})} \ar@{|->}[d] & {N =  
\begin{pmatrix} 1 & * \\ 0 & 1 \end{pmatrix}} & 7 \\
{X_0(7)} \ar_8[d] & {(E,\phi(\mu_7))} \ar@{|->}[d] & {B = \begin{pmatrix} *  
& * \\ 0 & * \end{pmatrix}} & 21 \\
{X(1)} & E & {G = \begin{pmatrix} * & * \\ * & * \end{pmatrix}} & 168 \\
}$$
\caption{Modular curves of level~$7$.}
\label{modulartower}
\bigskip
\hrule
\end{figure}
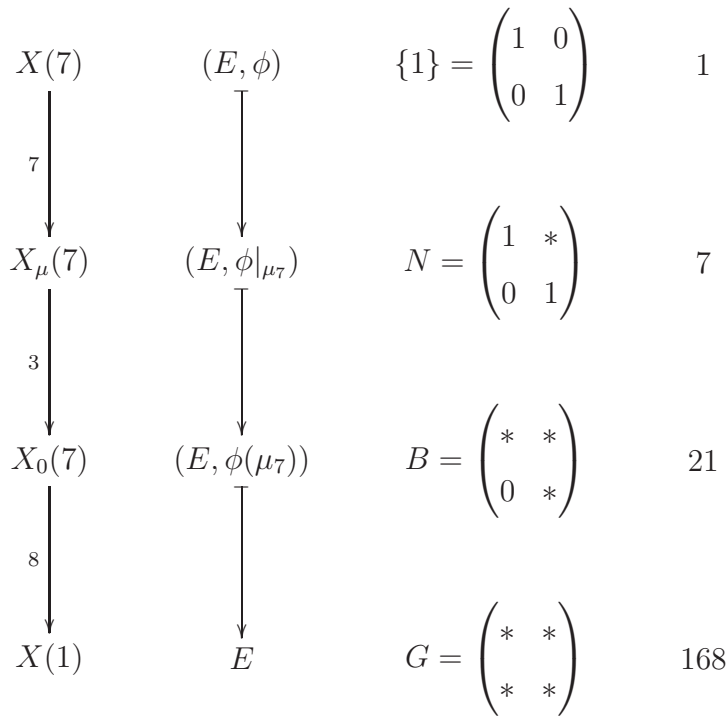

The columns of Figure~\ref{modulartower} show
\begin{itemize}
\item the morphisms relating the modular curves, labelled by their degrees,
\item the moduli interpretations of these morphisms,
\item the subgroup $\Aut(X(7)_\Qbar/Z_\Qbar)$ of
$\Aut(X(7)_\Qbar/X(1)_\Qbar)=G$ corresponding
to each intermediate curve $Z$, and
\item the orders of these subgroups.
\end{itemize}
The notation
$B = \left( \begin{smallmatrix} * & * \\ 0 & * \end{smallmatrix} \right)$
means that $B$ equals the image of
$\left\{ \left(\begin{smallmatrix} a & b \\ 0 & d \end{smallmatrix} \right)
\in \SL_2(\F_7) \right\}$
in $\PSL_2(\F_7) \isom G$.

\subsection{The Galois action on $G$ and its subgroups}

Since $M = \mu_7 \times \Z/7\Z$,
we have
        $$\End(M_\Qbar) \isom
        \begin{pmatrix}
                \End(\mu_7) & \Hom(\Z/7\Z,\mu_7) \\
                \Hom(\mu_7,\Z/7\Z) & \End(\Z/7\Z)
        \end{pmatrix}
        \isom
        \begin{pmatrix}
                \Z/7\Z & \mu_7 \\
                \mu_7^{\tensor -1} & \Z/7\Z
        \end{pmatrix},$$
which is the matrix ring $\M_2(\F_7)$ with a nontrivial $\GQ$-action.
This $\GQ$-action induces a $\GQ$-action on $\PSL_2(\F_7)$
making $\PSL_2(\F_7) \isom G$ an isomorphism of {\em $\GQ$-groups}
(groups equipped with a continuous action of $\GQ$).

In particular, $N \isom \mu_7$ as $\GQ$-groups,
and $B$ is a semidirect product of
the subgroup
$H := \left( \begin{smallmatrix} * & 0 \\ 0 & * \end{smallmatrix} \right)  
\isom \F_7^\times/\{\pm 1\} \isom \Z/3\Z$
(with trivial $\GQ$-action)
by $N$.
Thus we have a split exact sequence of $\GQ$-groups
\begin{equation}
\label{Bsequence}
        0 \rightarrow \mu_7 \rightarrow B \rightarrow \Z/3\Z \rightarrow 0.
\end{equation}

The generators
$g:=\left( \begin{smallmatrix} 1 & 1 \\ 0 & 1 \end{smallmatrix} \right)$
of $N \isom \mu_7$
and
$h:=\left( \begin{smallmatrix} 2 & 0 \\ 0 & 4 \end{smallmatrix} \right)$
of $H \isom \Z/3\Z$
satisfy $h^{-1}gh = g^2$.
Under the isomorphism in Section~\ref{S:X=X(7)},
$g$ and $h$ correspond to the automorphisms
$(x:y:z) \mapsto (\zeta^4 x: \zeta^2 y: \zeta z)$
and
$(x:y:z) \mapsto (y:z:x)$
of $X$,
according to~\cite{Elkies1999}.

\subsection{Nonabelian cohomology of $B$}
\label{nonabeliancohomology}

Let $E$ be an elliptic curve over $\Q$ such that $E[7]$ is reducible
as a $\GQ$-module.
Thus there is a $\GQ$-stable order-$7$ subgroup $C \subset E[7]$.
Choosing $\phi\colon M_\Qbar \wedgeisom E[7]_\Qbar$
such that $\phi(\mu_7)=C$
shows that the twist $X_E(7)$ of $X(7)$
arises from a cocycle taking values in the subgroup $B$ of $G$.
In this subsection, we calculate $\HH^1(\GQ,B)$;
in the next subsection, we use the results
to write down equations for all $X_E(7)$
corresponding to $E$ with reducible $E[7]$.

{}From the exact sequence~(\ref{Bsequence}) of $\GQ$-groups,
we obtain the exact sequence of pointed
sets~\cite{SerreGaloisCohomology}*{Proposition~38,~p.~51}
        $$\HH^1(\GQ,\mu_7) \rightarrow \HH^1(\GQ,B)
        \stackrel{s}\rightarrow \HH^1(\GQ,\Z/3\Z).$$
Define $s$ as shown.
Since the last map in~(\ref{Bsequence}) admits a section,
$s$ admits a section; in particular $s$ is surjective.

Suppose $\xi$ is a cocycle representing a class $[\xi] \in \HH^1(\GQ,\Z/3\Z)$.
Since $\Z/3\Z$ acts on $\mu_7$,
$\xi$ determines a twist ${}_\xi \mu_7$ of $\mu_7$.
A cohomologous cocycle would determine a (noncanonically) isomorphic twist.
Via the section of $B \rightarrow \Z/3\Z$, we may view $\xi$
also as a $B$-valued cocycle.
By~\cite{SerreGaloisCohomology}*{Corollary~2,~p.~52},
there is a natural surjective map
$\HH^1(\GQ,{}_\xi\mu_7) \rightarrow s^{-1}([\xi])$.
We now calculate $\HH^1(\GQ,{}_\xi\mu_7)$.

\bigskip

\noindent\underline{\em Case 1: $[\xi]=0$}.

Then $\HH^1(\GQ,{}_\xi\mu_7) = \HH^1(\GQ,\mu_7) \isom \Q^\times/\Q^{\times 7}$.

\bigskip

\noindent\underline{\em Case 2: $[\xi] \not=0$}.

Then $[\xi] \in \HH^1(\GQ,\Z/3\Z) = \Hom_{\conts}(\GQ,\Z/3\Z)$
corresponds to a cyclic cubic extension $K$ of $\Q$
with a choice of generator $\tau \in \Gal(K/\Q)$.
{}From the Hochschild-Serre
spectral sequence~\cite{SerreGaloisCohomology}*{p.~15}
we can extract the exact sequence
\begin{align}
\label{hochschildserre}
        & \HH^1(\Gal(K/\Q),({}_\xi\mu_7)^{\GK})
        \rightarrow \HH^1(\GQ,{}_\xi\mu_7)
        \rightarrow \HH^1(\GK,{}_\xi\mu_7)^{\Gal(K/\Q)} \\
\nonumber
        & \rightarrow \HH^2(\Gal(K/\Q),({}_\xi\mu_7)^{\GK}).
\end{align}
The first and last terms are zero,
since $\Gal(K/\Q)$ and ${}_\xi\mu_7$ have relatively prime orders.
In the third term,
$\HH^1(\GK,{}_\xi\mu_7)$ equals $K^\times/K^{\times 7}$ with a
twisted action of $\Gal(K/\Q)$;
thus $\HH^1(\GK,{}_\xi\mu_7)^{\Gal(K/\Q)}$,
instead of equalling the subgroup of $K^\times/K^{\times 7}$ {\em fixed} by
the usual action of $\tau$,
equals the subgroup $(K^\times/K^{\times 7})^{\tau=2}$ of $K^\times/K^{\times 7}$
consisting of elements on which $\tau$ has the same effect
as squaring.
(The~$2$ comes from the action of $\Z/3\Z$ on $\mu_7$,
that is, from the fact that $h^{-1}gh=g^2$.)
Thus~(\ref{hochschildserre}) gives
        $$\HH^1(\GQ,{}_\xi\mu_7)
        \isom \HH^1(\GK,{}_\xi\mu_7)^{\Gal(K/\Q)}
        \isom (K^\times/K^{\times 7})^{\tau = 2}.$$

\subsection{Twists of $X(7)$ corresponding to reducible $E[7]$}
\label{writingdowntwists}

We now write down twists of $X \isom X(7)$
whose classes in $\HH^1(\GQ,G)$
represent the image of $\HH^1(\GQ,B) \rightarrow \HH^1(\GQ,G)$.

\bigskip

\noindent\underline{\em Case 1: $[\xi]=0$}.

Given $a \in \Q^\times$, fix $\alpha \in \Qbar$ with $\alpha^7 = a$.
The substitution $(x,y,z) = (\alpha^4 \xbar, \alpha^2 \ybar, \alpha \zbar)$
transforms $X_\Qbar$ into the plane curve
        $$(\alpha^4 \xbar)^3 (\alpha^2 \ybar)
        + (\alpha^2 \ybar)^3 (\alpha \zbar)
        + (\alpha \zbar)^3 (\alpha^4 \xbar) = 0,$$
which, if we divide by $a$ and drop the bars, is the same as
\begin{equation}
\label{Case1twist}
        X'\colon a x^3 y + y^3 z + z^3 x = 0.
\end{equation}
Let $\beta\colon X_\Qbar \rightarrow X'_\Qbar$ be the isomorphism
we just described.
If we identify
the automorphism $(x:y:z) \mapsto (\zeta^4 x: \zeta^2 y : \zeta z)$
of $X_\Qbar$
with
$g:=\left( \begin{smallmatrix} 1 & 1 \\ 0 & 1 \end{smallmatrix} \right) \in N$ 
and with $\zeta \in \mu_7$,
then the cocycle $\sigma \mapsto \beta^{-1} ({}^\sigma \beta) \in N$
is identified with the cocycle
$\sigma \mapsto {}^\sigma \alpha^{-1} / \alpha^{-1} \in \mu_7$.
Thus the image of $a^{-1}$ under
        $$\Q^\times \rightarrow \Q^\times/\Q^{\times 7} \isom \HH^1(\GQ,\mu_7)
        \isom \HH^1(\GQ,N) \rightarrow \HH^1(\GQ,G)$$
corresponds to the twist $X'$ of $X$.
In particular, all twists of $X$ arising in Case~1 have
the form~(\ref{Case1twist}),
where $a$ runs through positive integers
representing the elements of $\Q^\times/\Q^{\times 7}$.

The same curves can be written in other ways.
Suppose that $X'$ is a curve given by
        $$ a x^3 y + b y^3 z + c z^3 x = 0 $$
for some $a,b,c \in \Z_{>0}$.
If $p$ is a prime such that $p^2$ divides $a$,
then multiplying by $p$ and substituting $x = \xbar/p$
leads to an isomorphic curve of the same type
but with a smaller value of $abc$.
Similar reductions apply if $p^2$ divides $b$ or $c$.
Hence eventually we can reduce to a curve where
$a$, $b$, and $c$ are squarefree.
Moreover we may assume that $\gcd(a,b,c)=1$
and that $a \ge b$ and $a \ge c$.

\bigskip

\noindent\underline{\em Case 2: $[\xi] \not=0$}.

We will write the action of $\tau$ on $K$ on the right,
so that we can use notation such as $a^{2\tau-3}:=(a^\tau)^2/a^3$
for $a \in K^\times$.
Let $\ell \in K \xbar + K \ybar + K \zbar$
be a linear form such that $\ell$, $\ell^\tau$, and $\ell^{\tau \tau}$
form a $K$-basis for $K \xbar + K \ybar + K \zbar$.
Let $a \in K^\times$ be such that its image in $K^\times/K^{\times 7}$
lies in $(K^\times/K^{\times 7})^{\tau = 2}$.
Thus $a^{\tau-2} = b^7$ for some $b \in K^\times$.
Applying $\tau^2+2\tau+4$ and using $\tau^3=1$ yields
$a^{-7} = b^{7 \tau^2 + 14 \tau + 28}$.
Since $\mu_7(K)=\{1\}$, it follows that
$a^{-1} = b^{\tau^2+2\tau+4}$.

Choose $\alpha \in \Qbar$ such that $\alpha^7=a$.
Define elements $\alpha_\tau := b \alpha^2$
and $\alpha_{\tau \tau} := b^{\tau+2} \alpha^4$ of $\Qbar$.
(There is no action of $\tau$ on $\alpha$, so we interpret
$\alpha_\tau$ and $\alpha_{\tau \tau}$ as new labels.)
Then $\alpha_\tau^7 =a^\tau$ and
$\alpha_{\tau \tau}^7 = a^{\tau \tau}$.

A calculation shows that the substitution
$(x,y,z)=(\alpha_{\tau \tau} \ell^{\tau \tau},
\alpha_\tau \ell^\tau , \alpha \ell)$
transforms $X_\Qbar$ into the plane curve
\begin{equation}
\label{Case2twist}
        X_{c,\ell} \colon c^{\tau \tau} \ell^{3 \tau \tau + \tau}
        + c^{\tau} \ell^{3\tau+1}
        + c \ell^{3 + \tau \tau}
        = 0,
\end{equation}
where $c = b^{\tau+2} a \in K^\times$.
By symmetry, the coefficients of the polynomial in $\xbar$, $\ybar$, $\zbar$
defining $X_{c,\ell}$ are fixed by $\tau$,
and hence are in $\Q$.
A chase through definitions (similar to that in Case~1, but more tedious)
shows that the $B$-valued cocycle defined by the
constructed isomorphism $X_\Qbar \rightarrow (X_{c,\ell})_\Qbar$
represents the element of $s^{-1}([\xi])$
coming from the image of $a^{-1}$ in $(K^\times/K^{\times 7})^{\tau = 2}$.

In summary,
all twists of $X$ arising in Case~2 have the form $X_{c,\ell}$
for some cyclic cubic extension $K$,
some choice of $\tau \in \Gal(K/\Q)$,
some choice of $\ell$,
and some choice of $c \in K^\times$.
Changing $\ell$ does not change $X_{c,\ell}$.

Changing $a$ to $\tilde{a} = d^7 a$ for some $d \in K^\times$,
so that the corresponding values of $b$ and $c$ are
$\tilde{b}=d^{\tau-2} b$ and $\tilde{c}=d^{\tau^2+3} c$,
does not change $X_{c,\ell}$.
Thus multiplying $c$ by an element of $K^{\times (\tau^2+3)}$
does not change $X_{c,\ell}$. (This amounts to changing $\ell$ into $d \ell$,
resulting in a linear substitution on $(\xbar, \ybar, \zbar)$ with rational
coefficients.)
Multiplying $c$ by an element of $K^{\times (\tau^2+\tau+1)} \subseteq \Q^\times$
simply multiplies the equation of $X_{c,\ell}$ by an element of $\Q^\times$,
so this also does not change $X_{c,\ell}$.
Since the difference of $\tau^2+\tau+1$ and $\tau^2+3$
equals $\tau-2$, and since $7$ is a multiple of $\tau-2$ in
$\Z[\tau]/(\tau^3-1)$,
multiplying $c$ by an element of $K^{\times 7}$
or more generally, an element of $K^{\times (\tau-2)}$,
does not change $X_{c,\ell}$.
Since the natural map $(K^\times/K^{\times 7})^{\tau=2} \to K^\times/K^{\times (\tau-2)}$
is an isomorphism, 
it suffices to let $c$ run through representatives of the elements of 
$(K^\times/K^{\times 7})^{\tau=2}$.

\subsection{Ramification conditions}

In our application, we know by Lemma~\ref{L:conductor}
that $E[7]$ is unramified outside $\{2,3,7\}$,
just as $M$ is.
Therefore any isomorphism $M_\Qbar \wedgeisom E[7]_\Qbar$
will be unramified outside $\{2,3,7\}$.
The choice of isomorphism gives the cocycle
defining the twist $X_E(7)$ of $X(7)$,
so the cocycle will be unramified outside $\{2,3,7\}$,
and then so will be all the cocycles in Section~\ref{nonabeliancohomology}.

\bigskip

\noindent\underline{\em Case 1: $[\xi]=0$}.

The twist has the form~(\ref{Case1twist}) for some $a \in \Z_{>0}$.
We need consider only $a$
whose image in $\Q^\times/\Q^{\times 7} \isom \HH^1(\GQ,\mu_7)$
is unramified outside $\{2,3,7\}$.
Equivalently, the field $\Q(a^{1/7})$ should be unramified
outside $\{2,3,7\}$.
Concretely, this means that the only primes occurring in
the factorization of $a$ are among $2$, $3$, and $7$.

\bigskip

\noindent\underline{\em Case 2: $[\xi] \not=0$}.

The number field $K$ must be unramified above all
finite places outside $\{2,3,7\}$.
By the Kronecker-Weber theorem, the only such cyclic cubic extensions
are the four degree-$3$ subfields of
the $(\Z/3\Z)^2$-extension
$\Q(\zeta+\zeta^{-1},\zeta_9+\zeta_9^{-1})$,
where $\zeta=e^{2\pi i/7}$ (as usual) and $\zeta_9=e^{2\pi i/9}$.
Let $\OO_K$ be the ring of integers of $K$.
The element $a \in K^\times$
must be such that $K(a^{1/7})$ is unramified above
finite places outside $\{2,3,7\}$.
Since each $K$ has class number~1,
it is possible to multiply $a$ by an element of $K^{\times 7}$
to assume that $a \in \OO_K[1/42]^\times$.
Then the element $c$ of Case~2 of Section~\ref{writingdowntwists}
lies in $\OO_K[1/42]^\times$, and we can finally reduce to the finite
group $(\OO_K[1/42]^\times/\OO_K[1/42]^{\times 7})^{\tau=2}$,
generators for which are easily computed for each $(K, \tau)$.

\section{Irreducible $7$-torsion: level lowering}

Now we list the possibilities for irreducible $E[7]$, up to quadratic twist.
Let $\calE$ be the following set of $13$ elliptic curves over $\Q$ in the
notation of~\cite{Cremona1997}:
\begin{gather*}
  \text{24A1}, \text{27A1}, \text{32A1}, \text{36A1}, \text{54A1}, 
  \text{96A1}, \text{108A1}, \text{216A1}, \text{216B1}, \text{288A1}, 
  \text{864A1}, \text{864B1}, \text{864C1}.
\end{gather*}

\begin{lemma}
\label{levelloweringlemma}
Suppose that $E$ is as in Lemma~\ref{L:conductor},
and that $E[7]$ is an irreducible $\GQ$-module.
Then there exists a quadratic twist $E'$ of some $E'' \in \calE$
such that $E[7] \isom E'[7]$ as $\GQ$-modules.
\end{lemma}

\begin{proof}
By~\cite{Breuil2001}, there is a weight-$2$ newform $f$ on $\Gamma_0(N)$
associated to $E$, where $N$ is the conductor of $E$.
By Lemma~\ref{L:conductor}, $N=2^r 3^s \prod_{p \in T} p$
for some $r \le 6$, $s \le 3$ and finite set $T$ of primes $\ge 5$.
Let $N'=2^r 3^s$.
Since $E[7]$ is irreducible,
Theorem~1.1 of \cite{Ribet1990} 
lets us lower the level of $f$ by
eliminating $7$ (if present) and then the other primes $\ge 5$:
more precisely, 
there is a weight-$2$ newform $f'$ on $\Gamma_0(N')$ such that
``$f' \equiv f \pmod{7}$.''
The congruence is to be interpreted as follows:
``The $\Z$-algebra generated by the Fourier coefficients of~$f'$
is an order $\OO_{f'}$ in some number field,
and there is a prime ideal $\pp \subset \OO_{f'}$ of norm~$7$
such that $a_p(f') \bmod {\mathfrak p} = a_p(f) \bmod 7$ 
for all but finitely many~$p$.''
Replacing~$f'$ by a quadratic twist does not change $\OO_{f'}$.

We now consult William Stein's Modular Forms Database~\cite{SteinTables}
(or use Magma \cite{Magma}),
and get all weight-2 newforms~$f'$ of level dividing $2^6 3^3$.
If two of them are quadratic twists of each other,
we can check this by looking at sufficiently many coefficients,
because the twisting factor is in $\{\pm1,\pm2,\pm3,\pm6\}$,
and the order of vanishing of a difference of modular forms 
at the cusp $\infty$ can be bounded in terms of its level.
We find that each $f'$ is a quadratic twist of one of $14$ newforms $f''$.
Now $13$ of these $14$ newforms have $\OO_{f''}=\Z$
and are associated to the elliptic curves $E'' \in \calE$.
(Warning: the labelling of isogeny classes of~\cite{SteinTables}
is different from that of~\cite{Cremona1997} at some levels.)
The $14^{\operatorname{th}}$ has $\OO_{f''}$ an order containing $\sqrt{13}$
in $\Q(\sqrt{13})$
(in fact it is $\Z[\sqrt{13}]$);
such an order has no prime of norm~$7$,
so we can discard this last newform.

Thus $f \equiv f' \pmod{7}$
where $f'$ is associated to an elliptic curve $E'$
having a quadratic twist $E'' \in \calE$.
By the Eichler-Shimura relation 
(see \cite{Eichler1954} for the relation in the context we need),
the congruence implies that for all but finitely many $p$,
the traces of the Frobenius automorphism at $p$ 
acting on $E[7]$ and $E'[7]$ are equal.
By the Chebotarev density theorem, the $\GQ$-representations
$E[7]$ and $E'[7]$ have the same trace.
They also have the same determinant, namely the cyclotomic character.
So the Brauer-Nesbitt theorem~\cite{Curtis-Reiner1988}*{Theorem~30.16} 
implies that the semisimplifications of $E[7]$ and $E'[7]$ are isomorphic.
But $E[7]$ is irreducible and hence semisimple,
so $E[7] \isom E'[7]$ as $\GQ$-modules.
\end{proof}

\begin{corollary}
\label{C:irreducibleE[7]}
  Suppose $(a,b,c) \in S(\Z)$ and $a,b,c$ are nonzero.
  If $E_{(a,b,c)}[7]$ is irreducible,
  then $j=1728 \, b^3/c^7 \in X(1)(\Q)$ is the image of a rational point
  on one of the $26$ curves $X_E(7)$ or $X^-_E(7)$ with $E \in \calE$.
\end{corollary}

\begin{proof}
Combine Lemmas \ref{L:isomorphism of Galois modules}
and~\ref{levelloweringlemma}.
\end{proof}

\section{Explicit equations}

Each twist of $X=X(7)$ is a nonhyperelliptic curve of genus~$3$,
and hence has a model as a smooth plane curve of degree~$4$.
In this section we find equations for the homogeneous forms
defining the twists that we need,
and also the equations for their canonical morphisms to $\PP^1=X(1)$.
The strategy is to exploit the $168$ symmetries.

\subsection{Covariants of ternary quartic forms}

The standard left action of $\GL_n(\C)$ on
$V = \C^n$ induces a right action of $\GL_n(\C)$ on the $\C$-algebra
$\C[x_1,\dots,x_n] = \Sym\, V^*$.
If $g \in \GL_n(\C)$ and $F \in \C[x_1,\dots,x_n]$,
let $F^g$ be the result of applying $g$ to $F$,
so $F^g(v)=F({}^g v)$ for all $v \in V$.
Let $\C[x_1,\dots,x_n]_d = \Sym^d V^*$
be the subspace of $\C[x_1,\dots,x_n]$
consisting of homogeneous polynomials of degree $d$.

Fix $n$ and $d$.
A {\em covariant} of {\em order} $j$ and {\em degree} $\delta$ is a function
$\Psi\colon \C[x_1,\dots,x_n]_d \rightarrow \C[x_1,\dots,x_n]_j$
such that
\begin{enumerate}
\item
The coefficient of a fixed monomial $x_1^{i_1}\dots x_n^{i_n}$
in $\Psi(F)$ depends polynomially on the coefficients of $F$,
as $F$ varies.
\item
For each $g \in \SL_n(\C)$, we have $\Psi(F^g) = \Psi(F)^g$.
\item
For each $t \in \C^\times$, we have $\Psi(tF)=t^\delta \Psi(F)$.
\end{enumerate}

In the same way, a {\em contravariant} of {\em order}~$j$
and {\em degree}~$\delta$ is a polynomial map 
$\Psi\colon \Sym^d V^* \to \Sym^j V$
that is equivariant with respect to
the right action of $\SL_n(\C)$ on the two spaces
and is homogeneous of degree~$\delta$ in the coefficients.
(The right action of $\SL_n(\C)$ on~$V$ is
the contragredient of the action on~$V^*$.)
In other words, for all $F \in \Sym^d V^*$
and $g \in \SL_n(\C)$ 
we have $\Psi(F^g) = \Psi(F)^{g^{-t}}$,
where $g^{-t}$ is the inverse transpose of $g$.
Our choice of basis for $V$ induces an isomorphism $V \to V^*$,
so we may express each contravariant as a polynomial
in the same variables $x_1,\dots,x_n$.

{}From now on, we take $n=3$ and $d=4$,
and use $x,y,z$ in place of $x_1,x_2,x_3$.
We will list some covariants, which we will call
$\Psi_0$, $\Psi_4$, $\Psi_6$, $\Psi_{14}$, and $\Psi_{21}$, and a  
contravariant, which we will call~$\Psi_{-4}$.
They will be of orders $0$, $4$, $6$, $14$, $21$, $4$
and of degrees $3$, $1$, $3$, $8$, $12$, $2$, respectively.

First of all, define $\Psi_4(F)=F$.
Then $\Psi_6(F)$ is essentially the Hessian of~$F$:
\[ \Psi_6(F) = -\frac{1}{54}\,\left|\begin{matrix}
       \frac{\partial^2 F}{\partial x^2} &
         \frac{\partial^2 F}{\partial x\,\partial y} &
         \frac{\partial^2 F}{\partial x\,\partial z} \\
       \frac{\partial^2 F}{\partial y\,\partial x} &
         \frac{\partial^2 F}{\partial y^2} &
         \frac{\partial^2 F}{\partial y\,\partial z} \\
       \frac{\partial^2 F}{\partial z\,\partial x} &
         \frac{\partial^2 F}{\partial z\,\partial y} &
         \frac{\partial^2 F}{\partial z^2} \end{matrix}\right|
\]
The constant $-1/54$ is there so that $\Psi_6$ has integer coefficients
with gcd~$1$.
Next, $\Psi_{14}(F)$ is obtained in a similar way:
\[ \Psi_{14}(F)  = \frac{1}{9}\,\left|\begin{matrix}
        \frac{\partial^2 F}{\partial x^2} &
         \frac{\partial^2 F}{\partial x\,\partial y} &
         \frac{\partial^2 F}{\partial x\,\partial z} &
         \frac{\partial \Psi_6(F)}{\partial x} \\
       \frac{\partial^2 F}{\partial y\,\partial x} &
         \frac{\partial^2 F}{\partial y^2} &
         \frac{\partial^2 F}{\partial y\,\partial z} &
         \frac{\partial \Psi_6(F)}{\partial y } \\
       \frac{\partial^2 F}{\partial z\,\partial x} &
         \frac{\partial^2 F}{\partial z\,\partial y} &
         \frac{\partial^2 F}{\partial z^2} &
         \frac{\partial \Psi_6(F)}{\partial z} \\
       \frac{\partial \Psi_6(F)}{\partial x} &
         \frac{\partial \Psi_6(F)}{\partial y} &
         \frac{\partial \Psi_6(F)}{\partial z} &
         0 \end{matrix}\right|
\]
And $\Psi_{21}(F)$ is essentially the Jacobian determinant of the other three:
\[ \Psi_{21}(F) = \frac{1}{14}\,\left|\begin{matrix}
       \frac{\partial F}{\partial x} &
         \frac{\partial F}{\partial y} &
         \frac{\partial F}{\partial z} \\
       \frac{\partial \Psi_6(F)}{\partial x} &
         \frac{\partial \Psi_6(F)}{\partial y} &
         \frac{\partial \Psi_6(F)}{\partial z} \\
       \frac{\partial \Psi_{14}(F)}{\partial x} &
         \frac{\partial \Psi_{14}(F)}{\partial y} &
         \frac{\partial \Psi_{14}(F)}{\partial z} \end{matrix}\right|
\]
The {\em invariant} $\Psi_0(F)$ can be obtained as follows. Consider the
following differential operator.
\[ D = \left|\begin{matrix}
       \frac{\partial}{\partial x_1} &
         \frac{\partial}{\partial y_1} &
         \frac{\partial}{\partial z_1} \\
       \frac{\partial}{\partial x_2} &
         \frac{\partial}{\partial y_2} &
         \frac{\partial}{\partial z_2} \\
       \frac{\partial}{\partial x_3} &
         \frac{\partial}{\partial y_3} &
         \frac{\partial}{\partial z_3} \end{matrix}\right|
\]
Then
\[ \Psi_0(F)
     = \frac{1}{5184} D^4 (F(x_1,y_1,z_1) F(x_2,y_2,z_2) F(x_3,y_3,z_3)) \,.
\]
Alternatively, one can find an explicit formula in
Salmon's book~\cite{Salmon1879}*{p.~264},
where he gives also a formula for the contravariant $\Psi_{-4}(F)$.
We do not reproduce it here, but in Section~\ref{equations}
we will give the result for the special curves we will have to deal with.

(The covariants $\Psi_j$ do not generate the $\C$-algebra of covariants,
but they are all we will need. In fact, they do generate the specialization
of the algebra of covariants for each ternary quartic 
$\GL_3(\C)$-equivalent to $x^3 y + y^3 z + z^3 x$.)

Let $F_0 = x^3 y + y^3 z + z^3 x$ be the polynomial defining $X$.
Each $g \in G \subset \GL_3(\C)$ must act on $F_0$
as multiplication by a scalar,
so there is a character $\chi\colon G \rightarrow \C^\times$ 
describing the action of $G$ on $F_0$.
But $G$ is simple and nonabelian, so $\chi$ is trivial.
In other words, $F_0$ belongs to $\C[x,y,z]^G$, the ring of $G$-invariants.
For $j \in \{0,4,6,14,21\}$, let $\Phi_j = \Psi_j(F_0) \in \C[x,y,z]_j$.
Since $\Psi_j$ is a covariant, $\Phi_j \in \C[x,y,z]^G$.
In fact, the $\Phi_j$ generate $\C[x,y,z]^G$ as a $\C$-algebra.
Moreover, $\Phi_0=1$,
the polynomials $\Phi_4$, $\Phi_6$, $\Phi_{14}$ are algebraically
independent, and $\Phi_{21}^2$ can be expressed as a polynomial in the others.
We give this relation modulo $\Phi_4$, since this is all we need:
\begin{equation}
\label{invariantrelation}
        \Phi_{21}^2 - \Phi_{14}^3 \equiv -1728\,\Phi_0 \Phi_6^7 \pmod{\Phi_4}.
\end{equation}

\begin{lemma}
\label{covariantrelation}
If $F = F_0^g$ for some $g \in \GL_3(\C)$,
then
        $$\Psi_{21}(F)^2 - \Psi_{14}(F)^3
                \equiv -1728\,\Psi_0(F) \Psi_6(F)^7  \pmod{F}.$$
\end{lemma}

\begin{proof}
Write $g=t g'$ with $t \in \C^\times$ and $g' \in \SL_3(\C)$,
and apply it to~(\ref{invariantrelation}).
The result follows, since the two sides of~(\ref{covariantrelation})
are covariants of the same degree~$24$ and order~$42$.
(This was the point of including the trivial factor $\Phi_0$
in~(\ref{covariantrelation}).)
\end{proof}

\begin{lemma}
\label{L:no common zeros}
Let $F = F_0^g$ for some $g \in \GL_3(\C)$.
No two of $\Psi_{6}(F)$, $\Psi_{14}(F)$, $\Psi_{21}(F)$
vanish simultaneously at a point of the curve $F=0$ in $\PP^2_\C$.
\end{lemma}

\begin{proof}
We may assume $F=F_0$, 
in which case the result is in \cite{Klein1879b}*{\S6}:
see \cite{Elkies1999}*{p.~66}.
\end{proof}

\subsection{Equations for $X_E(7)$ and $X_E^-(7)$} \label{equations}

In order to use Corollary~\ref{C:irreducibleE[7]},
we need explicit equations for the $26$ curves described there.
Halberstadt and Kraus~\cite{Halberstadt-Kraus2003}
find an explicit equation for $X_E(7)$ in terms of the coefficients of~$E$.
If $E$ is given as $Y^2 = X^3 + a\,X + b$, the equation is
\begin{equation}
\label{X_E(7)equation}
  a\,x^4 + 7b\,x^3 z + 3\,x^2 y^2 - 3a^2\,x^2 z^2 - 6b\,x y z^2 - 5ab\,x z^3
    + 2\,y^3 z + 3a\,y^2 z^2 + 2a^2\,y z^3 - 4b^2\,z^4 = 0 \,.
\end{equation}

The following result will let us deduce an equation for $X_E^-(7)$
from \eqref{X_E(7)equation}, using very little additional work.

\begin{proposition}
\label{P:equation for X_E^-(7)}
If $F$ is a ternary quartic form describing $X_{M'}(7)$, 
then $\Psi_{-4}(F)$ is a ternary quartic form describing $X_{M'}^-(7)$.
\end{proposition}

\begin{proof}
Given $M'$, the curve $X_{M'}(7)$ can be constructed as follows.
Fix a symplectic isomorphism 
$\phi\colon M_\Qbar \stackrel{\wedgeisom}{\to} M'_\Qbar$.
Then 
$\sigma \mapsto \xi_\sigma = \phi^{-1} \circ \ssigma \phi$ defines
a $1$-cocycle on~$\GQ$ with values in $\Aut_\wedge(M_\Qbar)$.
We have
\[
        \Aut_\wedge(M_\Qbar) \surjects G = \Aut(X_\Qbar) \injects \GL_3(\Qbar),
\]
in which the last homomorphism is the $3$-dimensional representation
$(V,\rho)$ of \cite{Elkies1999}*{p.~54} giving the
action of $\Aut(X_\Qbar)$ on 
the dual of the space of holomorphic differentials on~$X$,
and giving also the automorphisms 
of $\PP^2_\Qbar$ that restrict to the automorphisms of 
the canonically embedded curve $X_\Qbar$.
Under this composition, 
$\xi$ maps to a cocycle $\bar{\xi}\colon \GQ \to \GL_3(\Qbar)$,
which can be written as $\sigma \mapsto g^{-1} \cdot \ssigma g$ 
for some $g \in \GL_3(\Qbar)$,
since $\HH^1(\GQ,\GL_3)=0$.
If $\sim$ denotes ``equal up to scalar multiple'', then
\[
        \ssigma\left(F_0^{g^{-1}}\right) = \left( \ssigma F_0 \right)^{\left( \ssigma g^{-1} \right)} = F_0^{\bar{\xi}_\sigma^{-1} g^{-1}} \sim F_0^{g^{-1}}
\]
since $\bar{\xi}_\sigma^{-1}$ induces an automorphism of $X_\Qbar$.
Thus $F_0^{g^{-1}}=0$ defines a twist $X'$ of $X$ in $\PP^2$.
The cocycle $\GQ \to \Aut(X_\Qbar)$ defining this twist 
(obtained from the isomorphism $X_\Qbar \to X'_\Qbar$ given by $g$)
is the same as the cocycle defining $X_{M'}(7)$, 
so $X' \isom X_{M'}(7)$.
After replacing $F$ by $F^b$ for some $b \in \GL_3(\Q)$
(which we may do, without loss of generality),
we have $F_0^{g^{-1}} \sim F$.

The element 
$\iota = \left(\begin{smallmatrix} 1 & 0 \\ 0 & -1 \end{smallmatrix}\right)$
defines an anti-symplectic automorphism of $M$ over $\Q$.
A calculation based on the explicit description of $(V,\rho)$
in \cite{Elkies1999}*{p.~54} shows that
$\Aut M_\Qbar \to \GL_3(\Qbar)$ is equivariant
with respect to conjugation by $\iota$ on the left
and $g \mapsto g^{-t}$ on the right.
In particular, the image of $G$ in $\SL_3(\Qbar)$ is stable under
inverse transpose,
so the contravariance of $\Psi_{-4}$ implies that
$\Psi_{-4}(F_0)$ is a $G$-invariant ternary quartic form.
All such are multiples of $F_0$ \cite{Elkies1999}*{pp.55--56},
so $\Psi_{-4}(F_0) \sim F_0$.

Now we repeat the first paragraph of this proof 
using the anti-symplectic isomorphism $\phi \iota \colon M_\Qbar \to M'_\Qbar$
in place of $\phi$, to construct $X_{M'}^-(7)$.
The image of the cocycle 
\[
        \sigma \mapsto (\phi \iota)^{-1} \circ \ssigma(\phi \iota) = \iota^{-1} \xi_\sigma \iota
\]
under $\Aut M_\Qbar \to \GL_3(\Qbar)$ 
is
\[
        \sigma \mapsto \bar{\xi}_\sigma^{-t} = \left(g^{-1} \cdot \ssigma g \right)^{-t} = h^{-1} \cdot \ssigma h,
\]
where $h:=g^{-t}$,
so the polynomial defining $X_{M'}^-(7)$ is 
\[
        F_0^{h^{-1}} \sim \Psi_{-4}(F_0)^{h^{-1}} \sim \Psi_{-4}(F_0^{h^t}) = \Psi_{-4}(F_0^{g^{-1}}) \sim \Psi_{-4}(F).
\]
\end{proof}

A formula for $\Psi_{-4}$ is given in \cite{Salmon1879}*{p.~264}.
Applying it to~(\ref{X_E(7)equation}),
we obtain an equation for $X_E^-(7)$:
\begin{align}
\label{X_E^-(7)equation}
 -a^2\,x^4 &+ a (3 a^3+19 b^2)\,y^4 + 3\,z^4 + 6 a^2\,y^2 z^2 + 6 a\,z^2 x^2
 - 6 (a^3 + 6 b^2)\,x^2 y^2 - 12 a b\,y^2 z x \\
\nonumber
 {} &+ 18 b\,z^2 x y + 2 a b\,x^3 y
 - 12 b\,x^3 z - 2 (4 a^3 + 21 b^2)\,y^3 z + 2 a^2 b\,y^3 x - 8 a\,z^3 y = 0
\end{align}

\subsection{Equations for the degree-$168$ morphism to $\PP^1$}
\label{S:canonical morphism}

Recall that each twist $X'$ of $X=X(7)$ comes equipped with a canonical
degree-$168$ morphism to $\PP^1=X(1)$ that is invariant under $\Aut(X'_\Qbar)$.
The morphism with these properties is unique if we require that 
its branch points (with ramification indices $2$, $3$, $7$)
be at $1728$, $0$, $\infty$, respectively, on $\PP^1$.

\begin{lemma}
\label{L:canonical morphism}
Let $X'$ be a twist of $X$ defined by a ternary quartic form $F$.
Then the canonical morphism is 
\begin{align*}
        \pi' : \quad\qquad X' & \To \PP^1 \\
        (x:y:z) &\longmapsto \frac{\Psi_{14}(F)^{3}}{\Psi_0(F)\,\Psi_6(F)^7} \, ,
\end{align*}
where as usual we use the parameter $j$ on $\PP^1$.
\end{lemma}

\begin{proof}
We need only to check this for~$X$, 
in which case it is~(2.13) in \cite{Elkies1999}.
\end{proof}

\subsection{The local test} 
\label{S:local test}

Let $X'$ be a twist of the Klein quartic~$X$.
Then $X'$ is given by an equation $F=0$,
where $F = F_0^g \in \Q[x,y,z]_4$ for some $g \in \GL_3(\Qbar)$.
By Lemma~\ref{L:canonical morphism},
there is a morphism $\pi' \colon X' \to \PP^1$, given by
        $$(x : y : z) \mapsto
  \frac{\Psi_{14}(F)^{3}}{\Psi_0(F) \, \Psi_6(F)^7} \,.$$

We showed in Section~\ref{S:descent} that every $(a,b,c) \in S(\Z)$
maps to a point in $\pi'(X'(\Q))$ for some twist $X'$
whose class in $\HH^1(\GQ,G)$ is unramified outside $2$, $3$, and~$7$.
Now suppose we have such a twist.
Then $\pi'$ extends to a morphism $X'_R \rightarrow \PP^1_R$,
so $\pi'(X'(\Q)) = \pi'(X'(R)) \subseteq \PP^1_R(R)$,
and the corresponding solutions to $x^2 + y^3 = z^7$
will have $\gcd(x,y,z)$ divisible by at most primes in $\{2, 3, 7\}$.
For each of the remaining primes $p = 2, 3, 7$,
we search for residue classes that give rise to $p$-primitive solutions.
(We say that $(x,y,z)$ in $\Z^3$ or $\Z_p^3$ is {\em $p$-primitive}
if at least one of $x$, $y$,~$z$ is not divisible by $p$.)

Let us define a {\em residue class} in $\PP^2(\Q_p)$ to be a subset
that, after permuting the homogeneous coordinates, equals
$\{(1 : a + p^k u : b + p^l v) : u,v \in \Z_p\}$
for some fixed $a, b \in \Z_p$ and integers $k,l \ge 0$.
We have an initial partition of $\PP^2(\Q_p)$
into the three residue classes $(1 : \Z_p : \Z_p)$, $(p \Z_p : 1 : \Z_p)$
and $(p \Z_p : p \Z_p : 1)$. For each given residue class,
we rewrite $F$ in terms of $u$ and~$v$. If this expression has the form
$p^e(\alpha + p F_1(u,v))$ with $\alpha \in \Z_p^\times$ and $e \in \Z$,
then the residue class does not meet $X'(\Q_p)$, so we can discard it.
Otherwise, we compute the covariants $\Psi_6(F)$, $\Psi_{14}(F)$
and $\Psi_{21}(F)$ in terms of $u$ and~$v$. The minimum of the $p$-adic
valuations of the coefficients gives us a lower bound for the $p$-adic
valuation of the covariants on our residue class, and if the constant
term has the unique coefficient with minimal valuation, we even know
the $p$-adic valuation of the corresponding covariant exactly.
Let $W_6$, $W_{14}$ and $W_{21}$ be these lower bounds for
$v_p(\Psi_6(F))$, $v_p(\Psi_{14}(F))$ and $v_p(\Psi_{21}(F))$, respectively,
and let $w = v_p(1728 \Psi_0(F))$.

One of the solutions to $a^2 + b^3 = c^7$ corresponding to the point
$(x : y : z) \in X'(\Q)$ is given by
\[ a = (1728 \Psi_0(F))^3 \Psi_{21}(F)\,, \qquad
   b = -(1728 \Psi_0(F))^2 \Psi_{14}(F)\,, \qquad
   c = -1728 \Psi_0(F) \Psi_6(F) \,,
\]
evaluated at $(x,y,z)$. This can be scaled to a $p$-primitive solution
if and only if
\[ \min\{v_p(a)/21, v_p(b)/14, v_p(c)/6\} \in \Z \,.\]
So we define
\[ w_6 = (W_6 + 3w)/21\,, \qquad w_{14} = (W_{14} + 2w)/14\,, \qquad
   w_{21} = (W_{21} + w)/6 \,.
\]
If $\min\{w_6, w_{14}, w_{21}\}$ is an integer {\em and} the minimum comes
from a covariant whose valuation on the residue class is known exactly,
then there will be a $p$-adic point on~$X'$ giving rise to a $p$-primitive
solution (at least if the given residue class contains a $p$-adic point
on~$X'$, which we can check using Hensel's Lemma).
In the same way, if $\min\{w_6, w_{14}, w_{21}\}$ is not an integer
{\em and} the minimum comes from a covariant whose valuation on the residue
class is known exactly, then no $p$-adic point in this residue class can
produce a $p$-primitive solution. If the minimum comes only from a lower
bound, we split the residue class into $p$ (or~$p^2$) subclasses
and repeat. In this way, we can decide fairly quickly whether any given twist
gives rise to primitive solutions everywhere locally. If so, we say that
``$X'$ passes the local test.''

This algorithm terminates, since 
by Lemma~\ref{L:no common zeros}
no two of $\Psi_6(F)$, $\Psi_{14}(F)$
and $\Psi_{21}(F)$ vanish simultaneously on $X'(\Q_p)$.

\subsection{The list}

Putting together the twists coming from the reducible and irreducible cases,
respectively, we obtained an initial finite list of twists.
We then ran the local test on each of the curves.
There are exactly $10$ twists that pass the local test. 
They are given in Table~\ref{Table:curve equations}. 
We have put some effort into simplifying
these equations (trying to minimize the invariant $\Psi_0(F)$ and to
reduce the size of the coefficients); the equations given
are usually not those coming from formula~\eqref{X_E(7)equation} directly.

\begin{table}[h]
\begin{align*}
  C_1 &\colon 6 x^3 y + y^3 z + z^3 x = 0 \\
  C_2 &\colon 3 x^3 y + y^3 z + 2 z^3 x = 0 \\
  C_3 &\colon 3 x^3 y + 2 y^3 z + z^3 x = 0 \\
  C_4 &\colon 7 x^3 z + 3 x^2 y^2 - 3 x y z^2 + y^3 z - z^4 = 0 \\
  C_5 &\colon {-2} x^3 y - 2 x^3 z + 6 x^2 y z + 3 x y^3 - 9 x y^2 z
         + 3 x y z^2 - x z^3 + 3 y^3 z - y z^3 = 0 \\
  C_6 &\colon x^4 + 2 x^3 y + 3 x^2 y^2 + 2 x y^3 + 18 x y z^2 + 9 y^2 z^2 - 9 z^4 = 0 \\
  C_7 &\colon {-3} x^4 - 6 x^3 z + 6 x^2 y^2 - 6 x^2 y z + 15 x^2 z^2 - 4 x y^3
         - 6 x y z^2 - 4 x z^3 + 6 y^2 z^2 - 6 y z^3 = 0 \\
  C_8 &\colon 2 x^4 - x^3 y - 12 x^2 y^2 + 3 x^2 z^2 - 5 x y^3 - 6 x y^2 z + 2 x z^3 
         - 2 y^4 + 6 y^3 z + 3 y^2 z^2 + 2 y z^3 = 0 \\
  C_9 &\colon 2 x^4 + 4 x^3 y - 4 x^3 z - 3 x^2 y^2 - 6 x^2 y z + 6 x^2 z^2 - x y^3
         - 6 x y z^2 - 2 y^4 + 2 y^3 z \\
      &\qquad {} - 3 y^2 z^2 + 6 y z^3 = 0 \\
  C_{10} &\colon x^3 y - x^3 z + 3 x^2 z^2 + 3 x y^2 z + 3 x y z^2 + 3 x z^3 - y^4 
         + y^3 z + 3 y^2 z^2 - 12 y z^3 + 3 z^4 = 0
\end{align*}
\caption{Equations for the plane quartic curves $C_1$ through $C_{10}$.}
\label{Table:curve equations}
\medskip 
\hrule
\end{table}

Curves $C_1$, $C_2$, and $C_3$ come from the case that $E[7]$ is reducible.
As it happens, they all correspond to twists that involve only the
cyclic subgroup of order~$7$. 
Curves $C_4$ through $C_{10}$ 
are of the form $X_E(7)$ for $E = \text{27A1}$, 96A1, 288A1,
864A1, 864B1, 864C1, 54A2, respectively. 
(The curve $C_{10}$ is also $X_E^-(7)$ for $E = \text{54A1}$,
because 54A1 and 54A2 are $3$-isogenous, and $3$ is not a square
modulo $7$.)

\subsection{Local conditions for $C_5$}

If $P_1,P_2 \in \PP^2(\Q_p)$ then we write ``$P_1 \equiv P_2 \pmod{p}$''
to mean that the images of $P_1$ and $P_2$
under the reduction map $\PP^2(\Q_p) = \PP^2(\Z_p) \to \PP^2(\F_p)$
are equal.
Define
\begin{align*}
        C_5(\Q_2)_\subs &:= \{\,P \in C_5(\Q_2) : 
                        P \equiv (1:0:0) \text{ or } (1:1:1) \pmod{2}\,\} \\
        C_5(\Q_3)_\subs &:= \{\,P \in C_5(\Q_3) : 
                        P \equiv (0:1:0) \pmod{3}\,\} \\
        C_5(\Q)_\subs &:= \{\,P \in C_5(\Q) : P \in C_5(\Q_2)_\subs 
                        \text{ and } P \in C_5(\Q_3)_\subs \,\}.
\end{align*}

\begin{lemma}
\label{L:C5 constraints}
Let $p$ be $2$ or $3$.
If $P \in C_5(\Q_p)$ 
and the image of $P$ in $\PP^1(\Q_p)$ is in $j(S(\Z_p))$,
then $P \in C_5(\Q_p)_\subs$.
\end{lemma}

\begin{proof}
We do the computations for the ``local test'', keeping track of the
residue classes that are not excluded. For example, when $p = 3$,
the reduction mod~$3$ of the equation $F_5=0$ of~$C_5$ is
\[ (x - z)(x^2y + y^2z + xyz + xz^2 + yz^2) = 0 \,, \]
which describes a line and a cuspidal cubic intersecting at the cusp
$(0 : 1 : 0)$. The invariant $\Psi_0(F_5)$ equals $-24$, so 
in the notation of Section~\ref{S:local test}, $w=4$.
For all residue classes corresponding to smooth points on the
reduction, we find that the minimum of $w_6$, $w_{14}$ and $w_{21}$
is determined and non-integral. So $(0 : 1 : 0)$ is the only 
residue class that may contain points leading to $3$-primitive
solutions. The argument for $p = 2$ is similar.
\end{proof}

Because of Lemma~\ref{L:C5 constraints},
it will suffice to determine $C_5(\Q)_\subs$ instead of $C_5(\Q)$.
Eventually we will prove $C_5(\Q)_\subs = \emptyset$.

\section{Known rational points on the $10$ curves}

We can forget most of the paper up to now.
All we need to keep are the equations for $C_1$ through $C_{10}$,
the local conditions for $C_5$ given by Lemma~\ref{L:C5 constraints},
and the formula in Lemma~\ref{L:canonical morphism}.
(The latter is needed only to recover the solutions
to $x^2+y^3=z^7$ from the relevant rational points on the $10$ curves).

\begin{table}
\begin{center}
\renewcommand{\arraystretch}{1}
\begin{tabular}{c|c||cccc}
      & Rank & Points & $j$-invariant & Cremona & Primitive solutions \\ \hline \hline
$C_1$ & $1$  & $(1:0:0)$ & $\infty$ & & $(\pm 1, -1, 0)$ \\
      &      & $(0:1:0)$ & $\infty$ & & $(\pm 1, -1, 0)$ \\
      &      & $(0:0:1)$ & $\infty$ & & $(\pm 1, -1, 0)$ \\
      &      & $(1:-1:2)$ & $-13^3 221173^3 2^{-2} 3^{-1} 43^{-7}$ & 258F2 & 
                 ---
                 \\ \hline
$C_2$ & $1$  & $(1:0:0)$ & $\infty$ & & $(\pm 1, -1, 0)$ \\
      &      & $(0:1:0)$ & $\infty$ & & $(\pm 1, -1, 0)$ \\
      &      & $(0:0:1)$ & $\infty$ & & $(\pm 1, -1, 0)$ \\
      &      & $(1:1:-1)$ & $-7^4 2^{-1} 3^{-1}$ & 294A1 & --- \\
      &      & $(1:-2:-1)$ & $-13^3 26293^3 2^{-2} 3^{-1} 113^{-7}$ & 678D2 & 
                --- \\ \hline
$C_3$ & $1$  & $(1:0:0)$ & $\infty$ & & $(\pm 1, -1, 0)$ \\
      &      & $(0:1:0)$ & $\infty$ & & $(\pm 1, -1, 0)$ \\
      &      & $(0:0:1)$ & $\infty$ & & $(\pm 1, -1, 0)$ \\
      &      & $(1:1:-1)$ & $-7^3 44647^3 2^{-1} 3^{-1} 29^{-7}$ & 174B2 & 
                --- \\ \hline \hline
$C_4$ & $2$  & $(1:0:0)$ & $0$ & & $(\pm 1, 0, 1)$ \\
      &      & $(0:1:0)$ & $0$ & & $(\pm 1, 0, 1)$ \\
      &      & $(0:1:1)$ & $-2^{15} 3 \cdot 5^3$ & 27A2 & --- \\ \hline
$C_5$ & $3$  & $(1:0:0)$ & $2^6 13^3 3^{-1}$ & 96A2 & --- \\
      &      & $(0:1:0)$ & $2^3 23^3 3^{-4}$ & 96A4 & --- \\
      &      & $(0:0:1)$ & $2^3 97^3 3^{-1}$ & 96A3 & --- \\
      &      & $(1:1:1)$ & $2^6 7^3 3^{-2}$ & 96A1 & --- \\ \hline
$C_6$ & $2$  & $(0:1:0)$ & $2^6 3^3$ & & $\pm(0, 1, 1)$ \\
      &      & $(1:-1:0)$ & $2^6 3^3$ & &  $\pm(0, 1, 1)$ \\
      &      & $(0:1:1)$ & $2^9 3^3 11^3 421^3 113^{-7}$ & 32544* & $(\pm 15312283, 9262, 113)$ \\
      &      & $(0:1:-1)$ & $2^9 3^3 11^3 421^3 113^{-7}$ & 32544* & $(\pm 15312283, 9262, 113)$ \\ \hline
$C_7$ & $2$  & $(0:1:0)$ & $-2^3 3^3$ & 864A1 & --- \\
      &      & $(0:0:1)$ & $2^9 3^3 7^3 101^3 5^{-7} 13^{-7}$ & 56160* & $(\pm 2213459, 1414, 65)$ \\
      &      & $(0:1:1)$ & $-2^6 3^3 13^3 5867^3 17^{-7}$ & 14688D1 & $(\pm 21063928, -76271, 17)$ \\ \hline
$C_8$ & $2$  & $(0:0:1)$ & $-2^9 3^3$ & 864B1 & $(\pm 3,-2,1)$ \\
      &      & $(2:-1:0)$ & $2^3 3 \cdot 547^3 66029^3 977^{-7}$ & 
               844128* & --- \\ \hline
$C_9$ & $2$  & $(0:0:1)$ & $2^9 3$ & 864C1 & --- \\
      &      & $(1:1:0)$ & $2^3 3^9 163^3 8779^3 7^{-7} 79^{-7}$ & 
               477792* & --- \\ \hline
$C_{10}$&$2$ &$(1:0:0)$ & $-3^3 17^3 2^{-1}$ & 54A3 & $(\pm 71,-17,2)$ \\
      &      & $(1:1:0)$ & $-3 \cdot 73^3 2^{-9}$ & 54A2 & --- \\ \hline
\end{tabular}
\end{center}
\medskip
\caption{Known rational points on $C_1$ through $C_{10}$.
}
\label{Table:rational points}
\hrule
\end{table}

Table~\ref{Table:rational points} gives the rank of the
Mordell-Weil group of the Jacobian of each curve
(these numbers will be proved correct in Section~\ref{S:2-descent}),
and lists rational points 
on $C_1$ through $C_{10}$ discovered by a na\"{\i}ve search.
We will prove that the list is complete for each curve 
except possibly $C_5$.\footnote{
   Almost certainly the list is complete for $C_5$ too,
   since Noam Elkies ran an algorithm based on his paper~\cite{ElkiesANTS4}
   to prove that there are no additional points $(x:y:z)$ 
   with integer homogeneous coordinates satisfying $|x|,|y|,|z| \le 10^7$.
}
For $C_5$, we will prove instead that $C_5(\Q)_\subs = \emptyset$,
which will suffice for our purposes.
If we prove all this, then Theorem~\ref{goal} will be proved.

For each rational point,
Table~\ref{Table:rational points} next gives 
the corresponding value of $j\colon C_i \to \PP^1$ 
calculated by Lemma~\ref{L:canonical morphism}.
If $j \notin \{1728,0,\infty\}$, 
it gives also
the Cremona label for an elliptic curve of that $j$-invariant,
and all primitive solutions 
$(a,b,c)$ to $a^2+b^3=c^7$ with $1728\,b^3/c^7 = j$.
If the conductor $N$ is out of the range of Cremona's database,
we use $N*$ as a substitute for the Cremona label.

\begin{remark}
The last two points on $C_6$ are interchanged by the
involution $z \mapsto -z$ of $C_6$,
and hence map to the same $j$-invariant.
\end{remark}

\section{Overview of the strategy for determining the rational points}

In order to determine the sets of rational points on $C_1, \ldots , C_{10}$,
we will need to determine the Mordell-Weil ranks of their corresponding
Jacobians $J_1,\ldots ,J_{10}$ over $\Q$. 
The curves $C_1$, $C_2$ and~$C_3$
are $\mu_7$-twists of~$X$ rather than more general twists. So there are more
methods available for dealing with them than for the other curves.

First observe that each of these three curves is a Galois cover of $\PP^1$,
over $\Q$, with
Galois group~$\mu_7$; in fact, $C_1$, $C_2$ and~$C_3$ are birational
to the (singular) projective plane curves
\[ u^7 = v^2 (6\,v - w) w^4\,, \qquad
   u^7 = v^2 (18\,v - w) w^4 \text{\quad and\quad}
   u^7 = v^2 (12\,v - w) w^4\,,
\]
respectively. To see this, rescale and/or permute the variables to obtain
an equation
of the form $x^3 y + y^3 z + A\,z^3 x = 0$ (with $A = 6, 18, 12$); then
$(u : v : w) = (xyz : -z^3 : x^2 y)$ does the transformation.

The curves (and their Jacobians) are therefore amenable to the techniques
developed by Schaefer~\cite{Schaefer1998} and
Poonen and Schaefer~\cite{Poonen-Schaefer1997}: the Jacobians have
complex multiplication by $\Z[\zeta]$. We will do a $1-\zeta$ descent
for each of these Jacobians.

For $C_4, \ldots , C_{10}$, we need to develop a new method of determining
the Mordell-Weil ranks of their Jacobians over $\Q$. 
We do this in Section~\ref{S:2-descent}.

For the curves $C_1, \ldots , C_{10}$ other than $C_5$, 
the rank of $J_i(\Q)$ is less than $3$, the dimension of $J_i$. In each case,
a combination of Chabauty and Mordell-Weil sieve computations
determines the set of rational points (see Section~\ref{S:Chabautysection}).
For $C_5$, the rank of $J_5(\Q)$ is $3$,
so we will deal with that curve separately 
in Section~\ref{S:The strategy for $C_5$}.

\section{($1-\zeta$)-descent on $J_1$, $J_2$, and $J_3$}
\label{S:1-zeta-descent}

We now do a descent using the endomorphism $1-\zeta$ of~$J_i$ 
over $L:=\Q(\zeta)$ for $i = 1, 2, 3$. 
See \cite{Schaefer1998} and \cite{Poonen-Schaefer1997} for details.
The descent map is given by the function $f = v/w$ and takes values in
\[ H = L(\{2,3,7\}, 7)
      = \{\theta\in L^\times/L^{\times 7}
           : L(\sqrt[7]{\theta})/L \text{\ is unramified outside~\{2,3,7\}}\}\,.
\]
The group $H$ is an $\F_7$-vector space of dimension~$7$.

For each prime $\pp$ of $L$, let $\rho_{\pp}$ be the composition
$H \injects L^\times/L^{\times 7} \to L_{\pp}^\times/L_{\pp}^{\times 7}$.
Let $\pi= 1-\zeta \in L$.
Let $\frakq$ be a prime of $L$ above~2.
For each $i$, we found that 
\[
\rho_{\pi}^{-1} f\left(\frac{J_i(L_{\pi})}{(1-\zeta)J_i(L_{\pi})}\right)
\cap 
\rho_{\frakq}^{-1} 
 f\left(\frac{J_i(L_{\frakq})}{(1-\zeta)J_i(L_{\frakq})}\right) 
= \langle 2, 3 \rangle \subset H;
\]
in particular the intersection is of dimension~$2$. 
In addition, we found two independent elements 
in $J_i(L)/(1-\zeta)J_i(L)$.
Therefore, the Selmer group $S^{1-\zeta}(L,J_i)$ 
and $J_i(L)/(1-\zeta)J_i(L)$ are of dimension~$2$. 
Since $\dim J_i(L)[1-\zeta]=1$, the rank of $J_i(\Q)$ is~$1$.

For $i=1, 2, 3$ and $\pp \in \{\pi, \frakq\}$, 
we found a basis of $J_i(L_{\pp})/(1-\zeta)J_i(L_{\pp})$ 
consisting of divisor classes 
of the form $[P-\infty]$ or $[P + \overline{P}-2\infty]$
where $\infty$
is the single point with $w=0$ and $P$ is defined over~$L_{\pp}$, or $P$ is 
defined over
a quadratic extension of~$L_{\pp}$ and $\overline{P}$ is its conjugate.
In Table~\ref{Table:1-zeta descent} 
we describe each such divisor class by giving the value of $v/w$ at $P$.
The last column of Table~\ref{Table:1-zeta descent} does the same
for $J_i(L)/(1-\zeta)J_i(L)$.

\begin{table}
\begin{center}
\renewcommand{\arraystretch}{1}
\begin{tabular}{c|c|c|c|}
$i$      & $\dfrac{J_i(L_{\pi})}{(1-\zeta)J_i(L_{\pi})}$
 &  $\dfrac{J_i(L_{\frakq})}{(1-\zeta)J_i(L_{\frakq})}$ & 
$\dfrac{J_i(L)}{(1-\zeta)J_i(L)}$
\\[3mm] \hline 
1 & 0 & $-1/2$ & 0 \\
 & $(3+6\pi^3+2\pi^4)/6^5$ & & $-1/2$ \\
 & $(3+\pi^3+5\pi^4+\pi^5)/6^5$ & & \\ 
 &  $\pi^7(1+5\pi^6+5\pi^7)/6^5$ & & \\ \hline
2 & $1/27$ & $1/27$ & 0 \\
 & $\pi^7(1+3\pi^7)/864$ & & $1/27$ \\
 & $(4+\pi^4)/864$ & & \\
 & $(4+\pi^4+\pi^5)/864$ & &  \\ \hline
3 & 0 & $1/8$ & 0 \\
 &  $\pi^7(1+\pi^6)/1944$ & & $1/8$ \\
 & $(3+\sqrt{3}\pi^3+\pi^4/\sqrt{3})/1944$ & & \\
 & $(2+\pi^{3/2} + \pi^{5/2} + 3\pi^3 + \pi^{7/2}\quad$ & & \\
 & $\quad{}- \pi^4 + \pi^{9/2} + 2\pi^5 +2\pi^6)/1944$ &  & \\ \hline
\end{tabular}
\end{center}
\caption{Data needed for $(1-\zeta)$-descents for $C_1$, $C_2$ and $C_3$.}
\hrule
\label{Table:1-zeta descent}
\end{table}

For each $i$, we have
$J_i[1-\zeta] \subseteq J_i(\Q)$ and $\dim J(\Q)[7]=1$,
so $\dim J_i(\Q)/7J_i(\Q) =2$.
We will use this in Section~\ref{S:Results}.

\subsection{Side application to twists of the degree-$7$ Fermat curve}

Our calculation of $C_1(\Q)$, $C_2(\Q)$, $C_3(\Q)$
will have consequences for rational points on certain curves 
$c_1 X^7+ c_2 Y^7+ c_3 Z^7=0$ in $\PP^2$.
These are not needed for the proof of our main theorem,
but are of interest in their own right.

A rational point on a curve
\[ u^7 = v^2 (Av - w) w^4 \,, \]
can be represented by $u,v,w \in \Z$ with $\gcd(v,w)=1$.
Then $\gcd(v,Av-w)=1$; therefore there exist integers $a,b,r,s,t$ such that
\[ v = r^7\,, \quad Av - w = as^7\,, \quad w = bt^7, \]
and $a b^4$ is a $7^{\text{th}}$ power, and $a$ and~$b$ 
are divisible only by primes dividing~$A$. 
Therefore, we obtain a point on the covering curve
\[ A\,r^7 -  a\,s^7 - b\,t^7 = 0 \,. \]
Conversely, a rational point on this curve maps down to a rational point
on the original curve.

Section~\ref{S:Chabautysection} determines the rational points
on $u^7 = v^2 (Av - w) w^4$ for $A = 6, 12, 18$.
Assuming this, we can list the rational points
on the covering curves arising from these.
After simplifying their equations (by multiplying by constants, and making
the coefficients $7^{\operatorname{th}}$-power free),
and eliminating those that over some $\Q_p$ fail to have points,
we obtain the following:

\begin{corollary} 
The following curves in $\PP^2_\Q$ have $\Q_p$-points for all $p \le \infty$,
but the only rational points are the ones with $X, Y, Z \in \{-1, 0, 1\}$:
  \begin{gather*}
    X^7 + Y^7 + 12\,Z^7 = 0\,,\quad 
    X^7 + Y^7 + 18\,Z^7 = 0\,,\quad 
    X^7 + Y^7 + 48\,Z^7 = 0\,,\\
    X^7 + Y^7 + 144\,Z^7 = 0\,,\quad
    X^7 + Y^7 + 162\,Z^7 = 0\,,\quad
    X^7 + Y^7 + 324\,Z^7 = 0\,,\\
    X^7 + 2\,Y^7 + 3\,Z^7 = 0\,,\quad
    X^7 + 2\,Y^7 + 81\,Z^7 = 0\,,\quad
    X^7 + 3\,Y^7 + 4\,Z^7 = 0\,,\\
    X^7 + 3\,Y^7 + 16\,Z^7 = 0\,,\quad
    X^7 + 4\,Y^7 + 9\,Z^7 = 0\,,\quad 
    X^7 + 9\,Y^7 + 16\,Z^7 = 0\,,\\
    X^7 + 16\,Y^7 + 81\,Z^7 = 0.
  \end{gather*}
\end{corollary}


\section{$2$-descent on Jacobians of twists of $X$}
\label{S:2-descent}

\subsection{Theory}
\label{S:2-descent theory}
To compute the ranks of $J_4(\Q), \ldots ,J_{10}(\Q)$,
we will use the $2$-descent described in this section.
It is similar to the descent described in \cite{Poonen-Schaefer1997}. 
 
The Klein quartic $X$ has $24$ Weierstrass points,
each of weight~$1$ (ordinary flexes).
These can be partitioned into eight sets of three with the following property.
Choose one set of three and denote the points $W_1, W_2, W_3$. The
tangent line at $W_i$ intersects $X$ with multiplicity~$3$ 
at $W_i$ and multiplicity~$1$ at $W_{i+1}$
(with subscripts considered modulo $3$). 
We will call such a set of three Weierstrass points a {\em triangle}.
The set of eight triangles is a Galois-stable set.

By abuse of notation,
let $T_i$ denote both a triangle and the effective degree-$3$
divisor that is the sum of its points.
For $1\leq i<j\leq 8$, we can
find a function with divisor 
$2T_i - 2T_j$ whose numerator and denominator are cubics. 
Therefore the divisor class $[T_i - T_j]$ is killed by 2.
Since the points in 
each $T_i$ are not collinear,
an application of the Riemann-Roch theorem shows that
$T_i-T_j$ is not principal. Thus
the order of $[T_i - T_j]$ is exactly two. 

\begin{proposition}
\label{triangle}
The $[T_i - T_1]$ for $i=1, \ldots ,8$ sum to $0$, and any six
of them with $i\neq 1$ form a basis for $\Jac(X)[2]$.
Let $Q_1 = [\sum_{i=1}^8 a_i T_i]$ and $Q_2 = [\sum_{i=1}^8 b_i T_i]$
with $\sum a_i = \sum b_i = 0$. Let $e_2$ be the 2-Weil pairing.
Then $e_2(Q_1,Q_2)= (-1)^{\sum_{i=1}^8 a_i b_i}$.
\end{proposition}

\begin{proof}
Over some finite and totally ramified extension $K$ of $\Q_7$, 
the Klein quartic $X$ acquires good reduction
and has all its automorphisms defined.
So if $O_K$ is the valuation ring of $K$, 
then $X$ extends to a smooth curve over $O_K$.
The special fiber $X_{\text{s}}$ is the hyperelliptic curve branched above
the eight points in $\PP^1(\F_7)$.  
Then $\Aut (X)$ injects into $\Aut (X_{\text{s}})$ and injects
into the automorphism group of $\PP^1_{\F_7}$,
which is isomorphic to $\PGL_2(\F_7)$ (see \cite{Elkies1999}).
At a Weierstrass point $W_i$ of $X$,
there is a regular differential vanishing to order 3
(it corresponds to the tangent line in the canonical
model).  One can scale the differential by an element of $K$ so that it extends
to a regular differential on the model $X$ over $O_K$, and reduces
to a nonzero differential on $X_{\text{s}}$, vanishing at least 
to order 3 at the reduction
of the point.  This means that Weierstrass points of $X$ 
reduce to Weierstrass points
of the hyperelliptic curve~$X_{\text{s}}$.

The differential in characteristic 0 has divisor $3W_i+W_{i+1}$, where $W_i, 
W_{i+1}$ are two points in a
triangle.  The reduction must be $3W_i'+W_{i+1}'$ where 
$W_i', W_{i+1}'$ are Weierstrass points.
But the only possibility for $W_{i+1}'$ 
that makes this the divisor of a differential
on the hyperelliptic curve is $W_{i+1}'=W_i'$.
Repeating this argument shows that the points of a triangle all reduce
to the same Weierstrass point on $X_{\text{s}}$.
Applying the automorphisms, and using the fact that $\PSL_2(\F_7)$
acts transitively on the set $\PP^1(\F_7)$ of hyperelliptic branch points
for $X_{\text{s}}$, 
shows that the eight triangles reduce 
to the eight Weierstrass points on $X_{\text{s}}$.
Now the $2$-torsion point $[T_i-T_j]$
reduces to the $2$-torsion point $[3P_i-3P_j]=[P_i-P_j]$
where $P_i$ (for example) is the reduction of any point in $T_i$.
The map ${\rm Jac} (X)[2]\to {\rm Jac} (X_{\text{s}})[2]$ 
is an isomorphism of groups preserving the Weil pairing,
because $X$ has good reduction.

The desired statements for $X$ and the $T_i$
now follow from their analogues
for the genus-$3$ hyperelliptic curve $X_{\text{s}}$
and its Weierstrass points,
which are well known
(cf.~Propositions 6.2 and~7.1 in~\cite{Poonen-Schaefer1997}).
\end{proof}

This correspondence between $X$ 
with its $T_i$
and a genus-$3$ hyperelliptic curve with its Weierstrass points will enable
us to translate many of the results in \cite{Poonen-Schaefer1997} to our situation. 
That article involves descent on cyclic covers of the projective line
using ramification points of the covering. For hyperelliptic curves, these
are the Weierstrass points. In this section, for brevity, we will sometimes
cite an equivalent proof in \cite{Poonen-Schaefer1997} which is easily translated. 
First we need to find a function corresponding to the $x-T$ map in 
\cite{Poonen-Schaefer1997}. 

Let $C$ be a twist of $X$ with a rational point $P\in C(\Q)$. In addition,
assume that for each triangle $T_i$, there is a cubic whose intersection 
divisor with $C$ is $2T_i + 3P + R_i$ 
where $R_i$ is an effective divisor of degree~$3$, 
defined over $\Q$ and supported on three non-collinear $\Qbar$-points. 
A calculation shows that this holds for $C_4, \ldots , C_{10}$.
Let $J$ be the Jacobian of $C$.

Choose functions $f_1, \ldots , f_8$ with 
$\divv (f_i) = 2T_i + 3P + R_i - 3\Omega$ (where $\Omega$ is a canonical 
divisor defined over $\Q$) and such that if $\sigma\in \GQ$, and 
$\ssigma\, T_i = T_j$ then $\ssigma f_i = f_j$.
Then $\divv (f_i/f_j) = 2T_i - 2T_j +R_i - R_j$. Since $2T_i - 2T_j$
is principal, so is $R_i - R_j$. The $R_i$ are effective, of degree
3, and are supported on three non-collinear points. The Riemann-Roch theorem
then implies that $R_i=R_j$. 
Let $R$ be the common value of the $R_i$.

Let $T = \{ T_i \}_{i=1}^8$. Let $K$ be $\Q$ or $\Q_p$
(we allow $p=\infty$, in which case $\Q_p:=\R$).
Let $\Abar_K$ be the algebra of maps from
$T$ to $\Kbar$.
The group $\GK$ acts on $\Abar_K$ through its action on $T$ and $\Kbar$.
Let $A_K$ be the algebra of $\GK$-invariants in $\Abar_K$;
it is isomorphic as $K$-algebra 
to a product of finite extensions of $K$ corresponding
to the $\GK$-orbits of $T$.
We call $D$ a {\em good divisor} if 
it is a divisor of degree~$0$, defined over $K$,
whose support is disjoint from the support of
$\divv(f_i)$ for every $i$.
If $D$ is a good divisor, then we define $f(D) = (T_i\mapsto f_i(D))$.
Then $f$ is a homomorphism from the group of good divisors to $A_K^{\times}$.
The divisor of each $f_i$ is the sum of a double and a $\Q$-rational
divisor. 
So by Weil reciprocity, 
$f$ sends good principal divisors to $A_K^{\times 2}K^{\times}$.
Since $C$ has a rational point, every element of $J(K)$ is represented
by a good divisor. Therefore, $f$ induces a well-defined homomorphism 
$f\colon J(K)/2J(K)\to A_K^{\times}/A_K^{\times 2}K^{\times}$. In order to determine
the kernel and image of this map, it will help to find a cohomological
interpretation of it.

Let $\mu_2^T$ be the $2$-torsion in $\Abar_K^\times$; 
it equals the $\GK$-module of maps from $T$ to $\mu_2$. 
Let $\mu_2^T/\mu_2$ be the quotient by the set of constant maps.
Let $q$ be the canonical surjection $\mu_2^T \to \mu_2^T/\mu_2$.
Define $\epsilon \colon J[2]\to \mu_2^T/\mu_2$ by 
$\epsilon(Q)= (T_i \mapsto e_2(Q,[T_i-T_1]))$.
Let $N\colon \mu_2^T \to \mu_2$ be the restriction of the norm map
$\Abar_K \to \Kbar$.
A straightforward computation involving Proposition~\ref{triangle}
shows that the following is a commutative
diagram of $\GK$-modules with exact rows and columns.
\[
\xymatrix{ & & \mu_2 \ar[d] \\
           & & \mu_2^T \ar[d]^-{q} \ar[r]^N & \mu_2 \ar@{=}[d] \ar[r] & 1 \\
           0 \ar[r] & J[2] \ar[r]^{\epsilon}
             & \dfrac{\mu_2^T}{\mu_2} \ar[r]^N & \mu_2 \ar[r] & 1
         }
\]

We now compute the long exact sequences of $\GK$-cohomology for the rows
and columns. 
Since $\HH^1(\GK,\Kbar^{\times})$ and $\HH^1(\GK,\Abar_K^{\times})$ are
both 0, we have
$\HH^1(\GK,\mu_2)\isom K^{\times}/K^{\times 2}$ and 
$\HH^1(\GK,\mu_2^T) \isom A_K^{\times}/A_K^{\times 2}$ by Kummer isomorphisms.
We obtain the following diagram with exact rows and columns.
\[
\xymatrix{ & J(K)/2J(K) \ar[rdd]^{\delta'} \ar[rrd]^-f 
                & & K^{\times}/K^{\times 2} \ar[d] \\
           & & 
               & A_K^{\times}/A_K^{\times 2} \ar[d]^{q'} \ar[r]^N
                 & K^{\times}/K^{\times 2} \ar@{=}[d] \\
           \HH^0(\GK,\frac{\mu_2^T}{\mu_2}) \ar[r]^-N
           & \mu_2 \ar[r]^-{\delta}
             & \HH^1(\GK,J[2]) \ar[r]^{\epsilon}
               & \HH^1(\GK,\frac{\mu_2^T}{\mu_2}) \ar[d] \ar[r]^N
                 & K^{\times}/K^{\times 2} \\
           & & & \Br(K)[2]
         }
\]
The map $q'$ is the composition of the Kummer isomorphism of 
$A_K^{\times}/A_K^{\times 2}$ 
with
$\HH^1(\GK,\mu_2^T)$ and
the map induced by $q$. 
The map $\delta'$ comes from taking cohomology of
\[ 0 \To J[2] \To J \stackrel{2}{\To} J \To 0. \]

\begin{proposition}
\label{fcohom}
The maps $\epsilon\circ 
\delta'$ and $q'\circ f$ are the same as maps from $J(K)/2J(K)$
to $\HH^1(\GK,\mu_2^T/\mu_2)$ (i.e., the diagram above commutes).
\end{proposition}

\begin{proof}
Pick an element of $J(K)$.
Since $C$ has a rational point, 
that element is represented by a good divisor $D$.
Pick $D' \in \Div C_\Kbar$ whose support is disjoint from the
supports of the $\divv (f_i)$ and such that $2[D']=[D]$. 
Then  $\delta'([D])$ 
is the class of the cocycle $\sigma \mapsto [\ssigma D' - D']$
in $\HH^1(\GK,J[2])$. 
So $\epsilon\circ \delta'([D])$ is the class of the cocycle
$\xi_{\sigma} := (T_i\mapsto e_2([\ssigma D' - D'],[T_i - T_1]))$
in $\HH^1(\GK,\mu_2^T/\mu_2)$.
Let $g$ be a function, defined over $K(D')$, whose divisor is $2D'-D$. 
Then $\divv(\ssigma g / g)=2(\ssigma D' - D')$.  We also have
$\divv(f_i/f_1)=2(T_i-T_1)$. 
By definition,
\[
e_2([\ssigma D' - D'],[T_i - T_1])
= \frac{(f_i/f_1)(\ssigma D' - D')}{(\ssigma g/g)(T_i-T_1)}.\]
Therefore, we have 
\[
\xi_{\sigma} = 
\left( T_i\mapsto 
\frac{(f_i/f_1)(\ssigma D' - D')}{(\ssigma g/g)(T_i-T_1)}\right) .\] 
For each $i$, the divisor $E = 2T_i -\divv (f_i)$ is the same. Let $l=\sqrt{g(E)}$.
By Weil reciprocity, the square of 
$\eta_\sigma := f_1(\ssigma D' - D')(\ssigma g/g) (-T_1)\, \ssigma l/l$
is $1$, so $\eta_\sigma\in \mu_2$.
Since $\xi_{\sigma}$ is in the quotient $\mu_2^T/\mu_2$,
we can multiply the expression for $\xi_{\sigma}$ 
by the constant 
map $(T_i \mapsto \eta_\sigma)$ to get
\[
\xi_{\sigma} = \left( T_i\mapsto 
\frac{f_i(\ssigma D'-D')\ssigma l/l}{(\ssigma g/g)(T_i)}\right) = 
\left( T_i\mapsto \ssigma \nu/\nu \right)
\]
where $\nu = f_i(D')l/g(T_i)$. Now
\[
\nu^2 = \frac{f_i(2D')}{g(2T_i-E)} = \frac{f_i(\divv (g) + D)}{g(\divv (f_i))}
= f_i(D)\]
by Weil reciprocity.
It follows that $\xi_{\sigma}$ is the image of $f(D)$ under $q'$, since
$q'$ is the
composition of the Kummer
map $A_K^{\times}/A_K^{\times 2}\to \HH^1(\GK,\mu_2^T)$ and $q$, the map
induced by the quotient map.
\end{proof}

\begin{proposition}
\label{P:delta delta}
We have $\delta'([R-3P]) = \delta(-1)$.
\end{proposition}

\begin{proof}
By definition of $R=R_1$,
$2T_1+3P+R$ is the intersection of $C$ with a cubic,
and hence is linearly equivalent to the intersection of $C$
with any other cubic, 
which we may take to be the union of the tangent lines at the
three points of $T_1$.
Thus $[2T_1+3P+R]=[4T_1]$, so $[R-3P]=[2T_1-6P]$.
Therefore $\delta'([R-3P])$ is the class including
the cocycle $(\sigma \mapsto [\ssigma\, T_1 - T_1])$ in $\HH^1(\GK,J[2])$. 
{}From Proposition~\ref{triangle},
we have $\epsilon ([\ssigma\, T_1 - T_1])
= \ssigma M - M\in \mu_2^T$ where $M$ is the map sending $T_1$ 
to 1 and $T_i$ to $-1$ for $i>1$. 
Also $N(M)=-1$. 
Hence, by definition, $\delta(-1)$ is the class 
of the cocycle $(\sigma \mapsto [\ssigma\, T_1 - T_1])$ in $\HH^1(\GK,J[2])$. 
\end{proof}

\begin{proposition}
\label{mapker}
The kernel of $f \colon J(K)/2J(K) \to A_K^{\times}/A_K^{\times 2}K^{\times}$ 
is generated by $[R-3P]$. 
This element is trivial in $J(K)/2J(K)$ if and only if the
$\GK$-set $T$ has an orbit of odd size.
\end{proposition}

\begin{proof}
The maps $\delta' \colon J(K)/2J(K) \to \HH^1(\GK,J[2])$ and
$q' \colon A_K^{\times}/A_K^{\times 2}K^{\times} \to \HH^1(\GK,\mu_2^T/\mu_2)$ 
are injective.  
By Propositions \ref{fcohom} and~\ref{P:delta delta},
\[
\ker f = (\delta')^{-1}(\ker \epsilon) = (\delta')^{-1}(\delta(\mu_2)) = (\delta')^{-1}(\delta'(\langle [R-3P]\rangle)) = \langle [R-3P] \rangle,
\]

The proof of the second half of the proposition 
is essentially identical to the proof of
\cite{Poonen-Schaefer1997}*{Lemma~11.2}.
\end{proof}

To simplify notation, we let $A=A_{\Q}$ and $A_p=A_{\Q_p}$.
In Section~\ref{S:rank of J4} 
we will show how to compute the $\GQ$-orbits of $T$;
it will turn out that for each of $C_4, \ldots , C_{10}$, 
there is no orbit of odd size.
Therefore, the kernel of the map $f$ from $J(\Q)/2J(\Q)$ 
to $A^{\times}/A^{\times 2}\Q^{\times}$ has size 2 and is generated by
$[R-3P]$.

Let $A\isom\prod A_i$ 
where the $A_i$ are number fields. Let $p$ be a prime of $\Q$.
Let $a$ be an element of $A^{\times}$ and $a_i$ its image in $A_i$.
We say that $a\in A^{\times}/A^{\times 2}$ is unramified at $p$ if for each $i$,
the field extension $A_i(\sqrt{a_i})/A_i$ is unramified at primes over $p$. 
Let $S$ be a set of places of $\Q$ including $2$, the infinite prime, 
and all primes at which $J$ has bad reduction
excluding the odd primes at which the Tamagawa number is odd.
Let $(A^{\times}/A^{\times 2}\Q^{\times})_S$
be the image in $A^{\times}/A^{\times 2}\Q^{\times}$
of the elements of $A^{\times}/A^{\times 2}$ that are unramified
outside of primes in $S$.

Denote by $H$ the kernel of the norm from
$(A^{\times}/A^{\times 2}\Q^{\times})_S$ to $\Q^{\times}/\Q^{\times 2}$.

\begin{proposition}
The image of $f \colon J(\Q) \to A^{\times}/A^{\times 2}\Q^{\times}$ 
is contained in $H$.
\end{proposition}

\begin{proof}
The divisor of $f$ is the sum of the double of a divisor 
and a rational divisor. This was the condition necessary
for the proofs of the essentially identical results in 
\cite{Poonen-Schaefer1997}*{Props.\ 12.1 \& 12.4}. 
The exclusion of primes with odd Tamagawa numbers is explained in
\cite{Schaefer-Stoll2004}*{Prop.~3.2}.
\end{proof}

Let $p$ be a prime of $\Q$, finite or infinite.
We have the following commutative diagram
\[
\xymatrix{ J(\Q)/\ker(f) \ar[d] \ar[r]^-{f} & H \ar[d]^{\rho_p} \\
           J(\Q_p)/\ker(f_p) \ar[r]^-{f_p}
             & A_p^{\times}/A_p^{\times 2}\Q_p^{\times}
         }
\]
where $f_p$ is the map $f$ for the case $K=\Q_p$,
and $\rho_p$ is the composition
$H \injects A^{\times}/A^{\times 2}\Q^{\times} \to A_p^{\times}/A_p^{\times 2}\Q_p^{\times}$.

Recall that the $2$-Selmer group $\Sel^2(J,\Q)$ is
the set of elements in $\HH^1(\GQ,J[2])$, unramified outside $S$, which map
to the image of $J(\Q_p) \to \HH^1(\GQp,J[2])$ for all $p\in S$.
Define the {\em fake $2$-Selmer group} to be
\[ \Sel_{\fake}^2(J,\Q)=\{ \theta\in H 
        \mid \rho_p(\theta) \in f(J(\Q_p)) \text{ for all $p\in S$} \} \,.
\] 

\begin{proposition}
\label{P:real and fake}
The sequence $\mu_2\stackrel{\delta}{\To} 
\Sel^2(J,\Q)\stackrel{\epsilon}{\To} \Sel_{\fake}^2(J,\Q) \to 0$ is exact.
\end{proposition}

\begin{proof}
See the proof of \cite{Poonen-Schaefer1997}*{Thm.~13.2}.
\end{proof}

For the curves $C_4, \ldots ,C_{10}$,
we find that $\delta'([R-3P])$ is nonzero.
By Propositions \ref{P:delta delta} and~\ref{P:real and fake},
$\dim_{F_2} \Sel^2(J,\Q) = 1 + \dim_{\F_2} \Sel_{\fake}^2(J,\Q)$.

To compute $\Sel_{\fake}^2(J,\Q)$, we first find $A$. Let $\Lambda$ be a
subset of $\{ 1, \ldots , 8\}$ such that the set $\{ T_j \}_{j\in \Lambda}$
contains one representative of each $\GQ$-orbit of $T$. 
Let $A_j = \Q(T_j)$. Then we can find an 
isomorphism 
\[
A\isom \prod\limits_{j\in \Lambda} A_j.\] 
The composition of $f$ and this isomorphism is $\prod_{j\in \Lambda} f_j$.

Next we find a basis of
$(A^{\times}/A^{\times 2}\Q^{\times})_S$; an algorithm for doing this
is given in \cite{Poonen-Schaefer1997}*{\S 12}. 
To find a function $f_j$, we 
find a cubic form defined over $A_j$, 
with the property that the curve defined by the cubic meets $C$ at $P$
with multiplicity (at least) 3 and at each of the three points of $T_j$ with
multiplicity (at least) 2. The $j$th component of~$f$ can then be taken
to be this cubic divided by any cubic form defined over~$\Q$ (for example,
$z^3$). 

Then we find a basis for each $J(\Q_p)/\ker(f)$. 
It helps to know the dimension in advance.
For $p$ odd, we have $\# J(\Q_p)/2J(\Q_p) = \# J(\Q_p)[2]$.
We also have $\# J(\Q_2)/2J(\Q_2) = 2^3\# J(\Q_2)[2]$ and 
$\# J(\R)/2J(\R) = \# J(\R)[2]/2^3$ 
(cf.~\cite{Poonen-Schaefer1997}*{Lemma~12.10}).
We can use Proposition~\ref{mapker}
to determine whether $J(\Q_p)/\ker(f)$ has the same or half the size of 
$J(\Q_p)/2J(\Q_p)$. Then we determine the intersection of all 
$\rho_p^{-1}(f(J(\Q_p)))$ for $p\in S$. That intersection is equal to 
$\Sel_{\fake}^2(J,\Q)$, which in our cases is half the size of $\Sel^2(J,\Q)$.

We have 
\[ \dim\Sel^2(J,\Q)=\rank J(\Q) + \dim J(\Q)[2]+ \dim\Sha (J,\Q)[2] \,. \]
In each of our cases, the number of independent elements
we found so far in $J(\Q)/2J(\Q)$ equals $\dim \Sel^2(J,\Q)$,
so $\Sha(J,\Q)[2]=0$. 
Subtracting $\dim J(\Q)[2]$ gives us $\rank J(\Q)$.

\subsection{Results}
\label{S:Results}

In this section we use the $2$-descent developed in the previous section
to determine the ranks of $J_4(\Q),\ldots,J_{10}(\Q)$.
In each case, we find that
$\rho_2^{-1}(f(J_i(\Q_2)/\ker f))$ 
is the same as the image of the subgroup of $J_i(\Q)$ 
generated by the known rational points. 
Therefore, in each case,
$\Sel_{\fake}^2(J_i,\Q) = \rho_2^{-1}(f(J_i(\Q_2)/\ker f))$.
Recall that in each case,
the $\GQ$-set $T$ has no orbit of odd size,
so the dimension of $\Sel^2(J_i,\Q)\isom J_i(\Q)/2J_i(\Q)$ 
equals $1 + \dim \Sel_{\fake}^2(J_i,\Q)$.

\subsubsection{Computing the rank of $J_4(\Q)$}
\label{S:rank of J4}

In this section, we give a detailed description of the computation of
the rank of $J_4(\Q)$. For the curves $C_5, \ldots , C_{10}$,
we give only the data necessary to check the computations in
Table~\ref{Table:2 descent}.

Let $g(x,y,z)=0$ be the model of  $C_4$ given in Table~\ref{Table:curve
 equations}. For any smooth plane curve, the flex points 
are the points on the curve where the Hessian vanishes.
We dehomogenize both $g$ and the Hessian with respect to $z$ and
get two polynomials in $u=x/z$ and $v=y/z$.
We compute their resultant to eliminate $v$, and get $h_6(u)h_{18}(u)$
where 
\[ h_6(u) = 7\,u^6 - u^3 + 1 \] and
\[ h_{18}(u) = 343\,u^{18} + 23667\,u^{15} + 127743\,u^{12} + 72128\,u^9
                - 29379\,u^6 + 2184\,u^3 + 1
\]
are irreducible.
This tells us that the $\GQ$-orbits of triangles will
be of size $2$ and $6$. To determine a polynomial giving us the
field of definition of a triangle in the size-$2$ orbit,
we begin by taking a zero $u_1$ of $h_6$,
and solve the equations given by $g$ and its Hessian
to find an explicit $v_1 \in \Q(u_1)$
such that $(u_1,v_1)$ is a flex;
let $T_1$ be the triangle it belongs to.
We then construct the tangent line to $C_4$ at this flex 
to find the next point in $T_1$
and repeat once more to find the third point of $T_1$,
expressing all coordinates as rational functions of $u_1$.
The sum of the three $u$-coordinates must lie in the (quadratic) field
of definition $\Q(T_1)$ of $T_1$, 
and this sum is found to generate $\Q(\sqrt{-3})$,
so $\Q(T_1)=\Q(\sqrt{-3})$.
A similar computation with a zero of $h_{18}$ 
shows that if $T_2$ is a triangle in the other orbit,
then $\Q(T_2)=\Q(\sqrt[6]{189})$.
So $A\isom \Q(\alpha)\times \Q(\beta)$ where $\alpha = \sqrt{-3}$ and
$\beta = \sqrt[6]{189}$. Let $S= \{ \infty, 2, 3, 7\}$.
The subgroup of $\Q(\alpha)^{\times}/\Q(\alpha)^{\times 2}$ 
that is unramified outside $S$ has basis $\{-1, 2, n_1, n_2, n_3\}$ 
where 
\[ n_1 = \alpha\,, \quad n_2 = \frac{\alpha + 5}{2}\,,\quad
   n_3 = \frac{-\alpha + 5}{2} \,. 
\]
Their norms are $1, 4, 3, 7, 7$,
respectively. 
The subgroup of $\Q(\beta)^{\times}/\Q(\beta)^{\times 2}$ that is unramified outside $S$ has basis
$\{ -1, \varepsilon_2, \varepsilon_3, \varepsilon_4, m_1, m_2, m_3, m_4, m_5\}$
where
\begin{gather*}
 \varepsilon_2
  = \frac{1}{54}(\beta^5 - 3\beta^4 + 3\beta^3 + 9\beta^2 - 45\beta + 27)\,,
 \quad
 \varepsilon_3
  = \frac{1}{54}(\beta^5 + 3\beta^4 + 3\beta^3 - 9\beta^2 - 45\beta - 27)\,,
 \\
 \varepsilon_4 
  = \frac{1}{6}(\beta^3 + 15)\,, \quad
 m_1 = \frac{1}{18}(-\beta^4 - 6\beta^2 + 9\beta - 36)\,, \\
 m_2 = \frac{1}{54}(-\beta^5 - 3\beta^3 + 9\beta^2 + 18\beta + 81)\,, \quad
 m_3 = \frac{1}{54}(\beta^5 + 3\beta^4 + 3\beta^3 + 27\beta^2 + 9\beta + 81)\,,
 \\
 m_4 = \frac{1}{54}(\beta^5 + 3\beta^3 + 9\beta^2 + 36\beta - 27)\,, \quad
 m_5 = \frac{1}{18}(\beta^4 - 9\beta)\,.
\end{gather*}
Their norms are $1, 1, 1, 1, 64, 64, 64, -48, 7$, respectively.

Now let us find the image of $\langle -1,2,3,7 \rangle \subseteq \Q^{\times}$ 
in $A^{\times}/A^{\times 2}$.
In $\Q(\alpha)$ we have 
\[ 2 = 2\,, \qquad 3 = -\alpha^2\,, \qquad 7 = n_2 n_3 \,. \]
In $\Q(\beta)$ we have 
\[ 2 = -m_1^{-1} m_2 m_3\,, \quad 3 = \varepsilon_4^{-1} m_1^{-4} m_4^6 
   \text{\quad and\quad} 7 = \varepsilon_4^3 m_5^6 \,.
\]
So in $A^{\times}/A^{\times 2}\Q^{\times}$ we have 
\[ -1 \equiv (-1, -1)\,, \quad 2 \equiv (2, -m_1 m_2 m_3)\,, \quad 
   3 \equiv (-1, \varepsilon_4)\,, \quad 7 \equiv (n_2 n_3, \varepsilon_4)
\] 
are all trivial. 
So $(A^{\times}/A^{\times 2}\Q^{\times})_S$
has basis 
\[
       \left\{ (-1,1), (n_1,1), (n_2,1), (1,\varepsilon_2), (1,\varepsilon_3),
               (1,m_1), (1,m_2), (1,m_3), (1,m_4), (1,m_5) \right\} \,.
\]
These have norms $1, 3, 7, 1, 1, 64, 64, 64, -48, 7$,
respectively.
So the kernel $H$ 
of the norm from this group to $\Q^{\times}/\Q^{\times 2}$ has basis
\[
   \left\{ (-1,1), (n_2,m_5), (1,\varepsilon_2), (1,\varepsilon_3), (1,m_1),
           (1,m_2), (1,m_3) \right\}.
\]
This group $H$ contains the fake 2-Selmer group.

Now we must find the function $f$. We will express $f$ as a pair
$(f_2, f_6)$, where $f_2$ is defined over $\Q(\sqrt{-3})$ and $f_6$ is defined
over $\Q(\sqrt[6]{189})$.
By linear algebra,
we can find a cubic $f_2(x,y,z)$, defined over the field $\Q(u_1)$, with the
property that the intersection divisor of $f_2 = 0$ with $g = 0$
includes $2T_1 + 3(0:1:1)$. 
We can scale the solution in such a way that the
coefficients lie in the field $\Q(\sqrt{-3})$. 
We get 
\begin{align*}
f_2(x,y,z) &= 238\,x^3 + 84\,x^2 y + (9\alpha  + 3) x y^2 + (12\alpha  + 4) y^3 
             + (21\alpha  - 21) x^2 z - (57\alpha + 33) x y z \\ 
         &\quad{} + (-12\alpha  + 24) y^2 z + (42\alpha  - 42) x z^2 
                  + (-6\alpha  + 12) y z^2 + (6\alpha  - 40) z^3 \,.
\end{align*}
Doing the same for $T_2$ yields 
\begin{align*}
f_6(x,y,z) &= 2142\,x^3 + (-51 \beta^4 + 1071 \beta + 756) x^2 y 
          + (4 \beta^5 + 42 \beta^4 - 21 \beta^3 - 882 \beta + 189) x^2 z \\
  &\quad{}+ (-8 \beta^5 - 9 \beta^3 + 189) x y^2 
          + (8 \beta^5 + 6 \beta^4 + 57 \beta^3 - 126 \beta - 1323) x y z \\
  &\quad{}+ (2 \beta^5 + 6 \beta^4 - 42 \beta^3 - 126 \beta + 378) x z^2
          + (-12 \beta^3 + 252) y^3 \\
  &\quad{}+ (2 \beta^5 - 9 \beta^4 + 12 \beta^3 + 189 \beta) y^2 z 
          + (-6 \beta^5 + 6 \beta^4 + 6 \beta^3 - 126 \beta) y z^2 \\
  &\quad{}+ (4 \beta^5 + 3 \beta^4 - 6 \beta^3 - 63 \beta - 252) z^3
                \,.
\end{align*}

To compute
$\Sel_{\fake}^2(J_4,\Q)$, we first do global computations.
Since the $\GQ$-orbits of the $8$ triangles have
sizes $2$ and $6$, we have $\dim J_4(\Q)[2]=1$. 
Since no orbit has odd size,
Proposition~\ref{mapker} gives
\[
   \dim J_4(\Q)/\ker(f) = \dim J_4(\Q)/2J_4(\Q) - 1 = \dim J_4(\Q)/2J_4(\Q) - \dim J_4(\Q)[2] = \rank J_4(\Q).
\]
The image in $H$, under $f$, 
of the subgroup of $J_4(\Q)$ generated by the three
known rational points is $\langle (-n_2,m_5), (1,\varepsilon_2)\rangle$.
So $\rank J_4(\Q) \ge 2$.

The eight triangles break into four ${\mathcal G}_{\Q_2}$-orbits 
of size $2$. Each pair
is conjugate over $\Q_2(\sqrt{-3})$. Therefore, we have $\dim J_4(\Q_2)[2]=4$
and $\dim J_4(\Q_2)/2J_4(\Q_2)=7$. 
Since there is no orbit of odd size,
Proposition~\ref{mapker} gives $\dim J_4(\Q_2)/\ker f=6$. 
We need to find a basis.
Let $A_2 = A\otimes_{\Q}\Q_2$. The map $f$ induces a map
\[ J_4(\Q_2) \To A_2^{\times}/A_2^{\times 2}\Q_2^{\times}
    \isom (\Q_2(\sqrt{-3})^{\times}/\Q_2(\sqrt{-3})^{\times 2})^4
      \bigm/(\Q_2^{\times}/\Q_2^{\times 2}) \,,
\]
To find a basis of $J_4(\Q_2)/\ker f$, we look for
$\Q_2$-rational divisors of degree 0 and determine their images under $f$.
It is a straightforward exercise in algebraic number theory to determine
whether their images are independent in 
$(\Q_2(\sqrt{-3})^{\times}/\Q_2(\sqrt{-3})^{\times 2})^4/
 (\Q_2^{\times}/\Q_2^{\times 2})$.
The following form a basis of $J_4(\Q_2)/\ker f$:
\begin{gather*}
 [(0:1:0) - (0:1:1)]\,, \quad 
 [(1:0:0) - (0:1:1)]\,, \\
 [(2:-1:O(2))-(0:1:1)]\,, \quad
 [(1:0:1+O(2))-(0:1:1)]\,, \\
 [(4:-3:O(2))-(0:1:1)]\,, \quad
 [(5:4:1+O(2))-(0:1:1)]\,. 
\end{gather*}
We compute the image of $J_4(\Q_2)/\ker f$ in
$(\Q_2(\sqrt{-3})^{\times}/\Q_2(\sqrt{-3})^{\times 2})^4/
 (\Q_2^{\times}/\Q_2^{\times 2})$,
and compute that its inverse image under $\rho_2|_H$
has basis $\{(-n_2,m_5), (1,\varepsilon_2)\}$.
Thus
\[
    2 \ge \dim \Sel_{\fake}^2(J_4,\Q) \ge \dim J_4(\Q)/\ker f = \rank J_4(\Q) 
      \ge 2,
\]
and equality holds everywhere.
(We did not need the conditions on $\Sel_{\fake}^2(J_4,\Q)$
coming from other primes.)

\subsubsection{Table with $2$-descent data}

We computed the fake 2-Selmer groups for $C_5, \ldots, C_{10}$ in ways
similar to that explained in detail for $C_4$. 
Table~\ref{Table:2 descent}
can be used to check the fake 2-Selmer group computations
for $C_5, \ldots , C_{10}$. 
In the column
labelled {\it octic}, we present the coefficients of $t^8, \ldots , 1$
for a polynomial defining the $8$-dimensional algebra $A$. 
We also give the dimensions of $J(\Q_2)[2]$ and $J(\Q)[2]$. 
Basis elements of
$J(\Q_2)/\ker f$ are represented by $\Q_2$-rational degree 0 divisors.
Lastly, we give the rank of $J(\Q)$.

Regarding $\dim J(\Q_2)[2]$, note that $2$-torsion points
are represented by a partition of the set of $8$ flex triangles into
two sets of even cardinality. So we have the trivial $2$-torsion point~$0$,
corresponding to the trivial factorization of the octic polynomial
over~$\Q_2$, points corresponding to factors of degree~$2$, and finally
points corresponding to a factorization of the octic polynomial into
two factors of degree~4 over~$\Q_2$ 
{\em or into two factors of degree~$4$ that are conjugate over some quadratic extension of~$\Q_2$.} 
For example, $\dim J_5(\Q_2)[2] = 2$, even though the
octic polynomial is irreducible over~$\Q_2$. 
This phenomenon
came up already when dealing with~$C_4$ above: 
the fact that the octic polynomial defining the algebra 
split into four quadratic factors over~$\Q_2$ 
defining the same quadratic extension $\Q_2(\sqrt{-3})$
gives rise to an additional independent point in $J_4(\Q_2)[2]$.

\begin{table}
\begin{center}
\renewcommand{\arraystretch}{1}
\begin{tabular}{c|c|c|c|c|c|}
 & & dim & & dim & \\
$i$ & octic & $J_i(\Q_2)[2]$ & basis of $J_i(\Q_2)/\ker(f)$ 
& $J_i(\Q)[2]$ & rank \\ \hline
5 & $1, 2, 0, 0, -21,$ & 2 & $[(1:1:1) - (0:1:0)]$ & 0 & 3 \\ 
 & $ 0, 84, -30, -123$ & & $[(1 : 1 : 5 + O(2^3)) - (0:1:0)]$ & & \\ 
 & & & $[(2 : 2 : 3 + O(2^3)) - (0:1:0)]$ & & \\
 & & & $[(-1 : 1 : 3 + O(2^3)) - (0:1:0)]$ & & \\ \hline
6 & $1, 0, -42, 0, 567,$ & 1 & $[(0:1:0) - (0:1:1)]$ & 1 & 2 \\
 & $0, -1890, 0, -567$  &  & $[(-1:1:0) - (0:1:1)]$ & & \\
 & & & $[(4:1:-1+O(2^2)) - (0:1:1)]$ & & \\ \hline
7 & $1, 2, -14, -98, -217,$ & 2 & $[(0:0:1) - (0:1:0)]$ & 0 & 2 \\ 
 & $-182, 196, 548, 529$ & & $[(2 : 2 : -1+O(2^3)) - (0:1:0)]$ & & \\
 & & & $[(2 : -2 : 1 + O(2^3)) - (0:1:0)]$ & & \\
 & & & $[(2 : -1 : 3 + O(2^3)) - (0:1:0)]$ & & \\ \hline
8 & $1, 0, 0, 56, -210,$  & 2 & $[(2:-1:0) - (0:0:1)]$ & 0 & 2 \\ 
 & $336, -224, 24, 21$  & & $[(0 : 2 : 1+O(2^3)) - (0:0:1)]$ & & \\
 & & & $[(-2 : 1 : 10 + O(2^4)) - (0:0:1)]$ & & \\
 & & & $[(1 : -3 : 1+O(2^3)) - (0:0:1)]$ & & \\ \hline
9 & $1, -2, -14, 14, 119,$  & 2 & $[(1:1:0) - (0:0:1)]$ & 0 & 2 \\
 & $182, 112, 16, -11$  & & $[(0 : 2 : -3+O(2^3)) - (0:0:1)]$ & & \\
 & & & $[(4 : 2 : 1+O(2^3)) - (0:0:1)]$ & & \\
 & & & $[(2 : 1 : 6 + O(2^4)) - (0:0:1)]$ & & \\ \hline
10 & $1, -3, 0, 14, 0,$  & 0 & $[(-1 : 1 : 2+O(2^3)) - (1:1:0)]$ & 0 & 2 \\ 
 & $-84, 112, -24, -12$  & & $[(0 : 1 : -1+O(2^3)) - (1:1:0)]$ & & \\ 
 & & & $[(1 : -4 : 4 + O(2^3)) - (1:1:0)]$ & & \\
\hline
\end{tabular}
\end{center}
\smallskip
\caption{Data needed for $2$-descents for $C_5, \ldots, C_{10}$.}
\hrule
\label{Table:2 descent}
\end{table}


\section{Chabauty and Mordell-Weil sieve}
\label{S:Chabautysection}

We discuss the Mordell-Weil sieve first 
since this is often performed before the Chabauty computation.


\subsection{Mordell-Weil sieve theory}
\label{S:Mordell-Weil sieve}

Let $C$ be a curve over $\Q$ and $J$ be its Jacobian.
(The assumption that the ground field is $\Q$ is for simplicity only;
other global fields would do just as well.)
Suppose that we know an embedding $C \to J$ over $\Q$.
What we call the {\em Mordell-Weil sieve}
is a method introduced 
by Scharaschkin~\cite{Scharaschkin2004preprint}
that uses reduction modulo $p$
for several primes $p$ to show that certain points in $J(\Q)$
cannot lie on the image of $C$.
Suppose also that we know explicit generators for $J(\Q)$.
Then the desired set $C(\Q)$ is the set of points of $J(\Q)$ that lie on $C$.
This set is difficult to compute: indeed, it is not known whether there
exists an algorithm to solve this problem in general.

Therefore we pick a prime $p$ of good reduction for $C$ (and $J$),
and we approximate the condition that a point $P \in J(\Q)$ lie on $C$
with the weaker condition that the reduction of $P$ in $J(\F_p)$
lie in $C(\F_p)$.
For each $p$, this weaker condition allows us to ``sieve out''
certain cosets of a finite-index subgroup in $J(\Q)$.
The conditions at different $p$ can interact nontrivially,
because the orders of $J(\F_p)$ may share prime factors.
If we are lucky, after using these sieve conditions at a few primes,
no points on $J(\Q)$ remain;
in this case one concludes that $C(\Q)$ is empty.
In fact, there is a heuristic that predicts that when $C(\Q)=\emptyset$,
sieving with a suitable set of $p$ will rule out all points 
of $J(\Q)$ \cite{Poonen-heuristic2005preprint}.
This obstruction to the existence of rational points
is closely related to the Brauer-Manin obstruction: 
see~\cite{Scharaschkin2004preprint} and in particular~\cite{Stoll2005preprint}.

\begin{remarks}\hfill
\begin{enumerate}
\item
With care, one can obtain sieve conditions also at primes $p$ 
of bad reduction.
Namely, let $\calC$ be the minimal proper regular model of $C$ over $\Z_p$,
and let $\calJ$ be the N\'eron model of $J$ over $\Z_p$.
Suppose we have an Albanese embedding $C \hookrightarrow J$
associated to a rational point on $C$.
Let $\calC^\smooth$ be the smooth locus of the structure morphism
$\calC \to \Spec \Z_p$; thus $\calC^\smooth$ is obtained 
from $\calC$ by deleting 
the singular points on the special fiber $\calC_{\F_p}$.
By the N\'eron property, the morphism $C \hookrightarrow J$
extends to a $\Z_p$-morphism $\calC^\smooth \hookrightarrow \calJ$.
Points in $C(\Q)$ extend to $\Z_p$-points of $\calC^\smooth$,
so we need only examine the set of points $J(\Q)$ whose
reduction in $\calJ(\F_p)$ lies in the image of $\calC^\smooth(\F_p)$.
\item
Sometimes it is more convenient to work in a homomorphic image $\Phi$
of $\calJ(\F_p)$ instead of $\calJ(\F_p)$ itself.
\item
Working in $\calJ(\Z/p^n\Z)$ for $n>1$ also might give more
information than $\calJ(\F_p)$ alone.
\item
We assumed that we know generators for $J(\Q)$,
but it may be enough to know generators for a subgroup $\calQ$
of finite index.
For instance, 
if one can prove that the index $(J(\Q):\calQ)$ is prime to
the order of $J(\F_p)$,
then $\calQ$ and $J(\Q)$ will have the same image in $J(\F_p)$.
\item
If $C(\Q)$ is known {\em a priori} to be nonempty,
obviously no complete obstruction to rational points can be obtained.
In this case, 
one can try to use the sieve information as input to a Chabauty argument,
to restrict the possibilities for residue classes.
The Chabauty argument still requires that the $p$-adic closure of $J(\Q)$
be provably of dimension less than $\dim J$, however,
and this usually requires $\rank J(\Q) < \dim J$.
\end{enumerate}
\end{remarks}

The first examples of this method were given 
in~\cite{Scharaschkin2004preprint}.
The method was refined and implemented 
for genus~$2$ curves over~$\Q$ in \cite{Flynn2004}.

\subsection{Chabauty theory}
\label{S:ChabautyTheory}

Here we explain Chabauty's method \cite{Chabauty1941},
as made explicit by Coleman \cite{Coleman1985chabauty}.
Let $C$ be a smooth projective curve.
Let $J$ be its Jacobian.
Assume that $C(\Q)$ is nonempty (this is true for all of our $C_i$);
choose $P \in C(\Q)$.
Identify $C$ with a subvariety of $J$ 
by mapping each $T \in C$ to $T \mapsto [T - P]$.
Assume that $J(\Q)$ has rank $r < g$, 
where $g = \dim J$ is the genus of~$C$.
Let $p$ be a finite prime.
If $\omega$ is a holomorphic $1$-form on $J_{\Q_p}$ 
or its pullback to $C_{\Q_p}$,
then following~\cite{Coleman1985torsion}
we define a homomorphism $\lambda_\omega\colon J(\Q_p) \to \Q_p$ 
by $Q \mapsto \int_0^Q\omega$. 
(This integral can be defined on a neighborhood of~$0$ using the formal group,
and then extended linearly to all of $J(\Q_p)$.)
Say that $\omega$ kills an element or subset $S$ of $J(\Q_p)$
if $\lambda_\omega|_S=0$.
By linear algebra,
we can find (to any desired $p$-adic precision) 
at least $g-r > 0$ independent forms $\omega$
killing $J(\Q)$ and hence $C(\Q)$.
The method consists of bounding the number of common zeros of
the corresponding integrals $\lambda_\omega$ on $C(\Q_p)$,
and hoping that enough points in $C(\Q)$ are found to meet the bound,
in which case these points are all of them.
If the method fails, we can try again with a different $p$,
or combine the argument with a Mordell-Weil sieve.

Let us explain how to do a Chabauty computation in practice.
Let $C$ be the smooth projective model of 
the affine plane quartic curve $g(u,v) = 0$.
For the computation, we need $G_1, \ldots, G_r \in J(\Q)$ which generate a 
subgroup of finite index in $J(\Q)$. 
Usually we obtain this information when we determine $r$.
Choose a (small) prime $p$ at which $C$, and hence $J$, has good reduction.
Pick a uniformizing parameter~$t$ of~$C$ 
at~$P$ that is also a uniformizer at~$P$ modulo~$p$.
For each $i$,
we calculate $m_i \in \Z_{>0}$ such that $m_i G_i$ reduces to~$0$ mod~$p$,
and find a divisor of the form $D_i - 3P$ linearly equivalent to
$m_i G_i$, where $D_i = P_{i,1} + P_{i,2} + P_{i,3}$ is an effective divisor
of degree~$3$, defined over~$\Q$. 
If $P$ does not reduce to a Weierstrass point mod~$p$, 
then it follows that each $P_{i,j}$ is congruent to~$P$ modulo~$p$.

A basis for the space of holomorphic differentials on~$C$ is given by 
\[ \omega_1 = \frac{dv}{g_u}\,, \quad 
   \omega_2 = \frac{u\,dv}{g_u} \quad\text{and}\quad
   \omega_3 = \frac{v\,dv}{g_u}\,. 
\]
We express $\omega_j$ as an element of $\Q[\!\,[t]\,\!]\,dt$ and integrate
formally to obtain $\lambda_j = \int_0^t \omega_j \in t\,\Q[\!\,[t]\,\!]$.
Then for all $i$ and~$j$, we compute
\begin{align*}
  \lambda_{\omega_j}(m_i G_i) 
    &= \int_0^{m_i G_i} \iota_{\ast} \omega_j
     = \int_0^{[D_i-3P]} \iota_{\ast} \omega_j
     = \int_{P}^{P_{i,1}} \omega_j + \int_{P}^{P_{i,2}} \omega_j
        + \int_{P}^{P_{i,3}} \omega_j \\[2mm]
    &= \lambda_j(t(P_{i,1})) + \lambda_j(t(P_{i,2})) + \lambda_j(t(P_{i,3}))\,.
\end{align*}
The series converge, since the $t(P_{i,k})$ are of positive $p$-adic
valuation (this is because $P_{i,k}$ has the same reduction mod~$p$ as~$P$).
Computing the kernel of 
$\omega \longmapsto (\lambda_\omega(m_1 G_1), \dots, \lambda_\omega(m_r G_r))$,
we find $s = 3-r$ independent holomorphic differentials 
$\eta_1, \ldots ,\eta_s$
such that the corresponding $1$-forms on $J_{\Q_p}$ kill $J(\Q)$.
We rescale each so that it reduces to a nonzero differential 
$\tilde\eta_i$ modulo~$p$.
Let $Q \in C(\F_p)$ and let $\nu = \min_i v_Q(\tilde\eta_i)$. 
If $p > \nu+1$, there can be at most $\nu+1$ rational points on~$C$
reducing to~$Q$ mod~$p$: see, e.g., \cite{Stoll2004preprint}.

For $i = 1, 3, 7, 8, 9, 10$, 
we want to do Chabauty arguments using the prime~$7$,
at which $C_i$ has bad reduction. 
However, the given model of~$C_i$ in $\PP^2_{\Z_7}$ is 
a minimal proper regular model over~$\Z_7$, 
and the special fiber $\tilde{C}_i$ 
has just one component, which is of (geometric) genus~0.
Let $\calJ_i$ be the N\'{e}ron model of~$J_i$ over~$\Z_7$. 
Let $\tilde{J}_i = \calJ_i \times_{\Z_7} \F_7$ be the special fiber 
of~$\calJ_i$. Since $\tilde{C}_i$ has just one component, $\tilde{J}_i$
is connected.
Let $\tilde{C}_i^\smooth$ be the smooth part of $\tilde{C}_i$.
We write $J_i(\F_7)$ for $\tilde{J}_i(\F_7)$,
and write $C_i(\F_7)$ for $\tilde{C}_i^\smooth(\F_7)$, 
which is a subset of $J_i(\F_7)$.
In all the cases we consider, we find that 
$J_i(\F_7) \isom (\Z/7\Z)^3$. 
We can proceed as in the case of
good reduction; the differentials reduce to holomorphic differentials
on $\tilde{C}_i^\smooth$ in this case.

\subsection{Results}

For brevity, we give full detail for the Chabauty and
Mordell-Weil sieve arguments for just one curve. 
We chose $C_2$ since our argument for $C_2$ involves a shortcut
that apparently is not described elsewhere in the literature.
Afterwards we sketch the Chabauty arguments for the remaining curves.

Define the {\em known part} $J_i(\Q)_\known$ of $J_i(\Q)$
as the subgroup generated by the rational
points in Table~\ref{Table:rational points} 
and the element $[R-3P]$ mentioned in
Propositions \ref{P:delta delta} and~\ref{mapker}.

\subsubsection{Determining $C_2(\Q)$}

We will show that 
\[ C_2(\Q) = \{(1: 0: 0), (0: 1: 0), (0: 0: 1), (1: -2: -1), (1: 1: -1)\}\,. \]
We work with $p = 5$; we have $J_2(\F_5) \isom \Z/126\Z$.
The reduction of $J_2(\Q)_\known$ is the cyclic 
subgroup of order~$63$. We use $P = (1: -2: -1)$  
as a basepoint. The five known rational points of $C_2$ reduce to
different points in $C_2(\F_5)$, which has size $6$.

We first have to find the space of $5$-adic differentials that vanish
on $J_2(\Q)$. Dehomogenize by setting $z = 1$; then a basis
of the space of holomorphic differentials on~$C_2$ is given by
\[ \omega_1 = \frac{dy}{9\,x^2 y + 2}\,, \quad
   \omega_2 = \frac{x\,dy}{9\,x^2 y + 2}\,, \quad
   \omega_3 = \frac{y\,dy}{9\,x^2 y + 2}\,.
\]
(As line sections, we have $\omega_1 \leftrightarrow z$,
$\omega_2 \leftrightarrow x$, $\omega_3 \leftrightarrow y$.
This point of view is useful when determining the zeros of differentials.)
The point $G = [(1:0:0) - (1:-2:-1)] \in J_2(\Q)$ is of infinite 
order, and $9G$ is in the kernel of reduction mod~$5$. We find an effective
divisor $D$ of degree~3 such that $9G = [D - 3P]$. 
We choose $t = x+1$
as a uniformizer at~$P$; then we can express $x$, $y$, and the differentials
in terms of~$t$ as follows.
\begin{align*}
  x &= -1 + t \\
  y &= 2 - \frac{20}{3^2}\,t + \frac{226}{3^5}\,t^2 + \frac{1024}{3^8}\,t^3
         - \frac{7702}{3^{10}}\,t^4 + \dots \\
  \omega_1 &= \Bigl(-\frac{1}{3^2} - \frac{53}{3^5}\,t - \frac{478}{3^7}\,t^2
                     - \frac{6530}{3^{10}}\,t^3 - \frac{1135}{3^{14}}\,t^4 
                     + \dots\Bigr)\,dt \\
  \omega_2 &= \Bigl(\frac{1}{3^2} + \frac{26}{3^5}\,t + \frac{1}{3^7}\,t^2
                     - \frac{6376}{3^{10}}\,t^3 - \frac{527795}{3^{14}}\,t^4 
                     + \dots\Bigr)\,dt \\
  \omega_3 &= \Bigl(-\frac{2}{3^2} - \frac{46}{3^5}\,t - \frac{122}{3^7}\,t^2 
                     + \frac{2618}{3^{10}}\,t^3 + \frac{107380}{3^{14}}\,t^4
                     + \dots\Bigr)\,dt
\end{align*}
We obtain the logarithms $\lambda_{\omega_i}$ in terms of~$t$
by formal integration of the latter three power series:
\begin{align*}
  \lambda_1 &= -\frac{1}{3^2}\,t - \frac{53}{2 \cdot 3^5}\,t^2 
                - \frac{478}{3^8}\,t^3 - \frac{3265}{2 \cdot 3^{10}}\,t^4 
                - \frac{227}{3^{14}}\,t^5 + \dots \\
  \lambda_2 &= \frac{1}{3^2}\,t + \frac{13}{3^5}\,t^2 + \frac{1}{3^8}\,t^3
                - \frac{1594}{3^{10}}\,t^4 - \frac{105559}{3^{14}}\,t^5
                + \dots \\
  \lambda_3 &= -\frac{2}{3^2}\,t - \frac{23}{3^5}\,t^2 - \frac{122}{3^8}\,t^3 
                     + \frac{1309}{2 \cdot 3^{10}}\,t^4 
                     + \frac{21476}{3^{14}}\,t^5 + \dots
\end{align*}
If $t_1, t_2, t_3$ are the values of~$t$ at the three points in~$D$, then 
\[ \lambda_{\omega_i}(9G)
     = \lambda_i(t_1) + \lambda_i(t_2) + \lambda_i(t_3)
     = \sum_{j=1}^\infty \lambda_{i,j} (t_1^j + t_2^j + t_3^j) \,,
\]
where $\lambda_i(t) = \sum_{j \ge 1} \lambda_{i,j} t^j$.
We obtain the power sums $t_1^j + t_2^j + t_3^j$ from the coefficients of
the characteristic polynomial $(X - t_1)(X - t_2)(X - t_3) \in \Q[X]$.
The $5$-adic valuation of each $t_i$ is positive (at least $1/3$), and so
the series above converge in~$\Q_5$. We obtain
\[ \lambda_{\omega_1}(9G) = -2 \cdot 5^3 + O(5^4)\,, \quad
   \lambda_{\omega_2}(9G) = -16 \cdot 5 + O(5^4)\,, \quad
   \lambda_{\omega_3}(9G) = -27 \cdot 5 + O(5^4)\,.
\]
We see that 
\[ \tilde{\eta_1} = \tilde{\omega_1} = \frac{dy}{2 - x^2 y}
                  \leftrightarrow z
   \quad\text{and}\quad
   \tilde{\eta_2} = \tilde{\omega_2} + 2\tilde{\omega_3}
                  = \frac{(x + 2y)\,dy}{2 - x^2 y}
                  \leftrightarrow x + 2y
\]
are reductions mod~$5$ of differentials killing $J_2(\Q)$.
Since the common zero of $z$ and $x+2y$ in $\PP^2(\F_5)$ is not on $C$,
we have $\min\{v_Q(\tilde{\eta_1}),v_Q(\tilde{\eta_2})\}=0$ at
all $Q \in C(\F_5)$.
Thus each of the six residue classes contains at most one rational point.
Since five of the six do contain a rational point,
it remains to show that there is no rational point reducing
to the sixth, which is $(2 : 2 : 1) \in C_2(\F_5)$.

To that end, we pick $Q \in C_2(\Q_5)$ reducing to $(2 : 2 : 1)$; we
take $Q = (\frac{3 + \sqrt{-39}}{6} : 2 : 1)$, where we choose the square root
that is $-1$ mod~$5$. We first need to find $\lambda_\omega([Q - P])$
for a basis of differentials~$\omega$ killing $J_2(\Q)$.
In principle, we can find this by first computing $21[Q - P] = [D' - 3P]$
in the kernel of reduction as above, but in this case, there is a trick,
taught to us by Joseph Wetherell, that allows us to simplify the computation.

Note that $[(2 : 2 : 1) - \tilde{P}] \in J_2(\F_5)$ is in the reduction
of $J_2(\Q)_\known$; we find that 
$G' = [(0:1:0) + (1:1:-1) - (1:0:0) - P] \in J_2(\Q)$ has this image.
If $\omega$ kills $J_2(\Q)$, then 
$\lambda_\omega([Q-P]) = \lambda_\omega([Q-P] - G')$, where now
$[Q-P] - G'$ is in the kernel of reduction, and we can compute the
logarithm as before. With
\[ \eta_1 = \omega_1 - 50\,\omega_2 + O(5^3) \quad\text{and}\quad
   \eta_2 = \omega_2 - 33\,\omega_3 + O(5^3) \,,
\]
we obtain the values $\lambda_{\eta_1}([Q-P]) = -12 \cdot 5 + O(5^3)$
and $\lambda_{\eta_2}([Q-P]) = -14 \cdot 5 + O(5^3)$. We now compute
$\lambda_{\eta_i}([Q' - Q])$ for points $Q'$ in the residue class of~$Q$.
Taking $T = y - 2$ as a uniformizer at~$Q$, we obtain for points
$Q'$ such that $T(Q') = 5\tau$:
\begin{align*}
  \lambda_{\eta_1}([Q' - P])
     = \lambda_{\eta_1}([Q' - Q]) + \lambda_{\eta_1}([Q - P])
    &= -12 \cdot 5 + 4 \cdot 5\,\tau + 5^2\,\tau^2 + O(5^3) \\
  \lambda_{\eta_2}([Q' - P])
     = \lambda_{\eta_2}([Q' - Q]) + \lambda_{\eta_2}([Q - P])
    &= 11 \cdot 5 + 9 \cdot 5\,\tau - 3 \cdot 5^2\,\tau^2 + O(5^3)
\end{align*}
We see that
$\lambda_{\eta_1}([Q' - P]) - \lambda_{\eta_2}([Q' - P])
    = 2 \cdot 5 + O(5^2)$
and so never vanishes. This proves that there is no rational point
on~$C_2$ reducing to $(2 : 2 : 1)$.

\subsubsection{Determining $C_1(\Q)$}

We show that 
\[ C_1(\Q) = \{(0:0:1), (0:1:0), (0:0:1), (1: -1: 2)\}\,. \]
As described in Section~\ref{S:ChabautyTheory}, we
have $J_1(\F_7) \isom (\Z/7\Z)^3$. 
{}From Section~\ref{S:1-zeta-descent},
we see that $\dim_{\F_7} J_1(\Q)/7J_1(\Q) = 2$.
The reduction of $J_1(\Q)_\known$ has dimension~$2$,
so that is also the reduction of the whole $J_1(\Q)$.

Using $(1: -1: 2)$ as a basepoint, we find that there are exactly
four points of $C_1(\F_7)$ in the reduction of $J(\Q)$.
They are the reductions of the 
four known rational points on~$C_1$. So it suffices to show that there is
only one point of $C_1(\Q)$ in each of those residue classes.

Proceeding as for~$C_2$, 
we find reductions of differentials killing $J_1(\Q)$:
\[ \tilde{\eta_1} \leftrightarrow x + 2y\,, \qquad
   \tilde{\eta_2} \leftrightarrow x - z\,.
\]
The common zero of $x+2y$ and $z$ in $\PP^2(\F_7)$
is not on $C_1$, so there is at most one rational point per
residue class, as desired.

\subsubsection{Determining $C_3(\Q)$}

We show that 
\[ C_3(\Q) = \{(1: 0: 0), (0: 1: 0), (0: 0: 1), (1: 1: -1)\}\,. \]
We have $J_3(\F_7) \isom (\Z/7\Z)^3$. From
Section~\ref{S:1-zeta-descent}
we see that $\dim_{\F_7} J_3(\Q)/7J_3(\Q) = 2$.
The reduction of $J_3(\Q)_\known$ has dimension~$2$,
so that is also the reduction of the whole $J_3(\Q)$.

Using $(1: 1: -1)$ as a basepoint, we find that there are exactly
four points of $C_3(\F_7)$ in the reduction of $J_3(\Q)$.
Again, they are the reductions of the 
four known rational points on~$C_3$. So it suffices to show that there is
only one point of $C_3(\Q)$ in each of those residue classes.

Proceeding as before, we find reductions of differentials killing $J_3(\Q)$:
\[ \tilde{\eta_1} \leftrightarrow 2x - y\,, \qquad
   \tilde{\eta_2} \leftrightarrow x + 2z\,.
\]
The common zero of $2x-y$ and $x+2z$ in $\PP^2(\F_7)$ is not on $C_3$, 
so we are done.

\subsubsection{Determining $C_4(\Q)$}

We show that 
\[ C_4(\Q) = \{(1: 0: 0), (0: 1: 0), (0: 1: 1)\}\,. \]
We have $J_4(\F_5) \isom (\Z/6\Z)^3$.
The reduction of $J_4(\Q)_\known$ is isomorphic
to $(\Z/6\Z)^2$. Since $J_4(\Q)$ and $J_4(\Q)_\known$
are both isomorphic to $\Z^2 \times \Z/4\Z$, and $J_4(\Q)_\known$
has odd index in $J_4(\Q)$, the whole $J_4(\Q)$
must also reduce to the same $(\Z/6\Z)^2$ subgroup. We use $(1:0:0)$ as
a basepoint. Then the only elements in $C_4(\F_5)$ in the reduction 
of $J_4(\Q)$ are the reductions of the three known rational points.
So it suffices to show that there is
only one point of $C_4(\Q)$ in each of those residue classes.

Since the rank is~$2$, we have only one differential, 
$\tilde{\eta} \leftrightarrow x + y$. It does not vanish at any of the
three relevant points of $C_4(\F_5)$; this proves the claim.

\subsubsection{Determining $C_6(\Q)$}

We show that
\[ C_6(\Q) = \{(0:1:0), (1:-1:0), (0:1:1), (0:1:-1)\}\,. \]
We use $P = (0:1:-1)$ as a basepoint.
We have an isomorphism $\Z^2 \times \Z/4\Z \to J_6(\Q)_\known$
taking the standard generators $(1,0,0)$, $(0,1,0)$, $(0,0,1)$
to $[R-3P]$, $[(0:1:0)-P]$, $[(1:-1:0)-P]$, respectively.
{}From the $2$-descent, we know that $J_6(\Q)_\known$
has odd index in $J_6(\Q)$, which also is isomorphic
to $\Z^2 \times \Z/4\Z$.
We have $J_6(\F_{11}) \isom (\Z/16\Z) \times (\Z/8\Z)^2 \times (\Z/2\Z)$.
The reduction of $J_6(\Q)_\known$ is isomorphic
to $(\Z/8\Z)^2 \times (\Z/4\Z)$,
with $[R-3P]$, $[(0:1:0)-P]$, $[(1:-1:0)-P]$
mapping to the standard generators.
By the above, $J_6(\Q)$ has the same reduction as $J_6(\Q)_\known$.
There are five points in $C_6(\F_{11})$ lying in the reduction 
of $J_6(\Q)$; four of those are the reductions 
of the four known rational points. 
The fifth is the reduction of 
\[ m\,[R - 3P] + n\,[(0:1:0) - P] + p\,[(1:-1:0) - P] \]
with $(m,n,p)=(1,4,2)$.

We have $J_6(\F_{23}) \isom \Z/32\Z \times (\Z/16\Z)^2 \times \Z/2\Z$.
The mod $23$ reduction of $J_6(\Q)$ 
is isomorphic to $(\Z/16\Z)^2 \times \Z/4\Z$,
with $[R-3P]$, $[(0:1:0)-P]$, $[(1:-1:0)-P]$
mapping to the standard generators.
We find that if $(m,n,p) \in \Z^3$ 
satisfy $(m \bmod 8, n \bmod 8,p \bmod 4)=(1,4,2)$
then the mod $23$ reduction of
\[ m\,[R - 3P] + n\,[(0:1:0) - P] + p\,[(1:-1:0) - P] \]
is not in $C_6(\F_{23})$.
Because $J_6(\Q)_\known$ has odd index in $J_6(\Q)$,
the points of $J_6(\Q)$ reducing
to the fifth point of $C_6(\F_{11})$ of the previous paragraph
do not reduce mod $23$ to points in $C_6(\F_{23})$.
So it suffices to show that there is
only one point of $C_6(\Q)$ in each of the four residue classes
of the four known rational points modulo $11$.

The differential is $\tilde{\eta} \leftrightarrow x + 5y$; it does
not vanish at the four relevant points of $C_6(\F_{11})$. This proves
the claim.

\subsubsection{Determining $C_7(\Q)$}

We show that
\[ C_7(\Q) = \{(0: 1: 0), (0: 0: 1), (0: 1: 1)\}\,. \]
We have $J_7(\F_7) \isom (\Z/7\Z)^3$.
The group $J_7(\Q)$ has rank~$2$ and no $7$-torsion.
The reduction of $J_7(\Q)_\known$ is isomorphic to $(\Z/7\Z)^2$,
so that is also the reduction of the whole $J_7(\Q)$.

We use $(0:1:0)$ as a basepoint.
Then the only elements in $C_7(\F_7)$ in the reduction 
of $J_7(\Q)$ are the reductions of the three known rational points.
So it suffices to show that there is
only one point of $C_7(\Q)$ in each of those residue classes.

Our differential mod~$7$ is $\tilde{\eta} \leftrightarrow x-2z$.
It does not vanish at $(0 : 0 : 1)$ and at $(0 : 1 : 1)$, but it
has a simple zero at $(0 : 1 : 0)$. So we see that there is at most
one rational point in each of the two first residue classes, but there
may be two in the third. Since we are using $(0:1:0)$ as a basepoint,
the information over~$\Q_7$ tells us that any point of $C_7(\Q)$ in the
residue class of $(0:1:0)$ is in $7 J_7(\Q)$. 

We have $J_7(\F_{13}) \isom (\Z/14\Z)^3$. 
The reduction of $J_7(\Q)_\known$ is isomorphic to $\Z/14\Z \times \Z/7\Z$.
{}From the $2$-descent, we know that the index of $J_7(\Q)_\known$ 
in $J_7(\Q)$ is odd. Thus, the
reduction of the whole $J_7(\Q)$ is isomorphic to $\Z/14\Z \times \Z/7\Z$. 
Only one point of $C_7(\F_{13})$ has trivial image in
$J(\F_{13}) \otimes \Z/7\Z$, namely $(0:1:0)$. 
A holomorphic differential killing $J_7(\Q)$
reduces mod~$13$ to $\tilde{\omega} \leftrightarrow x-2y$, which does not
vanish at $(0 : 1 : 0)$. So there is only one point of $C_7(\Q)$ in this 
residue class. (Incidentally, $\tilde{\omega}$ vanishes at one of the other
interesting residue classes, namely $(0 : 0 : 1)$. So information obtained
from~$\Q_{13}$ alone would also not be sufficient.)

\subsubsection{Determining $C_8(\Q)$}

We show that
\[ C_8(\Q) = \{(0: 0: 1), (2: -1: 0)\}\,. \]
We have $J_8(\F_7) \isom (\Z/7\Z)^3$.
The group $J_8(\Q)$ has rank~2 and no 7-torsion.
The reduction of $J_8(\Q)_\known$ is isomorphic to $(\Z/7\Z)^2$,
so that is also the reduction of the whole $J_8(\Q)$.

We use $(0:0:1)$ as
a basepoint. Then the only elements in $C_8(\F_7)$ in the reduction 
of $J_8(\Q)$ are the reductions of the two known rational points.
So it suffices to show that there is
only one point of $C_8(\Q)$ in each of those residue classes.

Our differential mod~$7$ is $\tilde{\eta} \leftrightarrow x - y - 3z$.
It does not vanish at $(0 : 0 : 1)$ and at $(2 : -1 : 0)$.
This proves the claim.

\subsubsection{Determining $C_9(\Q)$}

We show that
\[ C_9(\Q) = \{(0 : 0 : 1), (1 : 1 : 0)\}\,. \]
We use $P=(0 : 0 : 1)$ as a basepoint. 
We have $J_9(\Q)_\known \isom \Z^2$, 
generated by $[(1:1:0)-P]$ and $[R-3P]$.
The group $J_{9}(\Q)$ has rank~$2$ and no $7$-torsion.

We have $J_9(\F_7) \isom (\Z/7\Z)^3$,
with $[(1:1:0)-P]$ and $[R-3P]$ reducing to $\F_7$-independent points,
so their reductions form a basis also for the reduction of $J_9(\Q)$.

There are four points of $C_9(\F_7)$ in the reduction of
$J_9(\Q)$: the reductions of 
$m [(1:1:0)-P] + n [R-3P]$
with $(m,n) = (0,0), (1,0), (2,1), (4,3)$.
Clearly $(m,n) = (0,0)$ and $(1,0)$
correspond to the reductions of the two known rational points. 
So we must
eliminate the possibilities of $(m,n) = (2,1)$ and $(4,3)$.
We do this by working over $\F_{13}$.
We have $J_9(\F_{13}) \isom \Z/7\Z \times \Z/170\Z$. 
The classes $(m,n) = (2,1)$ and $(4,3)$ in $(\Z/7\Z)^2$
map to classes in $J_9(\F_{13}) \otimes \Z/7\Z$. 
that are not in the image of $C_9(\F_{13})$ in $J_9(\F_{13}) \otimes \Z/7\Z$.

Our differential mod~$7$ is $\tilde{\eta} \leftrightarrow x - y - 2z$.
It does not vanish at $(0 : 0 : 1)$, so there is only one
rational point in this residue class.
But $\tilde{\eta}$ has a simple zero at the other point, $(1 : 1 : 0)$.
So we have to do another Chabauty computation 
to show that $(1:1:0)$ is the only point of $C_9(\Q)$ 
in $[(1:1:0)-P] + 7 J_9(\Q)$.

We will use $p = 11$.
We have $J_9(\F_{11}) \isom \Z/10\Z \times \Z/140\Z$. 
We know from the $2$-descent
that the index of $J_9(\Q)_\known$ in $J_9(\Q)$ 
is odd. The reduction modulo~$11$ of $J_9(\Q)_\known$ is
isomorphic to $\Z/5\Z \times \Z/70\Z$ and so must be the reduction of 
the whole $J_9(\Q)$. Only three points of $C_9(\F_{11})$ are
in the reduction of $J_9(\Q)$. Only one of them the same image in
$J_9(\F_{11}) \otimes \Z/7\Z$ as $[(1:1:0)-P] + 7 J_9(\Q)$,
and that is the reduction of $(1:1:0)$. 
So it suffices to show that there is only one point
of $C_9(\Q)$ with reduction equal to $(1:1:0)$ modulo~$11$.
The differential mod~$11$ is $\tilde{\omega} \leftrightarrow x - 5y + 4z$,
and it does not vanish at $(1 : 1 : 0)$.
So there is only one point of $C_9(\Q)$ in this mod~$11$ residue class. 

\subsubsection{Determining $C_{10}(\Q)$}

We show that
\[ C_{10}(\Q) = \{(1 : 0 : 0), (1 : 1 : 0)\}\,. \]
We have $J_{10}(\F_7) \isom (\Z/7\Z)^3$.
The group $J_{10}(\Q)$ has rank~2 and no 7-torsion.
The reduction of $J_{10}(\Q)_\known$ is isomorphic to $(\Z/7\Z)^2$,
so that is also the reduction of the whole $J_{10}(\Q)$.
We use $(1 : 0 : 0)$ as a basepoint. Then the only elements in $C_{10}(\F_7)$ 
in the reduction of $J_{10}(\Q)$ are the reductions of the two 
known rational points. So it suffices to show that there is
only one point of $C_{10}(\Q)$ in each of those residue classes.
The differential mod~7 is $\tilde{\eta} \leftrightarrow x + y + 3z$;
it does not vanish at either of the two relevant points in $C_{10}(\F_7)$,
proving the claim.

\section{Mordell-Weil sieve for $C_5$}
\label{S:Mordell-Weil sieve for C_5}

In this section we prove that $C_5(\Q)_\subs=\emptyset$.
The method used is the Mordell-Weil sieve described in
Section~\ref{S:Mordell-Weil sieve}.

\subsection{The strategy for $C_5$}
\label{S:The strategy for $C_5$}

Because $C_5(\Q)$ is nonempty and $\rank J_5(\Q) = \dim J_5$,
we cannot expect to determine $C_5(\Q)$ 
using only the methods described in Section~\ref{S:Mordell-Weil sieve}.
On the other hand, by Lemma~\ref{L:C5 constraints} we need only
the subset $C_5(\Q)_\subs$ of rational points lying 
in certain residue classes,
so there is hope of proving $C_5(\Q)_\subs = \emptyset$.

To succeed, we must use the $2$-adic and $3$-adic
conditions defining $C_5(\Q)_\subs$, since neither alone is enough
to exclude all the rational points on $C_5$.
This means that we must use the sieve information at the bad primes
$2$ and $3$.
In fact, we will use only the component groups of the N\'eron models
at $2$ and $3$.
We will discover that the only primes dividing the orders of these groups
are $2$ and $7$ (respectively),
so we will sieve also at primes $p$ of good reduction
chosen so that $J_5(\F_p)$ has large $2$-primary and/or $7$-primary parts.

Let $P_0$, $P_1$, $P_2$, $P_3$ (in this order)
be the four rational points on $C_5$ 
listed in Table~\ref{Table:rational points}.
Let $Q_i \in J(\Q)$ be the class of the divisor $P_i-P_0$.
Let $\calQ$ be the subgroup of $J(\Q)$ generated by the $Q_i$.
We will show that $(J_5(\Q):\calQ)$ is prime to $14$,
and we will use groups $\Phi$ of order divisible only by $2$ and $7$
in the sieve.
We will determine the set of $(n_1,n_2,n_3) \in \Z^3$
such that $n_1 Q_1 + n_2 Q_2 + n_3 Q_3$ satisfies the sieve conditions.

\subsection{Sieve information at $3$}

We will compute the minimal proper regular model $\calC$ 
of $C_5$ over $\Z_3$, and then the component group of the N\'eron model
of its Jacobian.
Our exposition will be terse,
since \cite{Flynn-et-al2001} contains a detailed account of the general method.
We will use the same label for a component and for its strict transform(s)
after blowing up.
Subscripts on components denote multiplicities, 
and a component without a subscript has multiplicity~$1$.

Start with the model in $\PP^2_{\Z_3}$ defined by the ternary quartic.
Its special fiber has two irreducible components $A$ and $B$,
and is regular except at the point $s$ where $A$ intersects $B$.

\begin{figure}
\begin{center}
  \Graph{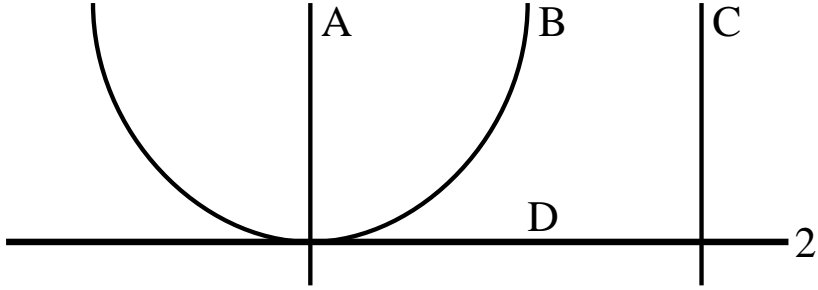}{0.7\textwidth}
\end{center}
\caption{The special fiber of $C_5$ at $3$.}
\label{F:C5 mod 3}
\bigskip
\hrule
\end{figure}

Blowing up $s$ gives a regular $\Z_3$-scheme $\calC$.
Its special fiber, shown in Figure~\ref{F:C5 mod 3},
has four irreducible components, all isomorphic to $\PP^1_{\F_3}$.

The intersection matrix is
\begin{center}
\begin{tabular}{c|cccc}
      & $A$   & $B$   & $C$   & $D_2$ \\ \hline
$A$   & $-3$  & $1$   & $0$   & $1$   \\
$B$   & $1$   & $-5$  & $0$   & $2$   \\
$C$   & $0$   & $0$   & $-2$  & $1$   \\
$D_2$ & $1$   & $2$   & $1$   & $-2$  \\
\end{tabular}
\end{center}
(The diagonal entries are computed using the condition
that the inner product of each row with the multiplicity vector $(1,1,1,2)$
is $0$.)

Let $\calJ$ be the N\'eron model of $J_5$ over $\Z_3$.
Then the component group $\Phi$ of $\calJ_{\Fbar_3}$ is
obtained as follows.
View the multiplicity vector as a $\Z$-linear map $\Z^4 \to \Z$.
Take the quotient of its kernel by the span of the rows of the intersection
matrix.
We get an isomorphism $\Phi \isom \Z/7\Z$.

We will use the homomorphic image $\Phi=\Z/7\Z$ of $\calJ(\F_3)$
for the sieve.
The points $P_0,P_1,P_2,P_3$ map to components $B,C,B,A$ respectively.
If $P \in C_5(\Q)$, then the corresponding point on the special fiber
lies on $A$, $B$, or $C$ (it must be a component of multiplicity~$1$),
and the class of $P-P_0$ in $J(\Q)$ 
maps to the element $3$, $0$, or $1$ in $\Z/7\Z$ accordingly
(for a particular choice of isomorphism $\Phi \isom \Z/7\Z$).
In particular, $Q_1,Q_2,Q_3$ map to $1,0,3$ in $\Z/7\Z$.

If moreover $P \in C_5(\Q)_\subs$, then the image of $P$ 
in $\calC(\F_3)$ is not on $A^\smooth$ or $B^\smooth$ in the regular model
since smooth points on these components in the regular model
map to points in the plane quartic model other than $s=(0:1:0)$.
Therefore the image of $P$ lies on $C$,
and $P-P_0$ maps to $1$ in $\Z/7\Z$.
Thus if $n_1 Q_1 + n_2 Q_2 + n_3 Q_3 = [P-P_0]$, then
\begin{equation}
\label{E:sieve condition at 3}
        n_1 + 3 n_3 \equiv 1 \pmod{7}.
\end{equation}

\subsection{Sieve information at $2$}

We repeat the computations of the previous section, but at the prime $2$.

The fiber at $2$ of the original plane model consists
of a conic $A$ and a line $B_2$.
They are tangent at the point $s_1:=(1:0:0)$.
Because of the multiplicity, no points on $B_2$ are smooth.
But the whole $\Z_2$-scheme is regular
except at the points $s_1$ and $s_2:=(1:1:1)$ on the special fiber.

Blowing up $s_1$ yields a model whose special fiber 
has new components $C$ and $D_2$, 
and a new non-regular point $s_3$ at the 
intersection of $A$, $B_2$, $C$, $D_2$.
Blowing up $s_2$ produces new components $E$ and $F$
and a new non-regular point $s_4$ at the intersection of $B_2$, $E$, $F$.
(For the sake of definiteness,
label $E$ and $F$ so that $P_3 \in C_5(\Q)$ reduces to a point on $F$.)
Blowing up $s_3$ produces new components $G_2$ and $H_4$,
and a new non-regular point $s_5$ at the intersection of $A$, $G_2$, $H_4$.
Blowing up $s_4$ produces a new component $I_2$,
and no new non-regular points.
Blowing up $s_5$ produces new components $J_3$ and $K_4$
and a new non-regular point $s_6$ at the intersection of $A$, $J_3$, $K_4$.
Blowing up $s_6$ produces new components $L_4$ and $M_4$
(let $L_4$ be the one that intersects $A$)
and a new non-regular point $s_7$ at the intersection of $L_4$ and $M_4$.
Blowing up $s_7$ produces new components $N_4$ and $O_4$
(let $N_4$ be the one that intersects $L_4$)
and a new non-regular point $s_8$ at the intersection of $N_4$ and $O_4$.
Blowing up $s_8$ produces a new component $P_4$ 
and no new non-regular points.

\begin{figure}
\begin{center}
  \Graph{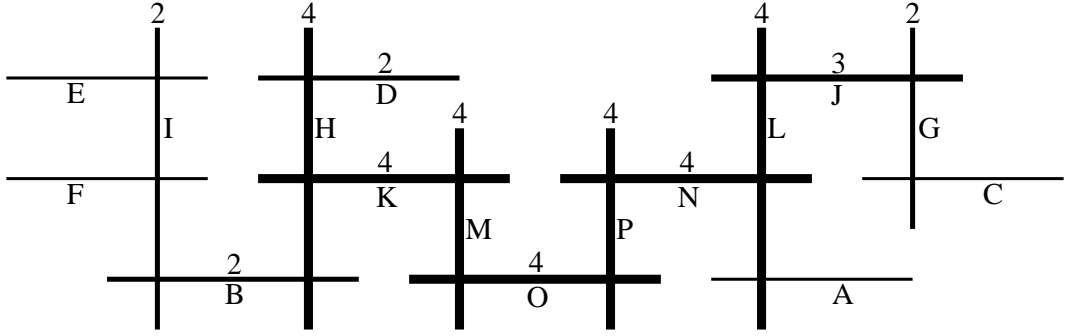}{0.9\textwidth}
\end{center}
\caption{The special fiber of $C_5$ at $2$.}
\label{F:C5 mod 2}
\bigskip
\hrule
\end{figure}

The resulting $\Z_2$-scheme $\calC$ is regular, and its special fiber has
$16$ components, all isomorphic to $\PP^1_{\F_2}$: 
see Figure~\ref{F:C5 mod 2}.
All the intersection numbers are $0$ or $1$,
except the self-intersection numbers,
which are $-4$ for $A$, $-3$ for $B_2$, and $-2$ for each other component.
The points $P_0$, $P_1$, $P_2$, $P_3$ map to components
$C$, $A$, $A$, $F$ respectively.

Let $\calJ$ be the N\'eron model of $J_5$ over $\Z_2$.
We find that the component group $\Phi$ of $\calJ_{\Fbar_2}$ 
is isomorphic to $\Z/4\Z \times \Z/4\Z$.
If $P \in C_5(\Q)$,
then the corresponding point on the special fiber lies on one of
$A$, $C$, $E$, $F$,
and the class of $P-P_0$ in $J(\Q)$ maps to
$(0,3)$, $(0,0)$, $(1,0)$, $(1,2)$ in $\Z/4\Z \times \Z/4\Z$
accordingly
(for a particular choice of isomorphism $\Phi \isom \Z/4\Z \times \Z/4\Z$).
In particular,
$Q_1$, $Q_2$, $Q_3$
map to
$(0,3)$, $(0,3)$, $(1,2)$ respectively.

If moreover $P \in C_5(\Q)_\subs$, then $P$ does not reduce to a point on $A$.
Thus
\begin{equation}
\label{E:sieve condition at 2}
        n_1 (0,3) + n_2 (0,3) + n_3 (1,2) \equiv (0,0) \text{ or } (1,0) \text{ or } (1,2) \quad \pmod{4}.
\end{equation}

\subsection{Sieve information at $23$}

Using Magma's \verb+ClassGroup+ function,
we find that 
\[
J_5(\F_{23}) \isom \frac{\Z}{2\Z} \times \frac{\Z}{2^4\Z} \times \frac{\Z}{2^4\Z} \times \frac{\Z}{2^5\Z}.
\]
Magma also lists the points of $C_5(\F_{23})$ and computes their
images in $J_5(\F_{23})$.
The condition that the image of $n_1 Q_1 + n_2 Q_2 + n_3 Q_3$ 
in $J_5(\F_{23})/4 J_5(\F_{23})$ be in the image of $C_5(\F_{23})$
forces $(n_1,n_2,n_3)$ to be congruent to one of
\begin{center}
$000$, $001$, $002$, $010$, $011$, $021$, $030$, $100$, $101$, $102$, $110$, $113$, $220$, $303$, $312$, $322$
\end{center}
modulo $4$ 
(we have omitted parentheses and commas within each triple, for readability).
Combining this with \eqref{E:sieve condition at 2}, we find that only
\[
        000, 001, 021, 220
\]
are possible.
In particular $n_1$ and $n_2$ are both even.

\subsection{Sieve information at $97$}

We have $J_5(\F_{97}) \isom (\Z/98\Z)^3$.
The condition that the image of $n_1 Q_1 + n_2 Q_2 + n_3 Q_3$ 
in $J_5(\F_{97})$ be in the image of $C_5(\F_{97})$
imposes congruence conditions on $(n_1,n_2,n_3)$ modulo $98$.
Combining these conditions with~\eqref{E:sieve condition at 3}
and the condition that $n_1$ and $n_2$ be even 
proves that $(n_1,n_2,n_3)$ is congruent modulo $14$ to one of the following:
\[
        (2,10,9), \quad (6,2,10), \quad (6,10,10), \quad (8,0,7).
\]

\subsection{Sieve information at $13$}

The group $J_5(\F_{13})$ is cyclic of order $2 \cdot 7 \cdot 157$.
The image of $C_5(\F_{13})$ in $J_5(\F_{13})/14 J_5(\F_{13})$ has size $6$,
and the resulting conditions on $(n_1,n_2,n_3)$ modulo $14$
contradict those in the previous paragraph.

\subsection{End of proof}

Finally we must verify that $(J_5(\Q):\calQ)$ is prime to $14$.
By the results in Table~\ref{Table:2 descent}, $J_5(\Q) \isom \Z^3$.
The image of $\calQ$ under $J_5(\Q) \to J_5(\F_{23})/2 J_5(\F_{23})$
has $\F_2$-dimension $3$,
so $2 \notdiv (J_5(\Q):\calQ)$.
The image of $\calQ$ under 
\[
        J_5(\Q) \to \frac{J_5(\F_{13})}{7 \, J_5(\F_{13})} \times 
                        \frac{J_5(\F_{97})}{7 \, J_5(\F_{97})}
\]
has $\F_7$-dimension $3$,
so $7 \notdiv (J_5(\Q):\calQ)$.
Thus $C_5(\Q)_\subs = \emptyset$.
This completes the proof of Theorem~\ref{goal}.

\section*{Acknowledgements}

M.S.\ thanks
the University of California at Berkeley for inviting him to visit the
Department of Mathematics for two months in 2002; 
many of the ideas in this paper evolved during this period. 
Thanks go also to the Deutsche
Forschungsgemeinschaft for providing the financial means to pay for this
visit through a Heisenberg Fellowship.
We thank Ken Ribet for suggesting some of the references in the proof of
Lemma~\ref{levelloweringlemma}.
B.P.\ and E.S.\ thank Nils Bruin 
and the Pacific Institute for the Mathematical Sciences
for their hospitality in the summer of 2004:
some calculations were doublechecked and some writing was done there.
B.P.\ thanks also the Isaac Newton Institute for a visit in
the summer of 2005.
B.P.\ was partially supported by 
        NSF grants DMS-9801104 and DMS-0301280, and a Packard Fellowship.
E.S.\ was partially supported by the NSA grant MDA904-03-1-0030.
Computations were done using the software Magma~\cite{Magma} 
and PARI/GP~\cite{PARI2}.



\begin{bibdiv}
\begin{biblist}

\bib{Bennett2004preprint}{article}{
  author={Bennett, Michael A.},
  title={The equation $x^{2n}+y^{2n}=z^5$},
  date={2004},
  note={Preprint},
}

\bib{Bennett-Skinner2004}{article}{
  author={Bennett, Michael A.},
  author={Skinner, Chris M.},
  title={Ternary Diophantine equations via Galois representations and modular forms},
  journal={Canad. J. Math.},
  volume={56},
  date={2004},
  number={1},
  pages={23\ndash 54},
  issn={0008-414X},
  review={MR2031121},
}

\bib{Beukers1998}{article}{
  author={Beukers, Frits},
  title={The Diophantine equation $Ax\sp p+By\sp q=Cz\sp r$},
  journal={Duke Math. J.},
  volume={91},
  date={1998},
  number={1},
  pages={61\ndash 88},
  issn={0012-7094},
  review={MR1487980 (99f:11039)},
}

\bib{Breuil2001}{article}{
  author={Breuil, Christophe},
  author={Conrad, Brian},
  author={Diamond, Fred},
  author={Taylor, Richard},
  title={On the modularity of elliptic curves over $\mathbf Q$: wild 3-adic exercises},
  journal={J. Amer. Math. Soc.},
  volume={14},
  date={2001},
  number={4},
  pages={843\ndash 939 (electronic)},
  issn={0894-0347},
  review={MR1839918 (2002d:11058)},
}

\bib{Bruin2004}{article}{
  author={Bruin, Nils},
  title={Visualising $\text {Sha[2]}$ in abelian surfaces},
  journal={Math. Comp.},
  volume={73},
  date={2004},
  number={247},
  pages={1459\ndash 1476 (electronic)},
  issn={0025-5718},
  review={MR2047096},
}

\bib{Bruin2003}{article}{
  author={Bruin, Nils},
  title={Chabauty methods using elliptic curves},
  journal={J. Reine Angew. Math.},
  volume={562},
  date={2003},
  pages={27\ndash 49},
  issn={0075-4102},
  review={MR2011330 (2004j:11051)},
}

\bib{Bruin2000}{article}{
  author={Bruin, Nils},
  title={On powers as sums of two cubes},
  booktitle={Algorithmic number theory (Leiden, 2000)},
  series={Lecture Notes in Comput. Sci.},
  volume={1838},
  pages={169\ndash 184},
  publisher={Springer},
  place={Berlin},
  date={2000},
  review={MR1850605 (2002f:11029)},
}

\bib{Bruin1999}{article}{
  author={Bruin, Nils},
  title={The Diophantine equations $x\sp 2\pm y\sp 4=\pm z\sp 6$ and $x\sp 2+y\sp 8=z\sp 3$},
  journal={Compositio Math.},
  volume={118},
  date={1999},
  number={3},
  pages={305\ndash 321},
  issn={0010-437X},
  review={MR1711307 (2001d:11035)},
}

\bib{Chabauty1941}{article}{
  author={Chabauty, Claude},
  title={Sur les points rationnels des courbes alg\'ebriques de genre sup\'erieur \`a l'unit\'e},
  language={French},
  journal={C. R. Acad. Sci. Paris},
  volume={212},
  date={1941},
  pages={882\ndash 885},
  review={MR0004484 (3,14d)},
}

\bib{Chen2004preprint}{article}{
  author={Chen, Imin},
  title={On the equation $s^2+y^{2p}=\alpha ^3$},
  date={2004-07-01},
  note={Preprint},
}

\bib{Chevalley-Weil1930}{article}{
  author={Chevalley, C.},
  author={Weil, A.},
  title={Un th\'{e}or\`{e}me d'arithm\'{e}tiques sur les courbes alg\'{e}briques},
  journal={Comptes Rendus Hebdomadaires des S\'{e}ances de l'Acad.\ des Sci., Paris},
  volume={195},
  date={1930},
  pages={570\ndash 572},
}

\bib{Coleman1985torsion}{article}{
  author={Coleman, Robert F.},
  title={Torsion points on curves and $p$-adic abelian integrals},
  journal={Ann. of Math. (2)},
  volume={121},
  date={1985},
  number={1},
  pages={111\ndash 168},
  issn={0003-486X},
  review={MR782557 (86j:14014)},
}

\bib{Coleman1985chabauty}{article}{
  author={Coleman, Robert F.},
  title={Effective Chabauty},
  journal={Duke Math. J.},
  volume={52},
  date={1985},
  number={3},
  pages={765\ndash 770},
  issn={0012-7094},
  review={MR808103 (87f:11043)},
}

\bib{Cremona1997}{book}{
  author={Cremona, J. E.},
  title={Algorithms for modular elliptic curves},
  edition={2},
  publisher={Cambridge University Press},
  place={Cambridge},
  date={1997},
  pages={vi+376},
  isbn={0-521-59820-6},
  review={MR1628193 (99e:11068)},
}

\bib{Curtis-Reiner1988}{book}{
  author={Curtis, Charles W.},
  author={Reiner, Irving},
  title={Representation theory of finite groups and associative algebras},
  series={Wiley Classics Library},
  note={Reprint of the 1962 original; A Wiley-Interscience Publication},
  publisher={John Wiley \& Sons Inc.},
  place={New York},
  date={1988},
  pages={xiv+689},
  isbn={0-471-60845-9},
  review={MR1013113 (90g:16001)},
}

\bib{Darmon1997Faltings}{article}{
  author={Darmon, H.},
  title={Faltings plus epsilon, Wiles plus epsilon, and the generalized Fermat equation},
  journal={C. R. Math. Rep. Acad. Sci. Canada},
  volume={19},
  date={1997},
  number={1},
  pages={3\ndash 14},
  issn={0706-1994},
  review={MR1479291 (98h:11034a)},
}

\bib{Darmon-Merel1997}{article}{
  author={Darmon, Henri},
  author={Merel, Lo{\"{\i }}c},
  title={Winding quotients and some variants of Fermat's last theorem},
  journal={J. Reine Angew. Math.},
  volume={490},
  date={1997},
  pages={81\ndash 100},
  issn={0075-4102},
  review={MR1468926 (98h:11076)},
}

\bib{Darmon-Granville1995}{article}{
  author={Darmon, Henri},
  author={Granville, Andrew},
  title={On the equations $z\sp m=F(x,y)$ and $Ax\sp p+By\sp q=Cz\sp r$},
  journal={Bull. London Math. Soc.},
  volume={27},
  date={1995},
  number={6},
  pages={513\ndash 543},
  issn={0024-6093},
  review={MR1348707 (96e:11042)},
}

\bib{Darmon1993CR}{article}{
  author={Darmon, Henri},
  title={The equation $x\sp 4-y\sp 4=z\sp p$},
  journal={C. R. Math. Rep. Acad. Sci. Canada},
  volume={15},
  date={1993},
  number={6},
  pages={286\ndash 290},
  issn={0706-1994},
  review={MR1260076 (94m:11037)},
}

\bib{Edwards2004}{article}{
  author={Edwards, Johnny},
  title={A complete solution to $X\sp 2+Y\sp 3+Z\sp 5=0$},
  journal={J. Reine Angew. Math.},
  volume={571},
  date={2004},
  pages={213\ndash 236},
  issn={0075-4102},
  review={MR2070150},
}

\bib{Eichler1954}{article}{
  author={Eichler, Martin},
  title={Quatern\"are quadratische Formen und die Riemannsche Vermutung f\"ur die Kongruenzzetafunktion},
  language={German},
  journal={Arch. Math.},
  volume={5},
  date={1954},
  pages={355\ndash 366},
  issn={0003-9268},
  review={MR0063406 (16,116d)},
}

\bib{ElkiesANTS4}{article}{
  author={Elkies, Noam D.},
  title={Rational points near curves and small nonzero $\vert x\sp 3-y\sp 2\vert $ via lattice reduction},
  booktitle={Algorithmic number theory (Leiden, 2000)},
  series={Lecture Notes in Comput. Sci.},
  volume={1838},
  pages={33\ndash 63},
  publisher={Springer},
  place={Berlin},
  date={2000},
  review={MR1850598 (2002g:11035)},
}

\bib{Elkies1999}{article}{
  author={Elkies, Noam D.},
  title={The Klein quartic in number theory},
  booktitle={The eightfold way},
  series={Math. Sci. Res. Inst. Publ.},
  volume={35},
  pages={51\ndash 101},
  publisher={Cambridge Univ. Press},
  place={Cambridge},
  date={1999},
  review={MR1722413 (2001a:11103)},
}

\bib{Ellenberg2004Fermat}{article}{
  author={Ellenberg, Jordan S.},
  title={Galois representations attached to $\mathbb Q$-curves and the generalized Fermat equation $A\sp 4+B\sp 2=C\sp p$},
  journal={Amer. J. Math.},
  volume={126},
  date={2004},
  number={4},
  pages={763\ndash 787},
  issn={0002-9327},
  review={MR2075481},
}

\bib{Faltings1983}{article}{
  author={Faltings, G.},
  title={Endlichkeitss\"atze f\"ur abelsche Variet\"aten \"uber Zahlk\"orpern},
  language={German},
  journal={Invent. Math.},
  volume={73},
  date={1983},
  number={3},
  pages={349\ndash 366},
  issn={0020-9910},
  review={MR718935 (85g:11026a)},
  translation={ title={Finiteness theorems for abelian varieties over number fields}, booktitle={Arithmetic geometry (Storrs, Conn., 1984)}, pages={9\ndash 27}, translator = {Edward Shipz}, publisher={Springer}, place={New York}, date={1986}, note={Erratum in: Invent.\ Math.\ {\bf 75} (1984), 381}, },
}

\bib{Flynn2004}{article}{
  author={Flynn, E. V.},
  title={The Hasse principle and the Brauer-Manin obstruction for curves},
  journal={Manuscripta Math.},
  volume={115},
  date={2004},
  number={4},
  pages={437\ndash 466},
  issn={0025-2611},
  review={MR2103661},
}

\bib{Flynn-et-al2001}{article}{
  author={Flynn, E. Victor},
  author={Lepr{\'e}vost, Franck},
  author={Schaefer, Edward F.},
  author={Stein, William A.},
  author={Stoll, Michael},
  author={Wetherell, Joseph L.},
  title={Empirical evidence for the Birch and Swinnerton-Dyer conjectures for modular Jacobians of genus 2 curves},
  journal={Math. Comp.},
  volume={70},
  date={2001},
  number={236},
  pages={1675\ndash 1697 (electronic)},
  issn={0025-5718},
  review={MR1836926 (2002d:11072)},
}

\bib{Ghioca2002}{misc}{
  author={Ghioca, Dragos},
  note={e-mail message},
  date={2002-09-01},
}

\bib{Halberstadt-Kraus2003}{article}{
  author={Halberstadt, Emmanuel},
  author={Kraus, Alain},
  title={Sur la courbe modulaire $X\sb E(7)$},
  language={French, with French summary},
  journal={Experiment. Math.},
  volume={12},
  date={2003},
  number={1},
  pages={27\ndash 40},
  issn={1058-6458},
  review={MR2002672 (2004m:11090)},
}

\bib{Katz-Mazur1985}{book}{
  author={Katz, Nicholas M.},
  author={Mazur, Barry},
  title={Arithmetic moduli of elliptic curves},
  series={Annals of Mathematics Studies},
  volume={108},
  publisher={Princeton University Press},
  place={Princeton, NJ},
  date={1985},
  pages={xiv+514},
  isbn={0-691-08349-5},
  isbn={0-691-08352-5},
  review={MR772569 (86i:11024)},
}

\bib{Klein1879b}{article}{
  label={Kle1879},
  author={Klein, Felix},
  title={\"Uber die Transformationen siebenter Ordnung der elliptischen Funktionen},
  language={German},
  journal={Math. Annalen},
  volume={14},
  date={1879},
  pages={428--471},
  reprint={ author={Klein, Felix}, title={Gesammelte mathematische Abhandlungen}, language={German}, note={Dritter Band: Elliptische Funktionen, insbesondere Modulfunktionen, hyperelliptische und Abelsche Funktionen, Riemannsche Funktionentheorie und automorphe Funktionen, Anhang verschiedene Verzeichnisse; Herausgegeben von R. Fricke, H. Vermeil und E. Bessel-Hagen (von F. Klein mit erg\"anzenden Zus\"atzen versehen) Reprint der Erstauflagen [Verlag von Julius Springer, Berlin, 1923]}, publisher={Springer-Verlag}, place={Berlin}, date={1973}, pages={ii+ix+774+36}, review={MR0389520 (52 \#10351)}, },
  translation={ author={Klein, Felix}, title={On the order-seven transformation of elliptic functions}, booktitle={{\it The eightfold way}}, series={Math. Sci. Res. Inst. Publ.}, volume={35}, pages={287\ndash 331}, publisher={Cambridge Univ. Press}, place={Cambridge}, date={1999}, review={MR1722419 (2001i:14042)}, },
}

\bib{Kraus1998}{article}{
  author={Kraus, Alain},
  title={Sur l'\'equation $a\sp 3+b\sp 3=c\sp p$},
  language={French},
  journal={Experiment. Math.},
  volume={7},
  date={1998},
  number={1},
  pages={1\ndash 13},
  issn={1058-6458},
  review={MR1618290 (99f:11040)},
}

\bib{Magma}{article}{
  author={Bosma, Wieb},
  author={Cannon, John},
  author={Playoust, Catherine},
  title={The Magma algebra system. I. The user language},
  note={Computational algebra and number theory (London, 1993). Magma is available at {\tt http://magma.maths.usyd.edu.au/magma/ }\phantom {m}},
  journal={J. Symbolic Comput.},
  volume={24},
  date={1997},
  number={3-4},
  pages={235\ndash 265},
  issn={0747-7171},
  review={MR1484478},
  label={Magma},
}

\bib{PARI2}{misc}{
  author={The PARI group},
  title={PARI/GP, version 2},
  date={2005},
  address={Bordeaux},
  note={available from {\tt http://pari.math.u-bordeaux.fr}},
  label={PARI2},
}

\bib{Poonen-heuristic2005preprint}{misc}{
  author={Poonen, Bjorn},
  title={Heuristics for the Brauer-Manin obstruction for curves},
  date={2005-07},
  note={Preprint},
}

\bib{Poonen-Schaefer1997}{article}{
  author={Poonen, Bjorn},
  author={Schaefer, Edward F.},
  title={Explicit descent for Jacobians of cyclic covers of the projective line},
  journal={J. Reine Angew. Math.},
  volume={488},
  date={1997},
  pages={141\ndash 188},
  issn={0075-4102},
  review={MR1465369 (98k:11087)},
}

\bib{Ribet1990}{article}{
  author={Ribet, K. A.},
  title={On modular representations of ${\rm Gal}(\overline {\bf Q}/{\bf Q})$ arising from modular forms},
  journal={Invent. Math.},
  volume={100},
  date={1990},
  number={2},
  pages={431\ndash 476},
  issn={0020-9910},
  review={MR1047143 (91g:11066)},
}

\bib{Salmon1879}{book}{
  author={Salmon, G.},
  title={A treatise on the higher plane curves},
  edition={3},
  publisher={Hodges, Foster, and Figgis},
  place={Dublin},
  date={1879},
  label={Sal1879},
}

\bib{Schaefer1998}{article}{
  author={Schaefer, Edward F.},
  title={Computing a Selmer group of a Jacobian using functions on the curve},
  journal={Math. Ann.},
  volume={310},
  date={1998},
  number={3},
  pages={447\ndash 471},
  issn={0025-5831},
  review={MR1612262 (99h:11063)},
}

\bib{Schaefer-Stoll2004}{article}{
  author={Schaefer, Edward F.},
  author={Stoll, Michael},
  title={How to do a $p$-descent on an elliptic curve},
  journal={Trans. Amer. Math. Soc.},
  volume={356},
  date={2004},
  number={3},
  pages={1209\ndash 1231 (electronic)},
  issn={0002-9947},
  review={MR2021618 (2004g:11045)},
}

\bib{Scharaschkin2004preprint}{article}{
  author={Scharaschkin, Victor},
  title={The Brauer-Manin obstruction for curves},
  date={2004-12},
  note={Preprint},
}

\bib{SerreGaloisCohomology}{book}{
  author={Serre, Jean-Pierre},
  title={Galois cohomology},
  series={Springer Monographs in Mathematics},
  edition={Corrected reprint of the 1997 English edition},
  note={Translated from the French by Patrick Ion and revised by the author},
  publisher={Springer-Verlag},
  place={Berlin},
  date={2002},
  pages={x+210},
  isbn={3-540-42192-0},
  review={MR1867431 (2002i:12004)},
}

\bib{Serre1987}{article}{
  author={Serre, Jean-Pierre},
  title={Sur les repr\'esentations modulaires de degr\'e $2$ de ${\rm Gal}(\overline {\bf Q}/{\bf Q})$},
  language={French},
  journal={Duke Math. J.},
  volume={54},
  date={1987},
  number={1},
  pages={179\ndash 230},
  issn={0012-7094},
  review={MR885783 (88g:11022)},
}

\bib{SilvermanAEC}{book}{
  author={Silverman, Joseph H.},
  title={The arithmetic of elliptic curves},
  series={Graduate Texts in Mathematics},
  volume={106},
  publisher={Springer-Verlag},
  place={New York},
  date={1992},
  pages={xii+400},
  isbn={0-387-96203-4},
  review={MR 95m:11054},
  note={Corrected reprint of the 1986 original},
}

\bib{Skorobogatov2001}{book}{
  author={Skorobogatov, Alexei},
  title={Torsors and rational points},
  series={Cambridge Tracts in Mathematics},
  volume={144},
  publisher={Cambridge University Press},
  place={Cambridge},
  date={2001},
  pages={viii+187},
  isbn={0-521-80237-7},
  review={MR1845760 (2002d:14032)},
}

\bib{SteinTables}{misc}{
  author={Stein, William A.},
  title={Modular forms database},
  note={Available at {\tt http://modular.fas.harvard.edu/Tables/index.html}\phantom {m}},
}

\bib{Stoll2004preprint}{misc}{
  author={Stoll, Michael},
  title={Independence of rational points on twists of a given curve},
  date={2004},
  note={Preprint},
}

\bib{Stoll2005preprint}{misc}{
  author={Stoll, Michael},
  title={Finite descent and rational points on curves},
  date={2005},
  note={Preprint},
}

\bib{Taylor-Wiles1995}{article}{
  author={Taylor, Richard},
  author={Wiles, Andrew},
  title={Ring-theoretic properties of certain Hecke algebras},
  journal={Ann. of Math. (2)},
  volume={141},
  date={1995},
  number={3},
  pages={553\ndash 572},
  issn={0003-486X},
  review={MR1333036 (96d:11072)},
}

\bib{Wiles1995}{article}{
  author={Wiles, Andrew},
  title={Modular elliptic curves and Fermat's last theorem},
  journal={Ann. of Math. (2)},
  volume={141},
  date={1995},
  number={3},
  pages={443\ndash 551},
  issn={0003-486X},
  review={MR1333035 (96d:11071)},
}

\bib{Wilson2001}{article}{
  author={Wilson, Robert A.},
  title={The Monster is a Hurwitz group},
  journal={J. Group Theory},
  volume={4},
  date={2001},
  number={4},
  pages={367\ndash 374},
  issn={1433-5883},
  review={MR1859175 (2002f:20022)},
}

\end{biblist}
\end{bibdiv}
\end{document}